\documentclass[a4paper, twoside,10pt]{article}
\usepackage{amsmath}
\usepackage{amssymb}
\usepackage{amsfonts}
\usepackage{graphicx}
\usepackage{subfigure}
\usepackage{eucal}
\usepackage{amsthm}
\usepackage{verbatim}
\usepackage{epsfig}
\usepackage{hyperref}

\newcommand{\brk}[1]{\left( #1 \right)}
\newcommand{\abrk}[1]{\left| #1 \right|}

\def\Rc{\mathbb{R}}
\def\Cc{\mathbb{C}}

\newcommand{\F}{\mathcal{F}}
\newtheorem{thm}{Theorem}
\newtheorem{prop}{Proposition}
\newtheorem{definition}{Definition}

\newtheorem*{goal}{A requirement for PSWFs-based quadrature rules}
\newtheorem{remark}{Remark}
\newtheorem{conjecture}{Conjecture}

\begin{document}

\title{
On the evaluation of prolate spheroidal wave functions and
associated quadrature rules}
\author{Andrei Osipov\footnote{This author's research was supported in part
 by the AFOSR grant \#FA9550-09-1-0241.}
\footnote{Yale University, 51 Prospect st, New Haven, CT 06511.
Email: andrei.osipov@yale.edu.
},
Vladimir Rokhlin\footnote{
This author's research was supported in part 
by the ONR grants \#N00014-10-1-0570, \#N00014-11-1-0718,
the AFOSR grant \#FA9550-09-1-0241,
and the ONRgrant
\#N00014-10-C-0176.}
\footnote{
This author has a significant financial interest
in the Fast Mathematical Algorithms and Hardware corporation
(FMAHc) of Connecticut.}
}

\maketitle

\begin{abstract}
As demonstrated by Slepian et. al. in a sequence of
classical papers (see \cite{SlepianComments},
\cite{ProlateSlepian1},
\cite{ProlateLandau1},
\cite{ProlateSlepian2},
\cite{SlepianAsymptotic}),
prolate spheroidal wave functions (PSWFs)
provide a natural and efficient tool
for computing with bandlimited functions defined on an interval.
Recently, PSWFs have been becoming increasingly popular
in various areas in which such functions occur - this
includes
physics (e.g. wave phenomena, fluid dynamics),
engineering (signal processing, filter design), etc.

To use PSWFs as a computational tool, one needs
fast and accurate numerical algorithms
for the evaluation of PSWFs and related quantities,
as well as for the construction of corresponding quadrature rules,
interpolation formulas, etc. 
During the last 15 years,
substantial progress has been made in
the design of such algorithms - see, for example, \cite{RokhlinXiaoProlate}
(see also \cite{Bouwkamp},
\cite{ProlateSlepian1}, \cite{ProlateLandau1}, \cite{ProlateSlepian2} 
for some classical results).

% ALTERNATIVE - start
% The existing algorithms are quite satisfactory
% for moderate values of band limit $c$ (e.g. $c \leq 10^3$),
% in terms of both speed and accuracy.
% Still, in the modern computational environment
% one might encounter problems with relatively large values of $c$
% (e.g. $c \geq 10^4$).
% ALTERNATIVE - end

% MAIN - start
The complexity of many of the existing algorithms, however,
is at least quadratic in the band limit $c$.
For example, the evaluation of the $n$th
eigenvalue of the prolate integral operator requires at least
$O(c^2)$ operations (see e.g. \cite{RokhlinXiaoProlate}); 
the construction of accurate quadrature
rules for the integration (and associated interpolation)
of bandlimited functions with band limit $c$
requires $O(c^3)$ operations (see e.g. \cite{ChengRokhlin}).
Therefore, while 
the existing algorithms are
satisfactory for moderate values of $c$
(e.g. $c \leq 10^3$), 
they tend to be relatively slow
when $c$ is large
(e.g. $c \geq 10^4$).
% MAIN - end

In this paper, we describe several numerical algorithms
for the evaluation of PSWFs and related quantities,
and design a class of PSWF-based quadratures for the integration of 
bandlimited functions. 
While the analysis is somewhat involved and will be published
separately (currently, it can be found in
\cite{Report3}, \cite{Report3Arxiv}), the resulting
numerical algorithms are quite simple and efficient in practice.
For example, the evaluation of the $n$th
eigenvalue of the prolate integral operator requires
$O(n+c \cdot \log c)$ operations;  
the construction of accurate quadrature
rules for the integration (and associated interpolation)
of bandlimited functions with band limit $c$
requires $O(c)$ operations.
All algorithms described in this paper produce results
essentially to machine precision.
Our results are illustrated via several numerical experiments.
\end{abstract}


\noindent
{\bf Math subject classification:} {
33E10, 34L15, 35S30, 42C10, 45C05, 54P05,
65D05, 65D15, 65D30, 65D32}
% 65D05 interpolation
% 65D15 algorithms for function approximation
% 65D30 numerical integration
% 65D32 quadrature and cubature formulas

%%%%%%%%%%%%%%%%%%%%%%%%%%%%%%%%%%%%%%%%%%%
% User defined macros
%%%%%%%%%%%%%%%%%%%%%%%%%%%%%%%%%%%%%%%%%%%

%%%%%%%%%%%%%%%%%%%%%%%%%%%%%%%%%%%%%%%%%%%%%
%%%%%%%%%%%%%%%%%%%%%%%%%%%%%%%%%%%%%%%%%%%%%
\section{Introduction}
\label{sec_intro}

The principal purpose of this paper is to describe several numerical
algorithms associated with bandlimited functions. While these
algorithms are quite simple and efficient in practice,
the analysis is somewhat involved, and will be published separately
(currently the proofs and additional details can be found in
\cite{Report3}, \cite{Report3Arxiv}, \cite{Report4}, \cite{Report4Arxiv}).

A function $f: \Rc \to \Rc$ is said to be 
bandlimited with band limit $c>0$ if there
exists a function $\sigma \in L^2\left[-1,1\right]$ such that
\begin{align}
f(x) = \int_{-1}^1 \sigma(t) e^{icxt} \; dt.
\label{eq_intro_f}
\end{align}
In other words, the Fourier transform of a bandlimited function
is compactly supported.
While \eqref{eq_intro_f} defines $f$ for all real $x$, 
one is often interested in bandlimited functions whose 
argument is confined to an interval, e.g. $-1 \leq x \leq 1$.
Such functions are encountered in physics (wave phenomena,
fluid dynamics), engineering (signal processing), etc.
(see e.g. \cite{SlepianComments}, \cite{Flammer}, \cite{Papoulis}).

About 50 years ago it was observed that the eigenfunctions of
the integral operator $F_c: L^2\left[-1,1\right] \to L^2\left[-1,1\right]$,
defined via the formula
\begin{align}
F_c\left[\varphi\right] \left(x\right) = \int_{-1}^1 \varphi(t) e^{icxt} \; dt,
\label{eq_intro_fc}
\end{align}
provide a natural tool for dealing with bandlimited functions defined
on the interval $\left[-1,1\right]$. Moreover, it
was observed 
(see \cite{ProlateSlepian1}, \cite{ProlateLandau1}, \cite{ProlateSlepian2})
that the eigenfunctions of $F_c$
are precisely the prolate spheroidal wave functions (PSWFs),
well known from the mathematical physics (see, for example,
\cite{PhysicsMorse}, \cite{Flammer}).

Obviously, to use PSWFs as a computational tool, one needs
fast and accurate numerical algorithms
for the evaluation of PSWFs and related quantities,
as well as for the construction of quadratures,
interpolation formulas, etc. 
For the last 15 years, substantial progress has been made in
the design of such algorithms - see, for example,
\cite{RokhlinXiaoProlate}
(see also \cite{Bouwkamp}, \cite{ProlateSlepian1}, 
\cite{ProlateLandau1}, \cite{ProlateSlepian2} for some classical results).

% ALTERNATIVE - start
% The existing algorithms are quite satisfactory
% for moderate values of band limit $c$ (e.g. $c \leq 10^3$),
% in terms of both speed and accuracy.
% Still, in the modern computational environment
% one might encounter problems with relatively large values of $c$
% (e.g. $c \geq 10^4$).
% ALTERNATIVE - end

% MAIN - start
The complexity of many of the existing algorithms, however,
is at least quadratic in the band limit $c$.
For example, the evaluation of the $n$th
eigenvalue of the prolate integral operator requires 
$O(c^2+n^2)$ operations (see e.g. \cite{RokhlinXiaoProlate}); 
also, the construction of accurate quadrature
rules for the integration (and associated interpolation)
of bandlimited functions with band limit $c$
requires $O(c^3)$ operations (see e.g. \cite{ChengRokhlin}).
Therefore, while 
the existing algorithms are
satisfactory for moderate values of $c$
(e.g. $c \leq 10^3$), 
they tend to be relatively slow
when $c$ is large
(e.g. $c \geq 10^4$).
% MAIN - end

In this paper, we describe several numerical algorithms
for the evaluation of PSWFs and related quantities,
and design a class of PSWF-based quadratures for the integration of 
bandlimited functions. 
While the analysis is somewhat involved and will be published
separately (currently, it can be found in
\cite{Report3}, \cite{Report3Arxiv}), the resulting
numerical algorithms are quite simple and efficient in practice.
For example, the evaluation of the $n$th
eigenvalue of the prolate integral operator requires
$O(n+c \log c)$ operations; 
also,  
the construction of accurate quadrature
rules for the integration of bandlimited functions with band limit $c$
requires $O(c)$ operations. In addition, the evaluation of the $n$th
PSWF is done in two steps. First, we carry out 
a certain precomputation, that requires
$O(n+c \log c)$ operations. Then, each subsequent evaluation of
this PSWF at a point in $[-1,1]$ requires $O(1)$ operations.

This paper is organized as follows.
Section~\ref{sec_overview} contains a brief overview. 
Section~\ref{sec_prel} contains mathematical and numerical
preliminaries to be used in the rest of the paper.
Section~\ref{sec_quad} contains the summary
of the principal analytical results of the paper.
Section~\ref{sec_num_algo} contains the description and analysis
of the numerical algorithms for the evaluation of the quadrature rules
and some related quantities.
In Section~\ref{sec_num_res}, we report some numerical results.
In Section~\ref{sec_num_ill}, we illustrate the analysis
via several numerical experiments.

\section{Overview}
\label{sec_overview}
In this section, we provide an overview of 
the paper. More specifically, Section~\ref{sec_outline_pswf}
is dedicated to the numerical evaluation of PSWFs and related quantities.
In Section~\ref{sec_outline_quad}, we discuss several existing
quadrature rules for the integration of bandlimited functions. 
In Section~\ref{sec_intuition}, we introduce a new class of
PSWFs-based quadrature rules and describe the underlying ideas.
In Section~\ref{sec_overview_analysis}, we outline the
analysis (further details can
be found in \cite{Report3}, \cite{Report3Arxiv}).

\subsection{Numerical Evaluation of PSWFs}
\label{sec_outline_pswf}

For any real $c>0$ and integer $n \geq 0$,
the corresponding PSWF $\psi_n$ can be expanded into
an infinite series of Legendre polynomials (see Section~\ref{sec_legendre}).
The coefficients of such expansions decay superalgebraically
(see e.g \cite{RokhlinXiaoProlate}); 
in particular,
relatively few terms of the Legendre series
are required to evaluate $\psi_n(x)$ to essentially the machine
precision, for any $-1 \leq x \leq 1$.
The use of this observation for the numerical evaluation of PSWFs
goes back at least to the classical Bouwkamp algorithm
\cite{Bouwkamp}
(see also Section~\ref{sec_legendre},
in particular Theorem~\ref{thm_tridiagonal} and Remark~\ref{rem_tridiagonal},
and \cite{RokhlinXiaoProlate} for more details).

Thus, the evaluation of PSWFs reduces to the evaluation of the
corresponding Legendre coefficients. 
For any integer $n \geq 0$, 
the Legendre coefficients of {\it all} the first $n$ PSWFs
$\psi_0, \psi_1, \dots, \psi_{n-1}$ can be obtained
via the solution of a certain symmetric tridiagonal eigenproblem
roughly of order $\max\left\{n,c\right\}$
(see Theorem~\ref{thm_tridiagonal} and Remark~\ref{rem_tridiagonal}
in Section~\ref{sec_legendre}, and also
\cite{RokhlinXiaoProlate} for more details about this algorithm). 
The corresponding eigenvalues $\chi_0, \chi_1, \dots, \chi_{n-1}$ of
the prolate differential operator (see Theorem~\ref{thm_prolate_ode}
in Section~\ref{sec_pswf}) are obtained as a by-product of this procedure.
On the other hand, 
additional computations are required to evaluate
the corresponding
eigenvalues $\lambda_0, \lambda_1, \dots, \lambda_{n-1}$
of the integral operator $F_c$ 
(see \eqref{eq_intro_fc} in Section~\ref{sec_intro}).
In practice, 
it is sometimes desirable to evaluate extremely small $\lambda_j$'s
(e.g. $\mbox{\text{\rm{1E-50}}}$), which presents a numerical challenge
(see Section~\ref{sec_pswf}).
To overcome this obstacle, the algorithm of 
\cite{RokhlinXiaoProlate} evaluates 
$\lambda_0, \lambda_1, \dots, \lambda_{n-1}$
by computing the ratios $\lambda_j/\lambda_{j+1}$, which
turns out to be a well-conditioned numerical procedure
(see \cite{RokhlinXiaoProlate} for more details).

% Roughly $2c/\pi$ of the first
% $\lambda_j$'s have essentially the same magnitudes, while the 
% subsequent ones decay superalgebraically with $j$ (see
% Section~\ref{sec_pswf}); in particular, some
% of the eigenvalues of interest can be extremely small 
% (e.g. $\mbox{\text{\rm{1E-100}}}$) in practice.

Suppose, on the other hand, that one is interested
in a single PSWF $\psi_n$ only (as opposed to all the first $n$ PSWFs).
Obviously, one can use the algorithm of \cite{RokhlinXiaoProlate}; however,
its cost is at least $O(n^2)$ operations (see Remark~\ref{rem_tridiagonal}).
Moreover, the cost of evaluating the corresponding eigenvalue $\lambda_n$
of the prolate integral operator $F_c$ (see \eqref{eq_intro_fc})
via the algorithm of \cite{RokhlinXiaoProlate} is
at least $O(n^2)$ operations, with a large proportionality constant.

In this paper, we describe more efficient algorithms for the numerical
evaluation of $\psi_n$ and associated quantities.
In particular, the cost of the evaluation of the Legendre coefficients
of $\psi_n$ via this algorithm is $O(n+c \log c)$ operations
(see Section~\ref{sec_evaluate_beta}). In addition, the cost
of the evaluation of the eigenvalue $\lambda_n$ is also
$O(n+c \log c)$ operations (see Section~\ref{sec_evaluate_lambda}).
On the other hand, this algorithm has the same accuracy
as that of \cite{RokhlinXiaoProlate}; in other words,
all of the quantities are evaluated to essentially the machine precision
(see Section~\ref{sec_num_algo} for more details). 
Since $\lambda_n$ can be extremely small,
the fact that it can be evaluated to high {\it relative} accuracy
(without computing the preceding $\lambda_j$'s) is, perhaps,
surprising (the related analysis is somewhat subtle, and will be
published separately; see 
\cite{Report4}, \cite{Report4Arxiv} for some preliminary results).

%%%%%%%%%%%%%%%%%%%%%%%%%%%%%%%%%%%%%%%%%%%%%%%%%%%%%%%
\subsection{Quadrature Rules for Bandlimited Functions}
\label{sec_outline_quad}
One of principal goals of this paper is a class of quadrature rules
designed for the integration of bandlimited functions
with a specified band limit $c > 0$ over the interval $[-1,1]$.
Suppose that $n > 0$ is an integer;
a quadrature rule of order $n$ is a pair
$ \left(t_1^{(n)}, \dots, t_n^{(n)}, W_1^{(n)}, \dots 
W_n^{(n)} \right)$ of finite sequences of length $n$,
where
\begin{align}
\label{eq_quadrature_nodes}
-1 < t_1^{(n)} < \dots < t_n^{(n)} < 1
\end{align}
are referred to as "the quadrature nodes", and
\begin{align}
\label {eq_quadrature_weights}
W_1^{(n)}, \dots, W_n^{(n)}
\end{align}
are referred to as "the quadrature weights".
If $f:[-1,1] \to \Rc$ is a bandlimited function
(see \eqref{eq_intro_f} in Section~\ref{sec_intro}),
we use the quadrature rule to
approximate the integral of $f$ over the interval $[-1,1]$ by
a finite sum; more specifically,
\begin{align}
\int_{-1}^1 f(t) \; dt \approx \sum_{j=1}^n W_j^{(n)} f\left(t_j^{(n)}\right).
\label{eq_quadrature_used}
\end{align}
The PSWFs constitute a natural basis for the bandlimited
functions with band limit $c>0$ over the interval $[-1,1]$
(see Section~\ref{sec_intro} above).
Therefore, when designing a quadrature rule for the integration
of such functions,
it is reasonable to require that this quadrature rule integrate several
first PSWFs with band limit $c$ to high accuracy.
To describe this property
in a more precise manner, we introduce the following definition.
%%%%%%%%%%%%%%%
\begin{definition}
Suppose that $c>0$ is a real number, and that $n>0$ is an integer.
Suppose also that a quadrature rule 
for the integration of bandlimited functions with band limit $c$
over $[-1,1]$
is specified via its
$n$ nodes and weights, 
as in \eqref{eq_quadrature_nodes}, \eqref{eq_quadrature_weights}.
Suppose furthermore that $\varepsilon > 0$ is a real number, and that
this quadrature rule integrates the first $n$ PSWFs of band limit $c$
to precision $\varepsilon$, in other words,
\begin{align}
\left| \int_{-1}^1 \psi_m(t) \; dt - 
       \sum_{j=1}^n W_j^{(n)} \psi_m\left(t_j^{(n)}\right) 
\right| \leq \varepsilon,
\label{eq_main_goal}
\end{align}
for every integer $m=0,1,\dots,n-1$, where $\psi_m:[-1,1]\to\Rc$
is the $m$th PSWF corresponding to band limit $c$. We refer to such
quadrature rules as
"quadrature rules of order $n$ to precision $\varepsilon$
(corresponding to band limit $c$)". We omit the reference to $c$
whenever the band limit is clear from the context.
\label{def_quad_n_eps}
\end{definition}
%%%%%%%%%%%%%%%
\begin{remark}
Obviously, if $\varepsilon$ is the machine precision
(e.g. $\varepsilon \approx \mbox{\text{\rm{1D-16}}}$ in
double precision calculations), then quadrature rules
of order $n$ to precision $\varepsilon$ 
(in the sense of Definition~\ref{def_quad_n_eps})
integrate the first $n$ PSWFs exactly,
for all practical purposes.
\label{rem_quad_n_eps}
\end{remark}
%%%%%%%%%%%%%%%

%%%%%%%%%%%%%%%%%%%%%%%%%%%%%%%%%%%%%%%%%%%%%%%%%%%%%%%%%%%%%%
\begin{comment}
With Definition~\ref{def_quad_n_eps} in hand, we formulate one
of principal objectives of this paper as follows.
Suppose that $c>0$ is a real number and $n>0$ is an integer.
We construct a quadrature rule of order $n$,
designed for the integration of bandlimited
functions of band limit $c$. Furthermore, we investigate
the error of these quadrature rules, in the sense of
Definition~\ref{def_quad_n_eps}.
for which values of $n$ these are quadrature rules
of precision $\varepsilon$.

%%%%%%%%%%%%%%%
\begin{goal}
Suppose that $c>0$ is a real number.
For every integer $n > 0$, 
we define a quadrature rule of order $n$ 
(for the integration of bandlimited functions
with band limit $c$ over $[-1,1]$) by specifying
$n$ nodes and 
$n$ weights (see \eqref{eq_quadrature_nodes}, \eqref{eq_quadrature_weights}).
Suppose also that $\varepsilon>0$.
We require that, for sufficiently large $n$,
the quadrature rule of order $n$ integrate
the first $n$ PSWFs of band limit $c$ up to the error $\varepsilon$.
More specifically, we find the integer $M=M(c,\varepsilon)$
such that, for every integer $n \geq M$
and all integer $m=0, 1, \dots, n-1$,
\begin{align}
\abrk{ \int_{-1}^1 \psi_m(t) \; dt - 
       \sum_{j=1}^n W_j^{(n)} \psi_m\brk{t_j^{(n)}} } \leq \varepsilon,
\label{eq_main_goal2}
\end{align}
where $\psi_m : [-1,1] \to \Rc$ is the $m$th PSWF
of band limit $c$
(see Section~\ref{sec_pswf}).
\end{goal}
%%%%%%%%%%%%%
\end{comment}
%%%%%%%%%%%%%%%%%%%%%%%%%%%%%%%%%%%%%%%%%%%%%%%%%%%%%%%%%%%%%%%%%%%

\begin{remark}
In practice, for a quadrature rule of order $n$ to precision $\varepsilon$
to be of any use for the integration of bandlimited functions with
band limit $c$, not only $\varepsilon$ should be "small", but
also $n$ has to be at least equal to $2c/\pi$. See Section~\ref{sec_pswf}
and \cite{RokhlinXiaoProlate} for more details.
\label{rem_useful}
\end{remark}

Quadrature rules for the integration of
bandlimited functions 
% which satisfy \eqref{eq_main_goal}
have already been discussed in the literature, for example:

{\it \bf{Generalized Gaussian Quadrature Rules.}}
Suppose that $n>0$ is an integer, and that $f_1, f_2, \dots, f_{2n}$
are $2n$ linearly independent functions defined on an interval. 
Under very mild conditions
on $f_1, \dots, f_{2n}$, there exists a quadrature rule of order $n$
that integrates these $2n$ functions exactly; moreover,
its weights are usually positive. Such quadrature
rules are referred to as "generalized Gaussian quadrature rules",
and their existence
was first observed more than 100 years ago
(see, for example, \cite{Karlin}, \cite{Krein}, \cite{Markov1},
\cite{Markov2}). Perhaps surprisingly, numerical algorithms for the design
of generalized Gaussian quadrature rules
were constructed only recently (see, for example, \cite{ChengRokhlin},
\cite{MaRokhlin},
\cite{YarvinRokhlin}). These algorithms tend to be
rather expensive (they require $O(n^3)$ operations with a large
proportionality constant).
Thus, the evaluation of the nodes and weights of 
a PSWF-based generalized Gaussian quadrature
rule for accurate integration of bandlimited functions
with band limit $c$ requires $O(c^3)$ operations
(see  Remark~\ref{rem_useful} above,
and also
\cite{RokhlinXiaoProlate} for more details).
\begin{remark}
We observe that a PSWF-based generalized Gaussian quadrature rule
of order $n$ integrates the first $2n$ PSWFs exactly; in other words,
\eqref{eq_main_goal} holds for every integer $m$ 
between $0$ and $2n-1$ with $\varepsilon=0$.
\label{rem_gauss}
\end{remark}

{\it \bf{Quadrature Rules from \cite{RokhlinXiaoProlate}.}}
Suppose now that $n>0$ is an integer, and that $\psi_n$ is the $n$th
PSWF corresponding to band limit $c$. Suppose also that
$t_1, \dots, t_n$ are the roots of $\psi_n$ in the interval
$(-1,1)$
(see Theorem~\ref{thm_pswf_main}
in Section~\ref{sec_pswf}).
Suppose furthermore that $W_1, \dots, W_n$ are real numbers, 
and that
\begin{align}
\sum_{i=1}^n \psi_m(t_i) \cdot W_i = \int_{-1}^1 \psi_m(t) \; dt,
\label{eq_weights_lin_system}
\end{align}
for every $m=0,\dots,n-1$.
% \begin{align}
% \psi_0(t_1) \cdot W_1 + \psi_0(t_2) \cdot W_2 + \dots + 
% \psi_0(t_n) \cdot W_n & = \int_{-1}^1 \psi_0(t) \; dt, \nonumber \\
% %\psi_1(t_1) \cdot W_1 + \psi_1(t_2) \cdot W_2 + \dots + 
% %\psi_1(t_n) \cdot W_n & = \int_{-1}^1 \psi_1(t) \; dt, \nonumber \\
% &\dots 
% \label{eq_weights_lin_system}
% \\
% \psi_{n-1}(t_1) \cdot W_1 + \psi_{n-1}(t_2) \cdot W_2 + \dots + 
% \psi_{n-1}(t_n) \cdot W_n & = \int_{-1}^1 \psi_{n-1}(t) \; dt. \nonumber
% %
% %\left\{ 
% %\sum_{j=1}^n \psi_m(t_j) W_j = \int_{-1}^1 \psi_m(t) \; dt.
% %\right\}_{m=0}^{n-1}
% \end{align}
Obviously, due to \eqref{eq_weights_lin_system}, 
the quadrature rule with nodes $t_1, \dots, t_n$
and weights $W_1, \dots, W_n$ integrates the first $n$ PSWFs
exactly (i.e. \eqref{eq_main_goal} holds for every $m=0,\dots,n-1$
with $\varepsilon=0$).
While this quadrature rule is clearly "sub-optimal" compared to
the generalized Gaussian quadrature rule of order $n$ (the latter
integrates the first $2n$ PSWFs exactly), it is somewhat less expensive
to evaluate. More specifically, the cost of evaluating
the roots $t_1,\dots,t_n$ of $\psi_n$ in $(-1,1)$ and
the weights $W_1,\dots,W_n$, defined via \eqref{eq_weights_lin_system},
is dominated by the cost of solving the dense $n$ by $n$
linear system \eqref{eq_weights_lin_system} for the unknowns
$W_1,\dots,W_n$
(see \cite{RokhlinXiaoProlate} for more details about the numerical
aspects of this procedure).
Thus, due to Remark~\ref{rem_useful} above,
the cost of evaluating the nodes and weights of this quadrature rule
for accurate integration of bandlimited functions
with band limit $c$ requires $O(c^3)$ operations.
\begin{remark}
The cost of the evaluation of the quadrature rule,
defined via \eqref{eq_weights_lin_system},
is $O(c^3)$ operations. The cost of the evaluation
of the generalized Gaussian quadrature rule is also
$O(c^3)$ operations,
but tends to have a larger proportionality constant.
\label{rem_same_cost}
\end{remark}
\begin{remark}
The quadrature rule defined via \eqref{eq_weights_lin_system}
is based on the PSWFs corresponding to band limit $c$.
It turns out, however, that this quadrature rule will also
integrate bandlimited functions
with band limit $2c$ to high accuracy. The reason for this is
that the classical Euclid algorithm for polynomial division
can be generalized 
to the PSWFs; the reader is referred
to \cite{RokhlinXiaoProlate} for further details.
\label{rem_euclid}
\end{remark}

In this paper, we describe another class of quadrature rules whose nodes
are the $n$ roots of $\psi_n$ in $(-1, 1)$. However, their weights 
differ slightly from those defined via \eqref{eq_weights_lin_system}.
In particular, strictly speaking, these quadrature rules do not integrate
the first $n$ PSWFs exactly, as opposed to the generalized
Gaussian quadrature rules and those defined 
via \eqref{eq_weights_lin_system} above. 
Nevertheless, for any $\varepsilon > 0$,
they {\it do} integrate the first $n$ PSWFs to precision $\varepsilon$,
provided that
\begin{align}
n > \frac{2c}{\pi} + 10 + \frac{2}{\pi^2} \cdot \left( \log c \right) \cdot 
                     \log \frac{1}{\varepsilon}
\label{eq_intro_er}
\end{align} 
(see Theorem~\ref{thm_quad_eps_large} from Section~\ref{sec_main_result}
and Conjectures~\ref{conj_quad_error}, \ref{conj_n1} 
from Section~\ref{sec_num_ill}
for more precise statements,
and Experiment 3 in Section~\ref{sec_exp14} for some numerical results).

Thus, provided that $\varepsilon$ is the machine precision and
that \eqref{eq_intro_er} holds, 
the quadrature rules of this paper
are, for all practical purposes, 
as accurate as those defined via \eqref{eq_weights_lin_system} above.
Also, their nodes and weights can be used as starting points
for an iterative scheme that computes the generalized Gaussian 
quadrature rule (see, for example,
\cite{ChengRokhlin},
\cite{MaRokhlin},
\cite{YarvinRokhlin}
for more details). Last but not least, the quadrature rules of this paper are
much faster to evaluate than those described above:
$O(c)$ operations are required
(see Sections~\ref{sec_evaluate_nodes}, \ref{sec_evaluate_weights}).

%%%%%%%%%%%%%%%%%%%%%%%%%%%%%%%%%%%%%%%%%%%%%%%%
\subsection{Intuition Behind Quadrature Weights}
\label{sec_intuition}

In this section, we describe the
quadrature rules of this paper,
and discuss the intuition behind them.

We start with a classical interpolation problem.
Suppose that $t_1, \dots, t_n$ are $n$ distinct points
on the interval $(-1,1)$. We need to find the real numbers
$W_1, \dots, W_n$ such that 
\begin{align}
\int_{-1}^1 p(t) \; dt = \sum_{i=1}^n W_i \cdot p(t_i),
\label{eq_quad_polynomial}
\end{align}
for all polynomials $p$ of degree at most $n-1$.
In other words, the quadrature rule with nodes $t_1,\dots,t_n$
and weights $W_1,\dots,W_n$ integrates all polynomials
of degree up to $n-1$ exactly
(see \eqref{eq_quadrature_nodes}, \eqref{eq_quadrature_weights},
\eqref{eq_quadrature_used}).

To this end, one constructs $n$ polynomials
$l_1, \dots, l_n$ of degree $n-1$ with the property
\begin{align}
\l_j(t_i) = \begin{cases}
0 & i \neq j, \\
1 & i = j
\end{cases}
\label{eq_pol_quad_l_prop}
\end{align}
for every integer $i, j = 1, \dots, n$ (see, for example, \cite{Isaacson}). 
It is easy to verify that, for every $j=1,\dots,n$, the polynomial $l_j$
is defined via the formula
\begin{align}
l_j(t) = \frac{ w_n(t) }{ w_n'(t_j) \cdot (t - t_j) },
\label{eq_pol_quad_l_formula}
\end{align}
for all real $-1 \leq t \leq 1$,
where $w_n$ is defined via the formula
\begin{align}
w_n(t) = (t-t_1) \cdot (t-t_2) \cdot \dots \cdot (t-t_n),
\label{eq_pol_quad_wn}
\end{align}
for all real $-1 \leq t \leq 1$ (in other words,
$w_n$ is
the polynomial of degree $n$ whose roots are
precisely $t_1, \dots, t_n$). The weights $W_1, \dots, W_n$
are defined via the formula
\begin{align}
W_j = \int_{-1}^1 l_j(t) \; dt
    = \frac{1}{w_n'(t_j)} \int_{-1}^1 \frac{w_n(t) \; dt}{t - t_j}, 
\label{eq_pol_quad_weights}
\end{align}
for every integer $j=1, \dots, n$.
% We observe that 
% a single function $w_n$ (see \eqref{eq_pol_quad_wn})
% is used to define all the $n$ weights; also, 
% $w_n$ is a polynomial of degree $n$, and hence does not belong to the space
% of the polynomials of degree up to $n-1$.

In our case, the basis functions are the PSWFs rather than polynomials.
We will consider the quadrature rule 
$\left(t_1,\dots,t_n,W_1,\dots,W_n\right)$,
with
$t_1, \dots, t_n$ the roots of $\psi_n$ on the interval
$(-1,1)$, and $W_1, \dots, W_n$ to be determined.
If we choose the weights $W_1,\dots,W_n$
such that the resulting quadrature rule integrates the 
first $n$ PSWFs exactly, this will lead to
the linear system \eqref{eq_weights_lin_system} 
from
Section~\ref{sec_outline_quad} (and hence 
to the corresponding quadrature rule). Instead, we define the weights
using $\psi_n$ in the same way we used $w_n$ in \eqref{eq_pol_quad_weights}.
More specifically, for every integer $j = 1, \dots, n$,
we define the function $\varphi_j : [-1,1] \to \Rc$ via
the formula
\begin{align}
\varphi_j(t) = \frac{ \psi_n(t) }{ \psi_n'(t_j) \cdot ( t-t_j ) },
\label{eq_quad_phi_first}
\end{align}
with $\psi_n$ the obvious analogue of $w_n$ 
in \eqref{eq_pol_quad_l_formula}.
We observe that, for every integer $i, j = 1, \dots, n$,
\begin{align}
\varphi_j(t_i) = 
\begin{cases}
0 & i \neq j, \\
1 & i = j,
\end{cases}
\label{eq_varphi_ij}
\end{align}
analogous to \eqref{eq_pol_quad_l_prop}. Viewed as a function
on the whole real line, each $\varphi_j$ is bandlimited
with the same band limit $c$ (see, for example,
\cite{Report3}, \cite{Report3Arxiv},
or
Theorem 19.3 in \cite{Rudin}). 
%On the other hand, $\varphi_j$
%does not belong to the span of $\psi_0, \psi_1, \dots, \psi_{n-1}$
%(see \cite{Report3}, \cite{Report3Arxiv} for more details).
We define the weights $W_1, \dots, W_n$ via the formula
\begin{align}
W_j = \int_{-1}^1 \varphi_j(t) \; dt,
\label{eq_quad_w_first}
\end{align}
for every $j = 1, 2, \dots, n$ 
(note the analogy with \eqref{eq_pol_quad_weights}). 
The weights $W_1, \dots, W_n$, defined
via \eqref{eq_quad_w_first}, are different
from the solution of the linear system \eqref{eq_weights_lin_system}.
Nevertheless,
the resulting quadrature rule 
turns out to satisfy \eqref{eq_main_goal}, provided
that $\varepsilon$ is of order $|\lambda_n|$
(see Theorem~\ref{thm_quad_simple} in Section~\ref{sec_quad_error}
for a more precise statement).

The analysis of this issue is somewhat long and involved;
the reader is referred to \cite{Report3}, \cite{Report3Arxiv}
for details and proofs. On the other hand, the underlying ideas
are relatively simple: Section~\ref{sec_overview_analysis} below
contains a short overview of this analysis.

%%%%%%%%%%%%%%%%%%%%%%%%%%%%%%%%%%%%%%%%%%%%%%%%%%%%%%%
\subsection{Overview of the Analysis}
\label{sec_overview_analysis}
The following observation lies at the heart of the analysis:
for any band limit $c>0$ and any integer $n > 0$,
the reciprocal of $\psi_n$ can be approximated by a rational
function with $n$ poles in $(-1,1)$ up to an error
of order $|\lambda_n|$, where $\lambda_n$ is the $n$th
eigenvalue of the integral operator $F_c$ (see \eqref{eq_intro_fc}
in Section~\ref{sec_intro}). 
In other words, the reciprocal of $\psi_n$ resembles
the reciprocal of a polynomial of order $n$, in the following sense.

If $P$ is a polynomial with $n$ simple roots $z_1, \dots, z_n$ in $(-1, 1)$,
then the function $z \to \left(P(z)\right)^{-1}$ 
is meromorphic in the complex plane;
moreover,
\begin{align}
\frac{ 1 }{ P(z) } = \sum_{j = 1}^n \frac{ 1 }{ P'(z_j) \cdot (z - z_j) },
\label{eq_one_over_p}
\end{align}
for all complex $z$ different from $z_1, \dots, z_n$ 
(this is a special case of the well known
Cauchy's integral formula: see, for example, \cite{Rudin}).
% The right-hand side of \eqref{eq_one_over_p}
% is referred to as ``partial fractions expansion of $P^{-1}$''.
Similarly, the function $z \to \left(\psi_n(z)\right)^{-1}$ is meromorphic;
however, it has
infinitely many poles, all of which are real and simple
(see Remark~\ref{rem_continuation} in Section~\ref{sec_pswf}),
and exactly
$n$ of which lie in $(-1, 1)$
(see Theorem~\ref{thm_pswf_main} in Section~\ref{sec_pswf}).
Suppose that the roots of $\psi_n$ in $(-1,1)$
are denoted by $t_1 < \dots < t_n$.
It turns out that 
\begin{align}
\frac{1}{\psi_n(t)} =
\sum_{j=1}^n \frac{ 1 }{ \psi_n'(t_j) \cdot (t-t_j) } + O(|\lambda_n|),
\label{eq_full_expansion}
\end{align}
for all real $-1 \leq t \leq 1$
(note the similarity between \eqref{eq_one_over_p} 
and \eqref{eq_full_expansion}).
In other words, \eqref{eq_full_expansion} means that 
the reciprocal of $\psi_n$
differs from a certain rational function with $n$ poles by a function
whose magnitude in the interval $[-1,1]$ is of order
$|\lambda_n|$.
A rigorous version of \eqref{eq_full_expansion}
is provided by Theorem~\ref{thm_complex} in Section~\ref{sec_pswf}
(its proof is somewhat involved; 
see \cite{Report3}, \cite{Report3Arxiv} for details).
More specifically, according to this theorem,
\begin{align}
\left| \frac{1}{\psi_n(t)} - \sum_{j=1}^n \frac{1}{\psi_n'(t_j)\cdot(t-t_j)}
\right| \leq 
% \text{const} \cdot 
|\lambda_n|
\left( 24 \cdot \log\left( \frac{1}{|\lambda_n|} \right)+ 
      130 \cdot (\chi_n)^{1/4} \right),
\label{eq_pf_final}
\end{align}
for all real $-1 \leq t \leq 1$, where $\chi_n$ is the $n$th
eigenvalue of the prolate differential operator
(see Theorem~\ref{thm_prolate_ode} in Section~\ref{sec_pswf}).

The identity \eqref{eq_full_expansion} is related to the quadrature,
discussed in Section~\ref{sec_intuition} above,
in the following way.
Multiplying both sides of \eqref{eq_full_expansion} by $\psi_n(t)$ and
using \eqref{eq_quad_phi_first},  we obtain
\begin{align}
1 = \varphi_1(t) + \dots + \varphi_n(t) + \psi_n(t) \cdot 
O\left(|\lambda_n|\right)
\label{eq_qw_1}
\end{align}
In other words, $\varphi_1, \dots, \varphi_n$ constitute a partition
of unity on the interval $[-1,1]$, up to an error of order
$|\lambda_n|$.
We integrate both sides of \eqref{eq_qw_1} over $[-1,1]$ 
and use Theorem~\ref{thm_pswf_main} in Section~\ref{sec_pswf} 
and \eqref{eq_quad_w_first} in Section~\ref{sec_intuition}
to obtain
\begin{align}
W_1 + \dots + W_n = 2 + O\left(|\lambda_n|\right),
\label{eq_qw_ws}
\end{align}
where $W_1,\dots,W_n$ are the weights of the quadrature rule
(see Section~\ref{sec_weights} for more details).

Suppose now that $m \neq n$ is an integer. We multiply both sides of
\eqref{eq_qw_1} by $\psi_m$ to obtain
\begin{align}
\psi_m(t) = \left(\sum_{j=1}^n \psi_m(t) \cdot \varphi_j(t)\right) +
\psi_m(t) \cdot \psi_n(t) \cdot
O\left(|\lambda_n|\right).
\label{eq_qw_psim}
\end{align}
A detailed analysis of a combination of
\eqref{eq_pf_final} and \eqref{eq_qw_psim}
leads to the conclusion that,
for all integer $0 \leq m < n$,
\begin{align}
\left|
\int_{-1}^1 \psi_m(t) \; dt - 
\sum_{j=1}^n \psi_m(t_j) \cdot W_j 
\right| \leq
|\lambda_n| \cdot \left(
24 \cdot \log \frac{1}{|\lambda_n|} + 6 \cdot \chi_n
\right)
\label{eq_qw_quad_error}
\end{align}
(see Theorem~\ref{thm_quad_simple} in 
Section~\ref{sec_quad_error}, and also
\cite{Report3}, \cite{Report3Arxiv}
for more details
).

According to \eqref{eq_qw_quad_error}, 
the quadrature rule of order $n$ integrates the first $n$ PSWFs
to precision of order $|\lambda_n|$ (see also \eqref{eq_main_goal}
in Section~\ref{sec_outline_quad}).
 It remains to establish for what
values of $n$ this error is smaller than a predetermined
$\varepsilon > 0$.
Theorem~\ref{thm_quad_eps_simple} 
from Section~\ref{sec_main_result}
provides an answer to this question:
namely, if
\begin{align}
n > \frac{2c}{\pi} +
\left(10 + \frac{3}{2} \cdot \log(c) + 
   \frac{1}{2} \cdot \log\frac{1}{\varepsilon}
\right) \cdot \log\left( \frac{c}{2} \right),
\label{eq_qw_main_n}
\end{align}
then
\begin{align}
\left|
\int_{-1}^1 \psi_m(t) \; dt - 
\sum_{j=1}^n \psi_m(t_j) \cdot W_j 
\right| \leq \varepsilon,
\label{eq_qw_main_error}
\end{align}
for all integer $0 \leq m < n$.

Numerical experiments seem to indicate that 
the situation is even better
in practice: namely, to achieve the accuracy $\varepsilon$
it suffices to pick the minimal $n$ such that
$|\lambda_n| < \varepsilon$, which occurs for
$n \approx 2c/\pi + 2 (\log c) \cdot (-\log \varepsilon) / \pi^2$
(see Section~\ref{sec_num_ill}, in particular,
Conjectures~\ref{conj_quad_error}, \ref{conj_n1}
and Experiment 3 in Section~\ref{sec_exp14}).

%%%%%%%%%%%%%%%%%%%%%%%%%%%%%%%%%%%%%%%%%%%%%
%%%%%%%%%%%%%%%%%%%%%%%%%%%%%%%%%%%%%%%%%%%%%
\section{Mathematical and Numerical Preliminaries}
\label{sec_prel}
In this section, we introduce notation and summarize
several facts to be used in the rest of the paper.

%%%%%%%%%%%%%%%%%%%%%%%%%%%%%%%%%%%%%%%%%%%%%
\subsection{Prolate Spheroidal Wave Functions}
\label{sec_pswf}

In this subsection, we summarize several facts about
the PSWFs. Unless stated otherwise, all these facts can be 
found in \cite{RokhlinXiaoProlate}, 
\cite{RokhlinXiaoApprox},
\cite{LandauWidom},
\cite{ProlateSlepian1},
\cite{ProlateLandau1},
\cite{ReportACHA},
\cite{ReportArxiv}.

Given a real number $c > 0$, we define the operator
$F_c: L^2\left[-1, 1\right] \to L^2\left[-1, 1\right]$ 
via the formula
\begin{align}
F_c\left[\varphi\right] (x) = \int_{-1}^1 \varphi(t) e^{icxt} \; dt.
\label{eq_pswf_fc}
\end{align}
Obviously, $F_c$ is compact. We denote its eigenvalues by
$\lambda_0, \lambda_1, \dots, \lambda_n, \dots$ and assume that
they are ordered such that $|\lambda_n| \geq |\lambda_{n+1}|$
for all natural $n \geq 0$. We denote by $\psi_n$ the eigenfunction
corresponding to $\lambda_n$. In other words,
\begin{align}
\label{eq_prolate_integral}
\lambda_n \psi_n(x) = \int_{-1}^1 \psi_n(t) e^{icxt} \; dt,
\end{align}
for all integer $n \geq 0$ and all real $-1 \leq x \leq 1$.
We adopt the convention\footnote{
This convention agrees with that of \cite{RokhlinXiaoProlate},
\cite{RokhlinXiaoApprox} and differs from that of \cite{ProlateSlepian1}.
}
that $\| \psi_n \|_{L^2\left[-1,1\right]} = 1$.
The following theorem describes the eigenvalues and eigenfunctions
of $F_c$.
%%%%%%%%%%%
\begin{thm}
Suppose that $c>0$ is a real number, and that the operator $F_c$
is defined via \eqref{eq_pswf_fc} above. Then,
the eigenfunctions $\psi_0, \psi_1, \dots$ of $F_c$ are purely real,
are orthonormal and are complete in $L^2\left[-1, 1\right]$.
The even-numbered functions are even, the odd-numbered ones are odd.
Each function $\psi_n$ has exactly $n$ simple roots in $(-1, 1)$.
All eigenvalues $\lambda_n$ of $F_c$ are non-zero and simple;
the even-numbered ones are purely real and the odd-numbered ones
are purely imaginary; in particular, $\lambda_n = i^n |\lambda_n|$,
for every integer $n \geq 0$.
\label{thm_pswf_main}
\end{thm}
%%%%%%%%%
\noindent
We define the self-adjoint operator
$Q_c: L^2\left[-1, 1\right] \to L^2\left[-1, 1\right]$ via the formula
\begin{align}
Q_c\left[\varphi\right] (x) =
\frac{1}{\pi} \int_{-1}^1 
\frac{ \sin \left(c\left(x-t\right)\right) }{x-t} \; \varphi(t) \; dt.
\label{eq_pswf_qc}
\end{align}
Clearly, 
\begin{align}
Q_c\left[\varphi\right] (x) = 
\chi_{\left[-1,1\right]}(x) \cdot
\F^{-1} \left[ 
  \chi_{\left[-c,c\right]}(\xi) \cdot \F\left[\varphi\right](\xi)
\right](x),
\label{eq_pswf_fourier}
\end{align}
where $\F:L^2(\Rc) \to L^2(\Rc)$ is
the Fourier transform, and
$\chi_{\left[-a,a\right]} : \Rc \to \Rc$ is the characteristic
function of the interval $\left[-a,a\right]$, defined via the formula
\begin{align}
\chi_{\left[-a,a\right]}(x) = 
\begin{cases}
1 & -a \leq x \leq a, \\
0 & \text{otherwise},
\end{cases}
\label{eq_char_function}
\end{align}
for all real $x$.
In other words,
$Q_c$ represents low-passing followed by time-limiting.
$Q_c$ relates to $F_c$, defined via \eqref{eq_pswf_fc}, by 
\begin{align}
Q_c = \frac{ c }{ 2 \pi } \cdot F_c^{\ast} \cdot F_c,
\label{eq_pswf_qc_fc}
\end{align}
and the eigenvalues $\mu_n$ of $Q_n$ satisfy the identity
\begin{align}
\mu_n = \frac{c}{2\pi} \cdot |\lambda_n|^2,
\label{eq_prolate_mu}
\end{align}
for all integer $n \geq 0$. Obviously,
\begin{align}
\mu_n < 1,
\label{eq_mu_leg_1}
\end{align}
for all integer $n \geq 0$, due to \eqref{eq_pswf_fourier}.
Moreover, $Q_c$ has the same eigenfunctions $\psi_n$ as $F_c$.
In other words,
\begin{align}
\mu_n \psi_n(x) = \frac{1}{\pi} 
      \int_{-1}^1 \frac{ \sin\left(c\left(x-t\right)\right) }{ x - t } 
  \; \psi_n(t) \; dt,
\label{eq_prolate_integral2}
\end{align}
for all integer $n \geq 0$ and all $-1 \leq x \leq 1$.
Also,  $Q_c$ is closely related to the operator
$P_c: L^2(\Rc) \to L^2(\Rc)$,
defined via the formula
\begin{align}
P_c\left[\varphi\right] (x) =
\frac{1}{\pi} \int_{-\infty}^{\infty}
\frac{ \sin \left(c\left(x-t\right)\right) }{x-t} \; \varphi(t) \; dt,
\label{eq_pswf_pc}
\end{align}
which is a widely known orthogonal projection onto the space
of functions of band limit $c > 0$ on the real
line $\Rc$.

The following theorem can be traced back to \cite{LandauWidom}:
%%%%%%%%%%%
\begin{thm}
Suppose that $c>0$ and $0<\alpha<1$ are positive real numbers,
and that the operator $Q_c: L^2\left[-1,1\right] \to L^2\left[-1,1\right]$
is defined via \eqref{eq_pswf_qc} above.
Suppose also that the integer $N(c,\alpha)$ is the number of 
the eigenvalues $\mu_n$ of $Q_c$ that are greater than $\alpha$. In
other words,
\begin{align}
N(c,\alpha) = \max\left\{ k = 1,2,\dots \; : \; \mu_{k-1} > \alpha\right\}.
\end{align}
Then,
\begin{align}
N(c,\alpha)
= \frac{2c}{\pi} + \left( \frac{1}{\pi^2} \log \frac{1-\alpha}{\alpha} \right)
    \log c + O\left( \log c \right).
\label{eq_mu_spectrum}
\end{align}
\label{thm_mu_spectrum}
\end{thm}
%%%%%%%%%
\noindent
According to \eqref{eq_mu_spectrum}, there are about $2c/\pi$
eigenvalues whose absolute value is close to one, order $\log c$
eigenvalues that decay rapidly, and the rest of them are
very close to zero. 

The eigenfunctions $\psi_n$ of $Q_c$ turn out to be the PSWFs, well
known from classical mathematical physics \cite{PhysicsMorse}.
The following theorem, proved in a more general form in
\cite{ProlateSlepian2},
formalizes this statement.
%%%%%%%%%%%
\begin{thm}
For any $c > 0$, there exists a strictly increasing unbounded sequence
of positive numbers $\chi_0 <  \chi_1 <  \dots$ such that, for
each integer $n \geq 0$, the differential equation
\begin{align}
(1 - x^2) \cdot \psi''(x) - 2 x \cdot \psi'(x) 
+ (\chi_n - c^2 x^2) \cdot \psi(x) = 0
\label{eq_prolate_ode}
\end{align}
has a solution that is continuous on $\left[-1, 1\right]$.
Moreover, all such solutions are constant multiples of 
the eigenfunction $\psi_n$ of $F_c$,
defined via \eqref{eq_pswf_fc} above.
\label{thm_prolate_ode}
\end{thm}
%%%%%%%%%
\begin{remark}
\label{rem_continuation}
For all real $c > 0$ and all integer $n \geq 0$,
\eqref{eq_prolate_integral} defines an analytic continuation
of $\psi_n$ onto the entire complex plane.
All the roots of $\psi_n$ are simple, real, and symmetric
about the origin. Moreover, $\psi_n$ has infinitely many
roots in $\left(1,\infty\right)$.
In addition, the ODE \eqref{eq_prolate_ode} is satisfied
for all complex $x$.
\end{remark}
%%%%%%%%%%%%

%%%%%%%%%%%%%%%%%
Many properties of the PSWF $\psi_n$ depend on
whether the eigenvalue $\chi_n$ of the ODE \eqref{eq_prolate_ode}
is greater than or less than $c^2$. 
In the following theorem from \cite{ReportACHA}, \cite{ReportArxiv}, we describe
a simple relationship 
between $c, n$ and $\chi_n$.
%%%%%%%%%%%%
\begin{thm}
Suppose that $n \geq 2$ is a non-negative integer.
\begin{itemize}
\item If $n \leq (2c/\pi)-1$, then $\chi_n < c^2$.
\item If $n \geq (2c/\pi)$, then $\chi_n > c^2$.
\item If $(2c/\pi)-1 < n < (2c/\pi)$, then either inequality is possible.
\end{itemize}
\label{thm_n_and_khi}
\end{thm}
%%%%%%%%%

In the following theorem, upper and lower bounds on $\chi_n$
in terms of $c$ and $n$ are provided.
\begin{thm}
Suppose that $c > 0$ is a real number,
and $n \geq 0$ is an integer. Then,
\begin{align}
\label{eq_khi_crude}
n \left(n + 1\right) < \chi_n < n \left(n + 1\right) + c^2.
\end{align}
\label{thm_khi_crude}
\end{thm}
%%%%%%%%%%%
It turns out that, for the purposes of this paper,
the inequality \eqref{eq_khi_crude} is insufficiently sharp.
Tighter bounds on $\chi_n$ are described in
the following theorem (see \cite{ReportACHA}, 
\cite{ReportArxiv}).
%%%%%%%%%%%%%%%%%%
\begin{thm}
Suppose that $n \geq 2$ is an integer, and that $\chi_n > c^2$. Then,
\begin{align}
n < & \; 
\frac{2}{\pi} \int_0^1 \sqrt{ \frac{\chi_n - c^2 t^2}{1 - t^2} } \; dt
% = \nonumber \\
% & \; \frac{2}{\pi} \sqrt{\chi_n} \cdot E \brk{ \frac{c}{\sqrt{\chi_n}} }
< n+3.
\label{eq_both_large_simple_prop}
\end{align}
% where the function $E:\left[0,1\right] \to \Rc$ is defined via
% \eqref{eq_E} in Section~\ref{sec_elliptic}.
\label{thm_n_khi_simple}
\end{thm}
%%%%%%%%%%%
%%%%%%%%%%%%%%%%
In the following theorem from \cite{Report2}, \cite{Report2Arxiv}, 
we provide an upper bound on $|\lambda_n|$ in terms of $n$ and $c$.
%%%%%%%%%%%%
\begin{thm}
Suppose that $c>0$ is a real number, and that
\begin{align}
c > 22.
\label{eq_c22}
\end{align}
Suppose also that $\delta>0$ is a real number, and that
\begin{align}
3 < \delta < \frac{\pi c}{16}.
\label{eq_delta_crude}
\end{align}
Suppose, in addition, that $n$ is a positive integer, and that
\begin{align}
n > \frac{2c}{\pi} + \frac{2}{\pi^2} \cdot \delta \cdot
    \log\left( \frac{4e\pi c}{\delta} \right).
\label{eq_n_crude}
\end{align}
Suppose furthermore that the real number $\xi(n,c)$ is defined
via the formula
\begin{align}
\xi(n,c) = 7056 \cdot c \cdot 
\exp\left[-\delta\left(1 - \frac{\delta}{2\pi c}\right) \right].
\label{eq_xi_n_c}
\end{align}
Then, 
\begin{align}
|\lambda_n| < \xi(n,c).
\label{eq_crude_inequality}
\end{align}
\label{thm_crude_inequality}
\end{thm}
%%%%%%%%%%%
%%%%%%%%%%%%%%%%%%%%%
In the following theorem, we provide a recurrence relation
between the derivatives of $\psi_n$ of arbitrary order
(see Lemma 9.1 in \cite{RokhlinXiaoProlate}).
\begin{thm}
Suppose that $c>0$ is a real number, and that $n \geq 0$
is an integer. Then, 
\begin{align}
& \left(1 - t^2\right) \psi_n'''(t) - 4t \psi_n''(t) +
\left(\chi_n - c^2 t^2 - 2\right) \psi_n'(t) - 2c^2 t \psi_n(t) = 0
\label{eq_dpsi_rec3}
\end{align}
for all real $t$. Moreover,
for all integer $k \geq 2$ and all real $t$,
\begin{align}
& \left(1 - t^2\right) \psi_n^{(k+2)}(t) 
  - 2 \left(k+1\right) t \psi_n^{(k+1)}(t)
  + \left(\chi_n - k\left(k+1\right) - c^2 t^2\right) \psi_n^{(k)}(t) \nonumber \\
& \quad \quad 
  -c^2 k t \psi_n^{(k-1)}(t) - c^2 k \left(k-1\right) \psi_n^{(k-2)}(t) = 0.
\label{eq_dpsi_reck}
\end{align}
\label{thm_dpsi_reck}
\end{thm}
%%%%%%%%%%%%%%%%%%%%%%

%
The following theorem asserts that, on the interval $[-1,1]$,
the difference between the reciprocal of $\psi_n$ and a certain
rational function with $n$ poles is of order $|\lambda_n|$.
Its proof can be found in \cite{Report3}, \cite{Report3Arxiv}.

\begin{thm}
Suppose that $c>30$ is a real number, that $n$ is a positive integer, and that
\begin{align}
n > \frac{2c}{\pi} + 7.
\label{eq_complex_1}
\end{align}
Suppose furthermore that $-1 < t_1 < \dots < t_n < 1$ are the roots
of $\psi_n$ in $(-1,1)$, and that the function $\delta:[-1,1] \to \Rc$
is defined via the formula
\begin{align}
\delta(t) = \frac{1}{\psi_n(t)} - 
\sum_{k=1}^n \frac{1}{\psi_n'(t_j) \cdot (t-t_j)},
\label{eq_complex_it}
\end{align}
for all real $-1 \leq t \leq 1$.
Then,
\begin{align}
|\delta(t)| \leq |\lambda_n| \cdot \left(
24 \cdot \log\left( \frac{1}{|\lambda_n|} \right) +
130 \cdot (\chi_n)^{1/4}
\right),
\label{eq_complex}
\end{align}
for all real $-1 \leq t \leq 1$.
\label{thm_complex}
\end{thm}
%%%%%%%%%%%%%%
\begin{remark}
Suppose that the 
function $\delta:[-1,1] \to \Rc$ is defined via \eqref{eq_complex_it}.
If $n$ is even, then $\delta$ is an even function.
If $n$ is odd, then $\delta$ is an odd function.
\label{rem_i_parity}
\end{remark}

%%%%%%%%%%%%%%%%%%%%%%%%%%%%%%%%%%%%%%%%%
\subsection{Legendre Polynomials and PSWFs}
\label{sec_legendre}
In this subsection, we list several well known
facts about Legendre polynomials and the relationship between
Legendre polynomials and PSWFs.
All of these facts can be found, for example, in 
\cite{Ryzhik},
\cite{RokhlinXiaoProlate},
\cite{Abramovitz}.

The Legendre polynomials $P_0, P_1, P_2, \dots$ are defined via
the formulae
\begin{align}
& P_0(t) = 1, \nonumber \\
& P_1(t) = t,
\label{eq_legendre_pol_0_1}
\end{align}
and the recurrence relation
\begin{align}
\left(k+1\right) P_{k+1}(t) = \left(2k+1\right) t P_k(t) - k P_{k-1}(t),
\label{eq_legendre_pol_rec}
\end{align}
for all $k = 1, 2, \dots$. Even Legendre
polynomials are even functions,
and odd Legendre polynomials are odd.
The Legendre polynomials $\left\{P_k\right\}_{k=0}^{\infty}$ constitute a 
complete orthogonal system in $L^2\left[-1, 1\right]$. The normalized
Legendre polynomials are defined via the formula
\begin{align}
\overline{P_k}(t) = P_k(t) \cdot \sqrt{k + 1/2}, 
\label{eq_legendre_normalized}
\end{align}
for all $k=0,1,2,\dots$. The $L^2\left[-1,1\right]$-norm of each
normalized Legendre polynomial equals to one, i.e.
\begin{align}
\int_{-1}^1 \left( \overline{P_k}(t) \right)^2 \; dt = 1.
\label{eq_legendre_normalized_norm}
\end{align}
Therefore, the normalized Legendre polynomials constitute
an orthonormal basis for $L^2\left[-1, 1\right]$. 
In particular, for every real $c>0$
and every integer $n \geq 0$, the prolate spheroidal
wave function $\psi_n$, corresponding
to the band limit $c$, can be expanded into the series
\begin{align}
\psi_n(x) = \sum_{k = 0}^{\infty} \beta_k^{(n)} \cdot \overline{P_k}(x)
          = \sum_{k = 0}^{\infty} \alpha_k^{(n)} \cdot P_k(x),
\label{eq_num_leg_exp}
\end{align}
for all $-1 \leq x \leq 1$,
where $\beta_0^{(n)}, \beta_1^{(n)}, \dots$ are defined via
the formula
\begin{align}
\beta_k^{(n)} = \int_{-1}^1 \psi_n(x) \cdot \overline{P_k}(x) \; dx,
\label{eq_num_leg_beta_knc}
\end{align}
and $\alpha_0^{(n)}, \alpha_1^{(n)}, \dots$ are defined via
the formula
\begin{align}
\alpha_k^{(n)} = \beta_k^{(n)} \cdot \sqrt{k + 1/2}
  = \left(k+1/2\right) \cdot \int_{-1}^1 \psi_n(x) \cdot P_k(x) \; dx,
\label{eq_num_leg_alpha_knc}
\end{align}
for all $k=0, 1, 2, \dots$. Due to the combination of
Theorem~\ref{thm_pswf_main} in Section~\ref{sec_pswf} with
\eqref{eq_legendre_normalized_norm},
\eqref{eq_num_leg_exp}, 
\eqref{eq_num_leg_beta_knc},
\begin{align}
\left( \beta^{(n)}_0 \right)^2 + 
\left( \beta^{(n)}_1 \right)^2 +
\left( \beta^{(n)}_2 \right)^2 + \dots = 1.
\label{eq_beta_unit_length}
\end{align}
For any integer $n \geq 0$, 
the sequence $\beta_0^{(n)}, \beta_1^{(n)}, \dots$ satisfies
the recurrence relation
\begin{align}
A_{0,0} \cdot \beta_0^{(n)} + 
A_{0,2} \cdot \beta_2^{(n)} & = \chi_n \cdot \beta_0^{(n)}, 
\nonumber \\
A_{1,1} \cdot \beta_1^{(n)} + 
A_{1,3} \cdot \beta_3^{(n)} & = \chi_n \cdot \beta_1^{(n)}, 
\nonumber \\
A_{k,k-2} \cdot \beta_{k-2}^{(n)} + 
A_{k,k} \cdot \beta_k^{(n)} +
A_{k,k+2} \cdot \beta_{k+2}^{(n)} & = \chi_n \cdot \beta_k^{(n)}, 
\label{eq_num_a_rec}
\end{align}
for all $k=2,3,\dots$, where $A_{k,k}$, $A_{k+2,k}$, $A_{k,k+2}$ are
defined via the formulae
\begin{align}
& A_{k, k} = k(k+1) + \frac{ 2k(k+1) - 1 }{ (2k+3)(2k-1) } \cdot c^2, 
  \nonumber \\
& A_{k, k+2} = A_{k+2, k} =
\frac{ (k+2)(k+1) }{ (2k+3) \sqrt{(2k+1)(2k+5)} } \cdot c^2,
\label{eq_num_a_matrix}
\end{align}
for all $k=0,1,2,\dots$.
In other words,
the infinite vector $\left( \beta_0^{(n)}, \beta_1^{(n)}, \dots \right)$
satisfies the identity
\begin{align}
\left(A - \chi_n I\right) \cdot 
\left( \beta_0^{(n)}, \beta_1^{(n)}, \dots \right)^T = 0,
\label{eq_num_beta}
\end{align}
where $I$ is the infinite identity matrix,
and the non-zero entries of the infinite symmetric matrix $A$ are given via
\eqref{eq_num_a_matrix}. 
 
The matrix $A$ naturally splits into two infinite 
symmetric tridiagonal matrices,
$A^{even}$ and $A^{odd}$, the former consisting of the elements of $A$
with even-indexed rows and columns, and the latter consisting of
the elements of $A$ with odd-indexed rows and columns.
Moreover, 
for every pair of integers $n, k \geq 0$,
\begin{align}
\beta^{(n)}_k = 0, \quad \text{if } k+n \text{ is odd},
\label{eq_beta_zero_parity}
\end{align}
due to the combination of
Theorem~\ref{thm_pswf_main} in Section~\ref{sec_pswf} and 
\eqref{eq_num_leg_beta_knc}. In the following theorem (that appears
in \cite{RokhlinXiaoProlate} in a slightly different form), we summarize
certain implications of these observations,
that lead to numerical algorithms for the evaluation of PSWFs.
\begin{thm}
Suppose that $c>0$ is a real number, and that the infinite
tridiagonal symmetric matrices $A^{even}$ and $A^{odd}$
are defined, respectively, via
\begin{align}
A^{even} =
\begin{pmatrix}
A_{0,0} & A_{0,2} &  &  &  \\
A_{2,0} & A_{2,2} & A_{2,4} &  &  \\
        & A_{4,2} & A_{4,4} & A_{4,6} & \\
        &         & \ddots  & \ddots & \ddots \\
\end{pmatrix}
\label{eq_a_even}
\end{align}
and
\begin{align}
A^{odd} =
\begin{pmatrix}
A_{1,1} & A_{1,3} &  &  &  \\
A_{3,1} & A_{3,3} & A_{3,5} &  &  \\
        & A_{5,3} & A_{5,5} & A_{5,7} & \\
        &         & \ddots  & \ddots & \ddots \\
\end{pmatrix},
\label{eq_a_odd}
\end{align}
where the entries $A_{k,j}$ are defined via \eqref{eq_num_a_matrix}. Suppose
also that the 
infinite vectors $\beta^{(n)}_{even} \in l^2$ 
and $\beta^{(n)}_{odd} \in l^2$ are 
defined, respectively, via 
the formulae
\begin{align}
\beta^{(n)}_{even} =
\left( \beta^{(n)}_0, \beta^{(n)}_2, \dots \right)^T, \quad
\beta^{(n)}_{odd} = 
\left( \beta^{(n)}_1, \beta^{(n)}_3, \dots \right)^T,
\label{eq_beta_n}
\end{align}
where $\beta^{(n)}_0, \beta^{(n)}_1, \dots$ are defined via
\eqref{eq_num_leg_beta_knc}. If $n$ is even, then
\begin{align}
A^{even} \cdot \beta^{(n)}_{even} = \chi_n \cdot \beta^{(n)}_{even}.
\label{eq_a_even_eig}
\end{align}
If $n$ is odd, then
\begin{align}
A^{odd} \cdot \beta^{(n)}_{odd} = \chi_n \cdot \beta^{(n)}_{odd}.
\label{eq_a_odd_eig}
\end{align}
\label{thm_tridiagonal}
\end{thm}
\begin{remark}
We define the infinite vector $\beta^{(n)} \in l^2$ to be equal
to $\beta^{(n)}_{even}$, if $n$ is even,
or to $\beta^{(n)}_{odd}$, if $n$ is odd. In this notation,
$\beta^{(0)}, \beta^{(2)}, \dots$ are the eigenvectors of $A^{even}$,
and 
$\beta^{(1)}, \beta^{(3)}, \dots$ are the eigenvectors of $A^{odd}$.
\label{rem_n_even_odd}
\end{remark}
\begin{remark}
While the matrices $A^{even}$ and $A^{odd}$ are infinite, and their entries
do not decay with increasing row or column number, the coordinates
of each eigenvector $\beta^{(n)}$ decay superexponentially fast
(see e.g. \cite{RokhlinXiaoProlate} for estimates of this decay).
In particular, suppose that we need to evaluate 
the first $n+1$ eigenvalues $\chi_0, \dots,
\chi_{n}$ and the corresponding eigenvectors 
$\beta^{(0)}, \dots, \beta^{(n)}$ numerically. Then, we can replace
the matrices $A^{even}, A^{odd}$ in \eqref{eq_a_even_eig}, 
\eqref{eq_a_odd_eig}, respectively, with their $N \times N$ upper left
square submatrices, where $N$ is of order $\max\left\{n,c\right\}$,
and solve the resulting symmetric tridiagonal
eigenproblem by any standard technique (see, for example,
\cite{Wilkinson}, \cite{Dahlquist}; see also \cite{RokhlinXiaoProlate}
for more details about this numerical algorithm).
The CPU cost of this procedure is $O(n^2)$ operations.
\label{rem_tridiagonal}
\end{remark}

The Legendre functions of the second kind $Q_0, Q_1, Q_2, \dots$ 
are defined via the formulae
\begin{align}
&  Q_0(t) = \frac{1}{2} \log \frac{1+t}{1-t}, \nonumber \\
& Q_1(t) = \frac{t}{2} \log \frac{1+t}{1-t} - 1,
\label{eq_legendre_fun_0_1}
\end{align}
and the recurrence relation
\begin{align}
\left(k+1\right) Q_{k+1}(t) = \left(2k+1\right) t Q_k(t) - k Q_{k-1}(t),
\label{eq_legendre_fun_rec}
\end{align}
for all $k = 1, 2, \dots$. 
We observe that the recurrence relation
\eqref{eq_legendre_fun_rec} is the same as the recurrence relation
\eqref{eq_legendre_pol_rec}, satisfied by the Legendre polynomials.
In addition, for every integer $k = 0, 1, 2, \dots$, the $k$th Legendre
polynomial $P_k$ and the $k$th Legendre function of the second
kind $Q_k$ are two independent solutions of the second order
differential equation
\begin{align}
(1 - t^2) \cdot y''(t) - 2t \cdot y'(t) + k(k+1) \cdot y(t) = 0.
\label{eq_legendre_ode}
\end{align}
\begin{remark}
Suppose that $-1 \leq x \leq 1$ is a real number, and that
$n \geq 0$ is an integer.
Combining  \eqref{eq_legendre_pol_0_1}, \eqref{eq_legendre_pol_rec}, 
\eqref{eq_legendre_fun_0_1}, \eqref{eq_legendre_fun_rec} gives
a numerical procedure for the evaluation of
$P_0(x), \dots, P_n(x)$ and
$Q_0(x), \dots, Q_n(x)$ to high precision. This procedure
is stable, and requires
$O(n)$ operations (see, for example,
\cite{Dahlquist} for more details).
\label{rem_legendre_evaluate}
\end{remark}

%%%%%%%%%%%%%%%%%%%%%%%%%%%%%%%%%%%%%%%%%%%%%
\begin{comment}
\subsection{Elliptic Integrals}
\label{sec_elliptic}
In this subsection, we summarize several facts about
elliptic integrals. These facts can be found, for example,
in section 8.1 in \cite{Ryzhik}, and in \cite{Abramovitz}.

The incomplete elliptic integrals of the first and second kind
are defined, respectively, by the formulae
\begin{align}
\label{eq_F_y}
& F(y, k) =  \int_0^y \frac{dt}{\sqrt{1 - k^2 \sin^2 t}}, \\
& E(y, k) = \int_0^y \sqrt{1 - k^2 \sin^2 t} \; dt,
\label{eq_E_y}
\end{align}
where $0 \leq y \leq \pi/2$ and $0 \leq k \leq 1$.
By performing the substitution $x = \sin t$, we can write 
\eqref{eq_F_y} and \eqref{eq_E_y} as
\begin{align}
& F(y, k) = \int_0^{\sin(y)}
  \frac{ dx }{ \sqrt{\brk{1 - x^2} \brk{1 - k^2 x^2} } },
\label{eq_F_y_2} \\
\nonumber \\
& E(y, k) = \int_0^{\sin(y)}
\sqrt{ \frac{1 - k^2 x^2}{1 - x^2} } \; dx.
\label{eq_E_y_2}
\end{align}
The complete elliptic integrals of the first and second kind are
defined, respectively, by the formulae
\begin{align}
\label{eq_F}
& F(k) = F\brk{\frac{\pi}{2}, k} = 
\int_0^{\pi/2} \frac{dt}{\sqrt{1 - k^2 \sin^2 t}}, \\
& E(k) = E\brk{\frac{\pi}{2}, k} =
\int_0^{\pi/2} \sqrt{1 - k^2 \sin^2 t} \; dt,
\label{eq_E}
\end{align}
for all $0 \leq k \leq 1$.
\end{comment}

%%%%%%%%%%%%%%%%%%%%%%%%%%%%%%%%%%%%%
\subsection{Pr\"ufer Transformations}
\label{sec_prufer}

The classical Pr\"ufer transformation 
of a second-order ODE 
is a well known analytical
tool for the study of the oscillatory properties of its solutions
(see, for example, \cite{Miller},\cite{Fedoryuk}).
Recently, a minor modification of Pr\"ufer transformation was
demonstrated to be also a convenient numerical tool
(see \cite{Glaser}). In the following theorem,
we summarize several properties of this transformation,
applied to the prolate ODE \eqref{eq_prolate_ode}
(see \cite{Glaser}, \cite{ReportACHA}, \cite{ReportArxiv} for details).

%%%%%%%%%%%
\begin{thm}
Suppose that $n \geq 2$ is an integer, and that
$\chi_n > c^2$.
Suppose also that the functions
$f, v: (-1,1) \to \Rc$ are defined, respectively, via the formulae
\begin{align}
f(t) = \sqrt{ \frac{\chi_n - c^2 t^2}{ 1-t^2 } }
\label{eq_jan_f}
\end{align}
and
\begin{align}
v(t) = 
\frac{1}{2} \left( \frac{t}{1-t^2} + \frac{c^2 t}{\chi_n-c^2 t^2 } \right),
\label{eq_jan_v}
\end{align}
for all real $-1<t<1$.
Suppose furthermore that $t_1$ is the minimal root of $\psi_n$
in $(-1,1)$, and that the function $\theta:(-1,1)\to\Rc$
is the solution of the differential equation
\begin{align}
\theta'(t) = f(t) + v(t) \cdot \sin( 2 \theta(t) )
\label{eq_prufer_theta_ode_old}
\end{align}
with the initial condition
\begin{align}
\theta(t_1) = \frac{\pi}{2}.
\label{eq_theta_at_t1_initial}
\end{align}
Then, $\theta$ has the following properties:
\begin{itemize}
\item $\theta$ extends continuously 
to the interval $\left[-1,1\right]$, and, moreover,
\begin{align}
\label{eq_theta_at_xn_old}
& \theta(-1) = 0, \\
\label{eq_theta_at_0_old}
& \theta(0) = \frac{\pi n}{2}, \\
\label{eq_theta_at_1_old}
& \theta(1) = \pi n.
\end{align}
\item For any real $-1 < t < 1$ such that $\psi_n(t) \neq 0$,
\begin{align}
\theta(t) = 
\text{atan}\left( -\sqrt{\frac{1-t^2}{\chi_n-c^2t^2}} \cdot
                        \frac{\psi_n'(t)}{\psi_n(t)} \right)
          + m(t) \cdot \pi,
\label{eq_prufer_theta_old}
\end{align}
where $m(t)$ is the number of the roots
of $\psi_n$ in the interval $(-1,t)$.
\item For each integer $i = 1, \dots, n$,  
\begin{align}
\label{eq_theta_at_t_old}
\theta(t_i) = \left(i-\frac{1}{2}\right) \cdot \pi,
\end{align}
where
$t_1,\dots,t_n$ are the roots of $\psi_n$ in $(-1,1)$.
\item For all real $-1 < t < 1$,
\begin{align}
\theta'(t) > 0.
\label{eq_dtheta_pos_old}
\end{align}
In other words, $\theta$ is monotonically increasing.
\end{itemize}
\label{thm_prufer_old}
\end{thm}
%%%%%%%%%%%
The following theorem is
closely related to Theorem~\ref{thm_prufer_old}
(see \cite{ReportACHA}, \cite{ReportArxiv} for more details).
% is a direct consequence
% of Theorem~\ref{thm_prufer_old}
% in Section~\ref{sec_prufer}.
\begin{thm}
Suppose that the function 
$\theta: \left[t_1,t_n\right] \to \Rc$ that of
Theorem~\ref{thm_prufer_old}. 
Suppose also that the function
$s: \left[\pi/2, \pi\cdot(n-1/2)\right] \to \left[t_1,t_n\right]$
is the inverse of $\theta$. Then, $s$ is well defined, monotonically
increasing and continuously differentiable. Moreover, for all
real
$\pi/2 < \eta < \pi\cdot(n-1/2)$, 
\begin{align}
s'(\eta) = \frac{1}{f\left(s(\eta)\right) + 
 v\left(s(\eta)\right)\cdot\sin(2\eta)},
\label{eq_ds_ode}
\end{align}
where the functions $f,v$ are defined, respectively,
via \eqref{eq_jan_f}, \eqref{eq_jan_v}.
In addition, for every integer $i=1,\dots,n$,
\begin{align}
s\left( \left(i-\frac{1}{2}\right) \cdot \pi \right) = t_i,
\label{eq_s_at_ti}
\end{align}
and also
\begin{align}
s\left( \frac{\pi n}{2} \right) = 0.
\label{eq_s_at_0}
\end{align}
\label{thm_prufer_inverse}
\end{thm}

%%%%%%%%%%%%%%%%%%%%%%%%%%%%
\subsection{Numerical Tools}
In this subsection, we summarize several numerical techniques
to be used in this paper.

%%%%%%%%%%%%%%%%%%%%%%%%%%%%%%%
\subsubsection{Newton's Method}
\label{sec_newton}
Newton's method solves the equation $f(x) = 0$ iteratively given
an initial approximation $x_0$ to the root $\tilde{x}$.
The $n$th iteration is defined by
\begin{align}
x_n = x_{n-1} - \frac{ f(x_{n-1}) }{ f'(x_{n-1}) }.
\label{eq_newton}
\end{align}
The convergence is quadratic provided that $\tilde{x}$ is a simple
root and $x_0$ is sufficiently close to $\tilde{x}$. More details
can be found e.g. in \cite{Dahlquist}.

%%%%%%%%%%%%%%%%%%%%%%%%%%%%%%%%%%%%%%%%%%%%%%%%%%%%%%%%%%%%%%%%%
\subsubsection{The Taylor Series Method for the Solution of ODEs}
\label{sec_taylor}
The Taylor series method for the solution of a linear second order
differential equation is based on the Taylor formula
\begin{align}
u(x+h) = \sum_{j = 0}^k \frac{ u^{(j)}(x) }{ j! } h^j + O(h^{k+1}).
\label{eq_taylor}
\end{align}
This method evaluates $u(x+h)$ and $u'(x+h)$ by using \eqref{eq_taylor}
and depends on the ability to compute $u^{(j)}(x)$ for $j = 0, \dots, k$.
When the latter satisfy a simple recurrence relation such as 
\eqref{eq_dpsi_reck}
and hence can be computed in $O(k)$ operations, this method is particularly
useful. The reader is referred to \cite{Glaser} for further details.

%%%%%%%%%%%%%%%%%%%%%%%%%%%%%%%%%%%%%%%%%%%%%%%%%
\subsubsection{A Second Order Runge-Kutta Method}
\label{sec_runge_kutta}
A standard second order Runge-Kutta Method
(see, for example, \cite{Dahlquist})
solves
the initial value problem
\begin{align}
y(t_0) = y_0, \quad y'(t) = f(t, y)
\label{eq_rk_ivp}
\end{align}
on the interval $t_0 \leq t \leq t_0 + L$ via the formulae
\begin{align}
& t_{i+1} = t_i + h, \nonumber \\
& k_{i+1} = h f\left(t_{i+1}, y_i + k_i\right), \nonumber \\
& y_{i+1} = y_i + \left(k_i + k_{i+1}\right)/2
\label{eq_runge_kutta}
\end{align}
with $i = 0, \dots, n$, where $h$ and $k_0$ are defined via the formulae
\begin{align}
h = \frac{ L }{ n }, \quad k_0 = f(t_0, y_0).
\label{eq_rk_h_k0}
\end{align}
This procedure requires exactly $n+1$ evaluations of $f$.
The global truncation
error is $O(h^2)$.

%%%%%%%%%%%%%%%%%%%%%%%%%%%%%%%%%%%%%%%%%%%%%%%
\subsubsection{Shifted Inverse Power Method}
\label{sec_power}
%The numerical method described in this subsection is widely known
%(see, for example, \cite{Dahlquist}).
Suppose that $n \geq 0$ is an integer, and that $A$ is an $n$ by $n$
real symmetric matrix. Suppose also that
$\sigma_1 < \sigma_2 < \dots < \sigma_n$ are the eigenvalues of $A$.
The Shifted Inverse Power Method iteratively
finds the eigenvalue $\sigma_k$ and the corresponding
eigenvector $v_k \in \Rc^n$, 
provided that an approximation $\lambda$ to $\sigma_k$
is given, and that
\begin{align}
|\lambda - \sigma_k| < \max \left\{
   |\lambda - \sigma_j| \; : \; j \neq k \right\}.
\label{eq_inverse_power_sigma}
\end{align}
Each Shifted Inverse Power iteration solves the linear system
\begin{align}
\left(A - \lambda_j I\right) \cdot x = w_j
\end{align}
in the unknown $x \in \Rc^n$, where $\lambda_j$ and $w_j \in \Rc^n$
are the approximations to $\sigma_k$ and $v_k$, respectively,
after $j$ iterations; 
the number $\lambda_j$ is usually referred to as "shift".
The approximations $\lambda_{j+1}$ and $w_{j+1} \in \Rc^n$ 
(to $\sigma_k$ and $v_k$, respectively) are evaluated
from $x$ via the formulae
\begin{align}
w_{j+1} = \frac{ x }{ \| x \|}, \quad
\lambda_{j+1} = w_{j+1}^T \cdot A \cdot w_{j+1}
\end{align}
(see, for example, 
\cite{Dahlquist}, \cite{Wilkinson} for more details).
\begin{remark}
\label{rem_power}
For symmetric matrices, the Shifted Inverse Power Method converges cubically
in the vicinity of the solution.
In particular, if
the matrix $A$ is tridiagonal, and the initial approximation $\lambda$
is sufficiently close to $\sigma_k$, the Shifted Inverse Power Method
evaluates $\sigma_k$ and $v_k$ essentially to machine precision
$\varepsilon$
in $O\left( -\log \log \varepsilon \right)$ 
iterations, and each iteration requires
$O(n)$ operations (see e.g
\cite{Wilkinson}, \cite{Dahlquist}).
\end{remark}

%%%%%%%%%%%%%%%%%%%%%%%%%%%%%%
\subsubsection{Sturm Bisection}
\label{sec_sturm}
In this subsection, we describe a well known algorithm for the evaluation
of a single eigenvalue of a real symmetric tridiagonal matrix.
This algorithm is based on the following theorem
that can be found, for example, in \cite{Wilkinson}, \cite{Wilkinson2}.
%%%%%%%%%%%
\begin{thm}[Sturm sequence]
Suppose that $n>0$ is an integer, that
\begin{align}
C = 
\begin{pmatrix} 
a_1 & b_2 & 0 & \cdots & \cdots & 0 \\
b_2 & a_2 & b_3 & 0 & \cdots & 0 \\
\vdots & \ddots & \ddots & \ddots & \ddots & \vdots \\
0 & \cdots & 0 & b_{n-1} & a_{n-1} & b_n \\
0 & \cdots & \cdots & 0 & b_n & a_n
\end{pmatrix}
\label{eq_sturm_c}
\end{align}
is an $n$ by $n$ 
symmetric tridiagonal matrix, and that none of
numbers $b_2, \dots, b_n$ is equal to zero.
Suppose also that
the polynomials $p_{-1}, p_0, \dots, p_n$ are defined via
the formulae
\begin{align}
p_{-1}(x) = 0, \quad p_0(x) = 1
\label{eq_sturm_p}
\end{align}
and
\begin{align}
p_k(x) = \left(a_k - x\right) p_{k-1}(x) - b^2_k p_{k-2}(x), 
\label{eq_sturm_rec}
\end{align}
for all real $x$ and every integer $k = 2, \dots, n$.
Suppose furthermore that $\sigma$ is a real number, and that
the integer $A(\sigma)$ is defined as the number of positive elements
in the finite sequence
\begin{align}
p_0(\sigma) p_1(\sigma), \; p_1(\sigma) p_2(\sigma), 
\; \dots, \; p_{n-1}(\sigma) p_n(\sigma).
\label{eq_sturm_p_seq}
\end{align}
Then,
the number of eigenvalues of $C$ that are strictly larger than $\sigma$
is precisely $A(\sigma)$.
\label{thm_sturm}
\end{thm}
%%%%%%%%%%%
\begin{remark}
Suppose now that $n>0$ is an integer, and $C$ is an $n \times n$ real
symmetric tridiagonal matrix,
such as \eqref{eq_sturm_c}. Theorem~\ref{thm_sturm}
yields a numerical scheme for the evaluation of the $k$th
smallest eigenvalue $\sigma_k$ of $C$. This scheme is known
in the literature as "Sturm Bisection". Provided that
two real numbers $x_0$ and $y_0$ are given such that
\begin{align}
x_0 < \sigma_k < y_0,
\end{align}
% For a predefined $\varepsilon > 0$,
Sturm Bisection requires
\begin{align}
O \left( n \cdot \log_2 
\left( \frac{y_0-x_0}{|\sigma_k|} \right) \right)
\end{align}
operations to evaluate $\sigma_k$ to machine precision
(see, for example, \cite{Wilkinson}, \cite{Wilkinson2} for more details).
\label{rem_sturm}
\end{remark}

\section{Analytical Apparatus}
\label{sec_quad}
The purpose of this section is to provide
the analytical apparatus to be used in the
rest of the paper.
More specifically, we define a PSWF-based quadrature
rule and list several of its properties.

The principal result of this section is Theorem~\ref{thm_quad_eps_simple}.
The reader is referred to \cite{Report3}, \cite{Report3Arxiv}
for the detailed analysis of all the tools listed in this section.

Throughout this section, the band limit $c > 0$ is
assumed to be a positive real number. Also, for any integer $n \geq 0$,
we denote by $\psi_n$ the $n$th PSWF corresponding to the band limit $c$
(see Section~\ref{sec_pswf}).
\begin{definition}
Suppose that $n > 0$ is an integer, and that
\begin{align}
-1 < t_1 < t_2 < \dots < t_n < 1
\label{eq_quad_t}
\end{align}
are the roots of $\psi_n$ in the interval $(-1, 1)$.
For each integer $j = 1, \dots, n$, we define
the function $\varphi_j: [-1,1] \to \Rc$ via the formula
\begin{align}
\varphi_j(t) = \frac{ \psi_n(t) }{ \psi_n'(t_j) \left( t-t_j \right) }.
\label{eq_quad_phi}
\end{align}
In addition, for each integer $j=1,\dots,n$, we define
the real number $W_j$ via the formula
\begin{align}
W_j = \int_{-1}^1 \varphi_j(s) \; ds =
\frac{1}{\psi_n'(t_j)} \int_{-1}^1 \frac{ \psi_n(s) \; ds }{ s - t_j }.
\label{eq_quad_w}
\end{align}
We refer to the pair of finite sequences
\begin{align}
S_n = \left(t_1,\dots,t_n,W_1,\dots,W_n\right)
\label{eq_quad_sn}
\end{align}
as the "PSWF-based quadrature rule of order $n$".
The points $t_1,\dots,t_n$ are referred to as the quadrature nodes, and 
the numbers $W_1,\dots,W_n$ are referred to as
the quadrature weights
(see \eqref{eq_quadrature_nodes}, \eqref{eq_quadrature_weights}
in Section~\ref{sec_outline_quad}).
We use $S_n$
to approximate the integral of a bandlimited function $f$ over
the interval $\left[-1,1\right]$ by a finite sum;
more specifically,
\begin{align}
\int_{-1}^1 f(t) \; dt \approx \sum_{j=1}^n W_j \cdot f(t_j).
\label{eq_quad_quad}
\end{align}
We refer to the number $\delta_n(f)$ defined via the formula
\begin{align}
\delta_n(f) = \left|
\int_{-1}^1 f(t) \; dt - \sum_{j=1}^n W_j \cdot f(t_j) \right|
\label{eq_quad_error_def}
\end{align}
as the "quadrature error".
\label{def_quad}
\end{definition}
%%%%%%%%%%%
%%%%%%%%%%%
%%%%%%%%%%%
\subsection{Quadrature Error and its Relation to $|\lambda_n|$}
\label{sec_quad_error}
Suppose now that $n$ is a positive integer, and that
$f:[-1,1] \to \Cc$ is an arbitrary bandlimited function
(with band limit $c$).
Suppose also that $S_n$ is the PSWF-based quadrature rule
of order $n$ (see \eqref{eq_quad_sn} in Definition~\ref{def_quad}).
One of the principal goals of this paper is to 
investigate the quadrature error $\delta_n(f)$ defined via
\eqref{eq_quad_error_def}.
The reader is referred to Section~\ref{sec_num_ill}
for the results of several related numerical experiments.
% The results of additional numerical experiments,
% in which this quadrature is used for the integration of
% certain functions, are 
% summarized in 
% Tables~\ref{t:test90},~\ref{t:test91} and
% Figures~\ref{fig:test92}, \ref{fig:test93log}
% (see Experiments 1 in Section~\ref{sec_exp12}).

The following theorem, illustrated in 
Table~\ref{t:test90}, provides an upper bound
on $\delta_n(\psi_m)$,
for any integer $m=0,\dots,n-1$.
%%%%%%%%%%%
This theorem is illustrated in
Table~\ref{t:test91} and in Figure~\ref{fig:test92}
(see Experiment 2 in Section~\ref{sec_exp12});
see also Conjecture~\ref{conj_quad_error}
and Remark~\ref{rem_conj} in Section~\ref{sec_exp12}.
%%%%%%%%%%%
\begin{thm}
Suppose that $c$ is a positive real number, and that
\begin{align}
c > 30.
\label{eq_quad_simple_c30}
\end{align}
Suppose also that $n>0$ and $0 \leq m \leq n-1$ are integers,
and that
\begin{align}
n > \frac{2c}{\pi} + 5.
\label{eq_quad_simple_n2c}
\end{align}
Suppose further that $\delta_n(\psi_m)$ is defined
via \eqref{eq_quad_error_def}.
Then,
\begin{align}
\delta_n(\psi_m) = 
\left| \int_{-1}^1 \psi_m(s)\;ds - \sum_{j=1}^n W_j\cdot\psi_m(t_j)\right| \leq
|\lambda_n| \cdot \left(
24 \cdot \log\left( \frac{1}{|\lambda_n|} \right) +
6 \cdot \chi_n
\right),
\label{eq_quad_simple_thm}
\end{align}
where
$\lambda_n, \chi_n$ are those of \eqref{eq_prolate_integral},
\eqref{eq_prolate_ode} in Section~\ref{sec_pswf}, respectively.
\label{thm_quad_simple}
\end{thm}
%%%%%%%%%%%%
%%%%%%%%%%%%
%%%%%%%%%%%%%
%%%%%%%%%%%%
%%%%%%%%%%%
\subsection{Quadrature Error and its Relation to $n$ and $c$}
\label{sec_main_result}
In Theorem~\ref{thm_quad_simple}, we established an upper bound
on the quadrature error $\delta_n(\psi_m)$
(see \eqref{eq_quad_error_def}
and \eqref{eq_quad_simple_thm} in Theorem~\ref{thm_quad_simple}).
However, this
bound depends on $\chi_n$ and $\lambda_n$. In particular, it is not obvious
how large $n$ should be to make sure that the quadrature error
does not exceed a prescribed $\varepsilon>0$. In this subsection,
we eliminate this inconvenience.

The following theorem is illustrated in Table~\ref{t:test178}
(see Experiment 3 in Section~\ref{sec_exp14}).
%%%%%%%%%%%%
\begin{thm}
Suppose that $c, \varepsilon$ are positive real numbers such that
\begin{align}
c > 30
\label{eq_quad_eps_large_c30}
\end{align}
and
\begin{align}
0 < \log \frac{1}{\varepsilon} < 
\frac{5 \cdot \pi}{4\sqrt{6}} \cdot c - 3 \cdot \log(c) - \log(6^5 \cdot 14340).
\label{eq_quad_eps_large_eps}
\end{align}
Suppose also that the real numbers $\alpha, \nu(\alpha)$ are defined
via the formulae
\begin{align}
\alpha = \frac{4\sqrt{6}}{\pi} \cdot
\left(
\log \frac{1}{\varepsilon} + 3 \cdot \log(c) + \log(6^5 \cdot 14340)
\right)
\label{eq_quad_eps_large_alpha}
\end{align}
and 
\begin{align}
\nu(\alpha) = \frac{2c}{\pi} + \frac{\alpha}{2\pi} \cdot
   \log\left( \frac{16ec}{\alpha} \right),
\label{eq_quad_eps_large_nu}
\end{align}
respectively.
Suppose furthermore that
$n>0$ and $0 \leq m \leq n-1$ are integers such that
\begin{align}
n > \nu(\alpha),
\label{eq_quad_eps_large_n_nu}
\end{align}
and that $\delta_n(\psi_m)$ is defined
via \eqref{eq_quad_error_def}.
Then,
\begin{align}
\delta_n(\psi_m) = 
\left| \int_{-1}^1 \psi_m(s) \; ds - \sum_{j=1}^n \psi_m(t_j) W_j \right|
<
\varepsilon.
\label{eq_quad_eps_large_thm}
\end{align}
\label{thm_quad_eps_large}
\end{thm}
%
%%%%%%%%%%%
The following theorem is a direct consequence of
Theorem~\ref{thm_quad_eps_large}. This theorem
is one of the principal results of the paper.
It is illustrated in Table~\ref{t:test178}
(see Experiment 3 in Section~\ref{sec_exp14}).
See also Conjecture~\ref{conj_quad_error}
in Section~\ref{sec_exp12}.
\begin{thm}
Suppose that $c,\varepsilon$ are positive real numbers such that
\begin{align}
c > 60
\label{eq_quad_eps_simple_c60}
\end{align}
and
\begin{align}
0 < \varepsilon < 1.
\label{eq_quad_eps_simple_eps}
\end{align}
Suppose also that $n>0$ and $0 \leq m < n$ are integers, 
and that
\begin{align}
n > \frac{2c}{\pi} +
\left(10 + \frac{3}{2} \cdot \log(c) + 
   \frac{1}{2} \cdot \log\frac{1}{\varepsilon}
\right) \cdot \log\left( \frac{c}{2} \right).
\label{eq_quad_eps_simple_n}
\end{align}
Suppose furthermore that $\delta_n(\psi_m)$ is defined
via \eqref{eq_quad_error_def} in Definition~\ref{def_quad}. Then,
\begin{align}
\delta_n(\psi_m) = 
\left| \int_{-1}^1 \psi_m(s) \; ds - \sum_{j=1}^n \psi_m(t_j) W_j \right|
<
\varepsilon.
\label{eq_quad_eps_simple_thm}
\end{align}
\label{thm_quad_eps_simple}
\end{thm}

%%%%%%%%%%%%%%%%%%%%%%%%%%%%%%%%%%%%%%%%%%%%%%%%
\subsection{Quadrature Weights}
\label{sec_weights}
In this subsection, we analyze the
weights of the quadrature rule $S_n$
(see \eqref{eq_quad_w}, \eqref{eq_quad_sn} 
in Section~\ref{sec_quad}).
This analysis has two principal purposes. On the one hand,
it provides the basis for a fast algorithm for the evaluation
of the weights. On the other hand, it provides an explanation
of some empirically observed properties of the weights.

The results of this subsection are 
illustrated in Table~\ref{t:test96} and in Figure~\ref{fig:test96}
(see Experiment 4 in Section~\ref{sec_exp15}).

The following theorem is instrumental for the evaluation
of the quadrature weights $W_1,\dots,W_n$ 
(see \eqref{eq_quad_w} in Definition~\ref{def_quad}).
\begin{thm}
Suppose that $n \geq 0$ is an integer, and that
the function $\tilde{\Phi}_n: (-1,1) \to \Rc$
is defined via the formula
\begin{align}
\tilde{\Phi}_n(t) = \sum_{k = 0}^{\infty} \alpha_k^{(n)} \cdot Q_k(t),
\label{eq_num_tilde_phi_def}
\end{align}
where $Q_k(t)$ and $\alpha_k^{(n)}$ are defined, respectively,
via \eqref{eq_legendre_fun_0_1}, \eqref{eq_legendre_fun_rec}
and \eqref{eq_num_leg_alpha_knc}
in Section~\ref{sec_legendre}
(compare to \eqref{eq_num_leg_exp} in Section~\ref{sec_legendre}).
Then, for every integer $j=1,\dots,n$,
\begin{align}
W_j = -\frac{2}{\psi_n'(t_j)} 
\sum_{k = 0}^{\infty} \alpha_k^{(n)} \cdot Q_k(t_j)
=
- 2 \cdot \frac{ \tilde{\Phi}_n(t_j) }{ \psi_n'(t_j) },
\label{eq_tilde_phi_w}
\end{align}
where $t_1,\dots,t_n$ and $W_1,\dots,W_n$ are, respectively,
the nodes and weights 
of the quadrature rule $S_n$ in Definition~\ref{def_quad}.
\label{thm_tilde_phi}
\end{thm}
Theorem~\ref{thm_tilde_phi} is illustrated in
Table~\ref{t:test96}. We observe that Theorem~\ref{thm_tilde_phi}
describes a connection
between the weights $W_1, \dots, W_n$ and the values
of $\tilde{\Phi}_n$ at $t_1, \dots, t_n$, where the function $\tilde{\Phi}_n$
is defined via \eqref{eq_num_tilde_phi_def}.

The following theorem states that $\tilde{\Phi}_n$ satisfies
a certain second-order non-homogeneous ODE, closely related
to the prolate ODE \eqref{eq_prolate_ode} in Section~\ref{sec_pswf}.
In particular, a recurrence
relation between the derivatives of $\tilde{\Phi}_n$ of
arbitrary order is established (compare to Theorem~\ref{thm_dpsi_reck} in
Section~\ref{sec_pswf}).
%%%%%%%%%%%
\begin{thm}
Suppose that $n \geq 0$ is an integer, and that
the function $\tilde{\Phi}_n : (-1,1) \to \Rc$
is defined via \eqref{eq_num_tilde_phi_def}. 
Suppose also that
the real numbers $\alpha_0^{(n)}, \alpha_1^{(n)}$ are 
defined via \eqref{eq_num_leg_alpha_knc} in Section~\ref{sec_legendre}.
Then, 
\begin{align}
%L_n \left[ \tilde{\Phi}_n \right] (t) = 
(1-t^2) \cdot \tilde{\Phi}_n''(t) - 2t \cdot \tilde{\Phi}_n'(t) + 
(\chi_n-c^2 t^2) \cdot \tilde{\Phi}_n(t) =
-c^2 \left( \alpha_0^{(n)} t + \alpha_1^{(n)} / 3 \right),
\label{eq_num_tilde_phi_ode}
\end{align}
for all real $-1 < t < 1$. Also,
\begin{align}
\left(1 - t^2\right) \cdot  \tilde{\Phi}_n'''(t) - 
4t \cdot \tilde{\Phi}_n''(t) +
\left(\chi_n - c^2 t^2 - 2\right) \cdot \tilde{\Phi}_n'(t) - 
2c^2 t \cdot \tilde{\Phi}_n(t) = 
-c^2 \alpha_0^{(n)},
\label{eq_num_dtilde_phi_3}
\end{align}
for all real $-1 < t < 1$. Finally,
\begin{align}
& \left(1 - t^2\right) \tilde{\Phi}_n^{(k+2)}(t) 
  - 2 \left(k+1\right) t \tilde{\Phi}_n^{(k+1)}(t)
  + \left(\chi_n - k\left(k+1\right) - c^2 t^2\right) 
 \tilde{\Phi}_n^{(k)}(t) \nonumber \\
& \quad \quad 
  -c^2 k t \tilde{\Phi}_n^{(k-1)}(t) 
  -c^2 k \left(k-1\right) \tilde{\Phi}_n^{(k-2)}(t) = 0,
\label{eq_dphi_reck}
\end{align}
for every integer $k \geq 2$ and all real $-1 < t < 1$
(compare to \eqref{eq_dpsi_reck} in Section~\ref{sec_pswf}).
\label{lem_tilde_phi_ode}
\end{thm}

In the following theorem, we establish the positivity
of the weights of the quadrature rule $S_n$
in Definition~\ref{def_quad}.
%%%%
\begin{thm}
Suppose that $c$ is a positive real number, and that
\begin{align}
c > 30.
\label{eq_positive_w_c30}
\end{align}
Suppose also that $n$ is a positive integer, and that
\begin{align}
n > \frac{2c}{\pi} + 5 \cdot \log(c) \cdot \log\left(\frac{c}{2}\right).
\label{eq_positive_w_n}
\end{align}
Suppose further that $W_1, \dots, W_n$ are defined
via \eqref{eq_quad_w}. 
Then, for all integer $j=1,\dots,n$,
\begin{align}
W_j > 0.
\label{eq_positive_w_thm}
\end{align}
\label{thm_positive_w}
\end{thm}
\begin{remark}
Extensive numerical experiments 
(see e.g. Table~\ref{t:test96} and Figure~\ref{fig:test96})
seem to indicate that the assumption
\eqref{eq_positive_w_n} is unnecessary. In other words,
the weights $W_1, \dots, W_n$ are always positive,
even for small values of $n$
(at the present time we do not have the proof of this fact).
\label{rem_w_always_pos}
\end{remark}
\begin{remark}
It was observed in \cite{Report3}, \cite{Report3Arxiv}
that, if $1 \leq j,k \leq n$
are integers, then
\begin{align}
\left( \psi_n'(t_j) \right)^2 \cdot (1-t_j^2) \cdot W_j =
\left( \psi_n'(t_k) \right)^2 \cdot (1-t_k^2) \cdot W_k +
O\left( |\lambda_n| \right)
\label{eq_w_approx}
\end{align}
(see also Experiment 4 in Section~\ref{sec_exp15}).
We observe that as $c \to 0$
the quadrature rule in Definition~\ref{def_quad} 
converges to the well known Gaussian quadrature rule,
whose nodes are the roots $t_1,\dots,t_n$ of the Legendre polynomial $P_n$
(see Section~\ref{sec_legendre}), and whose weights are defined via
the formula
\begin{align}
W_j = \frac{ 2 }{ \left(P_n'(t_j)\right)^2 \cdot \left(1 - t_j^2\right) },
\label{eq_gauss_weights}
\end{align}
for every $j=1,\dots,n$
(see e.g. \cite{Abramovitz}, Section 25.4). 
Thus, \eqref{eq_w_approx} is not
surprising. 
\label{rem_w_approx}
\end{remark}

%%%%%%%%%%%%%%%%%%%%%%%%%%%%%%%%%%%%%%%%%%%%%
%%%%%%%%%%%%%%%%%%%%%%%%%%%%%%%%%%%%%%%%%%%%%
\section{Numerical Algorithms}
\label{sec_num_algo}
In this section, we describe several numerical algorithms
for the evaluation of the PSWFs, certain related quantities,
and the quadrature rules
defined in % Definition~\ref{def_quad} in
Section~\ref{sec_quad}.
Throughout this section, the band limit $c > 0$ is a real number,
and the prolate index $n \geq 0$ is
a non-negative integer.

%%%%%%%%%%%%%%%%%%%%%%%%%%%%%%%%%%%%%%%%%%%%%%%%%%%%%
\subsection{Evaluation of $\chi_n$ and 
$\psi_n(x)$, $\psi_n'(x)$ for $-1 \leq x \leq 1$}
\label{sec_evaluate_beta}

The use of the expansion of $\psi_n$ into a Legendre series
(see \eqref{eq_num_leg_exp} in Section~\ref{sec_legendre})
for the evaluation of $\psi_n$
in the interval $[-1,1]$ goes back at least to 
the classical Bouwkamp algorithm (see \cite{Bouwkamp}).
More specifically, the coefficients 
$\beta_0^{(n)}, \beta_1^{(n)}, \dots$
of the Legendre expansion are precomputed first 
(see \eqref{eq_num_leg_beta_knc},
\eqref{eq_num_leg_alpha_knc} in Section~\ref{sec_legendre}).
These coefficients decay superalgebraically; in particular,
relatively few terms of the infinite sum \eqref{eq_num_leg_exp}
are required to evaluate $\psi_n$ to essentially
machine precision (see Section~\ref{sec_legendre},
in particular Theorem~\ref{thm_tridiagonal} and Remark~\ref{rem_tridiagonal},
and also \cite{RokhlinXiaoProlate} for more details).

%%%%%%%%%%%%%%%%%%%%%%%%%%%%%%
\subsubsection{Evaluation of $\chi_n$ and 
$\beta^{(n)}_0, \beta^{(n)}_1, \dots$}
\label{sec_evalaute_chi_beta}
Suppose now that $n \geq 0$, and one is interested in evaluating
the coefficients $\beta^{(m)}_0, \beta^{(m)}_1, \dots$ 
in \eqref{eq_num_leg_exp},
% of the Legendre expansion \eqref{eq_num_leg_exp} of $\psi_m$, 
for every
integer $0 \leq m \leq n$. This can be achieved by solving two 
$N \times N$ symmetric
tridiagonal eigenproblems, where $N$ is of order $n$
(see Theorem~\ref{thm_tridiagonal} and Remark~\ref{rem_tridiagonal}
in Section~\ref{sec_legendre}, and also
\cite{RokhlinXiaoProlate} for more details about this algorithm). 
In addition, this algorithm evaluates $\chi_0, \dots, \chi_n$.
Once this precomputation 
is done, for every integer $0 \leq m \leq n$ and for every 
real $-1 \leq x \leq 1$ one can evaluate $\psi_m(x)$ in $O(n)$ operations,
by computing the sum \eqref{eq_num_leg_exp}
(see, however, Remark~\ref{rem_o1_interpolation} below).

Suppose, on the other hand, that we are interested in 
a single PSWF only (as opposed to all the first $n$ PSWFs).
Obviously, we can use the algorithm above; however,
its cost is $O(n^2)$ operations (see Remark~\ref{rem_tridiagonal}
in Section~\ref{sec_legendre}).
In the rest of this subsection, we describe a procedure
for the evaluation of $\beta_0^{(n)}, \beta_1^{(n)}, \dots$ and $\chi_n$,
whose cost is $O(n + c \log(c))$ operations.

This algorithm is also based on Theorem~\ref{thm_tridiagonal}
in Section~\ref{sec_legendre}. It consists of two principal steps.
First, we compute a low-accuracy approximation $\tilde{\chi}_n$
of $\chi_n$, by means of Sturm Bisection
(see Section~\ref{sec_sturm}, \eqref{eq_a_even_eig},
\eqref{eq_a_odd_eig} and Remark~\ref{rem_tridiagonal}
in Section~\ref{sec_legendre},
and also \cite{Wilkinson2}).
Second, we compute $\chi_n$ and $\beta^{(n)}$ 
(see \eqref{eq_beta_n} and Remark~\ref{rem_n_even_odd}
in Section~\ref{sec_legendre})
by means of the Shifted Inverse Power Method (see Section~\ref{sec_power},
and also \cite{Wilkinson}, \cite{Dahlquist}).
The Shifted Inverse Power Method requires an initial approximation to
the eigenvalue;
for this purpose we use $\tilde{\chi}_n$.

Below is a more detailed description of these two steps.
\paragraph{Step 1 (initial approximation $\tilde{\chi}_n$ of $\chi_n$).}
Suppose that the infinite symmetric tridiagonal matrices $A^{even}$
and $A^{odd}$ are defined, respectively, via 
\eqref{eq_a_even}, \eqref{eq_a_odd} in Section~\ref{sec_legendre}.
Suppose also that $A^{(n)}$ is the $N \times N$ upper left square
submatrix of $A^{even}$, if $n$ is even, or of $A^{odd}$, if $n$ is odd. \\
{\bf Comment.} $N$ is an integer of order $n$ 
(see Remark~\ref{rem_tridiagonal} in Section~\ref{sec_legendre}). The choice
\begin{align}
N = 1.1 \cdot c + n + 1000
\label{eq_n_choice}
\end{align}
is sufficient for all practical purposes.
\begin{itemize}
\item use Theorems~\ref{thm_n_and_khi}, \ref{thm_khi_crude}
and \ref{thm_n_khi_simple} in Section~\ref{sec_pswf}
to choose real numbers $x_0 < y_0$ such that
\begin{align}
x_0 < \chi_n < y_0.
\label{eq_khi_brackets}
\end{align}
{\bf Comment.} For a more detailed discussion of lower and upper 
bounds
on $\chi_n$, see, for example, \cite{ReportACHA}, \cite{ReportArxiv}.
See also Remark~\ref{rem_sturm_cost} below.
\item use Sturm Bisection (see Section~\ref{sec_sturm}) with
initial values $x_0, y_0$ to compute $\tilde{\chi}_n$.
On each step of Sturm Bisection, the Sturm sequence
(see \eqref{eq_sturm_p_seq} in Theorem~\ref{thm_sturm}) is computed based on
the matrix $A^{(n)}$ (see above). \\
{\bf Comment.} 
In principle, Sturm Bisection can be used to evaluate $\chi_n$ 
to machine precision. However, the convergence rate of Sturm
Bisection is linear, and each iteration requires order $n$ 
operations (see Remark~\ref{rem_sturm} in Section~\ref{sec_sturm}).
On the other hand, the convergence rate of the Shifted Inverse Power Method
is cubic in the vicinity of the solution, while each iteration
requires also order $n$ operations
(see Remark~\ref{rem_power} in Section~\ref{sec_power}).
Thus, we use Sturm Bisection
to compute a low-order approximation $\tilde{\chi_n}$ to $\chi_n$,
and then refine it by the Shifted Inverse Power Method to obtain
$\chi_n$ to machine precision.

\end{itemize}
\begin{remark}
The use of Sturm Bisection as a tool to compute the eigenvalues of
a symmetric tridiagonal matrix goes back at least to 
\cite{Wilkinson2}; in the context of PSWFs, it appears in \cite{Hodge}.
\label{rem_sturm_first}
\end{remark}

The cost analysis of Step 1
relies on the following observation
based on 
Theorems~\ref{thm_prolate_ode},
\ref{thm_n_and_khi},
\ref{thm_khi_crude},
\ref{thm_n_khi_simple} in Section~\ref{sec_pswf}.

{\bf Observation 1.} Suppose that $n \geq 0$ is an integer.

If $0 \leq n < 2c/\pi$, then %(e.g. it seems $0 \leq \chi_0 \leq c$)
\begin{align}
\chi_{n+1} - \chi_n = O(c).
\label{eq_delta_khi_small}
\end{align} 

If $n > 2c/\pi$, then
\begin{align}
\chi_{n+1} - \chi_n = O(n).
\label{eq_delta_khi_large}
\end{align}
\begin{remark}
Due to Theorems~\ref{thm_n_and_khi},
\ref{thm_khi_crude} in Section~\ref{sec_pswf},
the inequality
\begin{align}
n \cdot(n+1) < \chi_n < c^2
\label{eq_khi_small_bounds}
\end{align}
holds for any real $c>0$ and all integer $0 \leq n < 2c/\pi$.
In this case, Step 1 requires
$O(c \cdot \log(c))$ operations,
due to the combination of 
\eqref{eq_delta_khi_small}, \eqref{eq_khi_small_bounds}
and Remark~\ref{rem_sturm} in Section~\ref{sec_sturm}.
On the other hand, if $n > 2c/\pi$, then 
the cost of Step 1 is $O(n)$ operations, due to the combination of
Theorems~\ref{thm_n_and_khi}, \ref{thm_n_khi_simple},
Remark~\ref{rem_sturm} in Section~\ref{sec_sturm} and
\eqref{eq_delta_khi_large}.
\label{rem_sturm_cost}
\end{remark}

\paragraph{Step 2 (evaluation of $\chi_n$ and $\beta^{(n)}$).}
Suppose now that $\tilde{\chi}_n$ is an approximation to $\chi_n$
evaluated in Step 1.
Suppose also that the integer $N$ is defined via
\eqref{eq_n_choice} above
(see also Remark~\ref{rem_tridiagonal} in Section~\ref{sec_legendre}).
%
% and that $\beta^{(n)} \in \Rc^N$ is that of
% Remark~\ref{rem_n_even_odd}
% in Section~\ref{sec_legendre} (see also \eqref{eq_beta_n}
% and Remark~\ref{rem_tridiagonal} in Section~\ref{sec_legendre}).
%
\begin{itemize}
\item generate a pseudorandom 
vector $\tilde{\beta} \in \Rc^N$ of unit length. \\
{\bf Comment.}
We use $\tilde{\chi}_n$ and $\tilde{\beta}$ as initial approximations
to the eigenvalue $\chi_n$ and the corresponding eigenvector, respectively,
for the Shifted Inverse Power Method (see Section~\ref{sec_power}).
\item conduct Shifted 
Inverse Power Method iterations until $\chi_n$ is evaluated
to machine precision. The corresponding eigenvector of unit length
is denoted by $\hat{\beta}^{(n)}$.\\
{\bf Comment.} Each Shifted Inverse Power iteration 
costs $O(N)$ operations,
and essentially $O(1)$ iterations are required (see 
Remark~\ref{rem_power} in Section~\ref{sec_power} for
more details).
In practice, in double precision calculations 
the number of iterations is usually between three and five.
\end{itemize}
\begin{remark}
Clearly, the cost of Step 2 is $O(n)$ operations
(see Remark~\ref{rem_tridiagonal} in Section~\ref{sec_legendre}
and Remark~\ref{rem_power} in Section~\ref{sec_power}).
\label{rem_inverse_cost}
\end{remark}
\begin{remark}
Suppose that the coordinates of the vector $\beta^{(n)} \in \Rc^N$ 
are defined via \eqref{eq_beta_n} 
(see also 
Remark~\ref{rem_n_even_odd} in Section~\ref{sec_legendre}). 
Then, $\hat{\beta}^{(n)}$ (evaluated
in Step 2 above) approximates $\beta^{(n)}$ to essentially
machine precision
(this is a well known property of the Inverse Power Method;
see Section~\ref{sec_power}, and also
\cite{Wilkinson}, \cite{Dahlquist} for more details).
In other words, 
\begin{align}
\| \hat{\beta}^{(n)} - \beta^{(n)} \| \leq 
\varepsilon \cdot \| \beta^{(n)} \| = \varepsilon,
\label{eq_beta_vector_acc}
\end{align}
where $\varepsilon$ is the machine accuracy
(e.g. 
$\varepsilon \approx \mbox{\text{\rm{1D-16}}}$ for
double precision calculations).
In addition, the eigenvalue $\chi_n$ is also evaluated
to relative accuracy $\varepsilon$.
\label{rem_beta_absolute}
\end{remark}

\subsubsection{Evaluation of $\psi_n(x)$, $\psi_n'(x)$ for $-1\leq x \leq 1$,
given $\chi_n$ and $\beta^{(n)}_0, \beta^{(n)}_1, \dots$}
\label{sec_evaluate_psi}
Suppose now that $\chi_n$ and the coefficients 
$\beta^{(n)}_0, \beta^{(n)}_1, \dots$ 
defined via
\eqref{eq_num_leg_beta_knc}
have already been evaluated.
Suppose also that the integer $N$ is defined via \eqref{eq_n_choice} above.

For any real $-1 \leq x \leq 1$, 
we evaluate $\psi_n(x)$ via the formula
\begin{align}
\psi_n(x) = \sum_{k=0}^{2N} P_k(x) \cdot \alpha^{(n)}_k
          = \sum_{k=0}^{2N} P_k(x) \cdot \beta^{(n)}_k \cdot \sqrt{k+1/2}
\label{eq_evaluate_psi}
\end{align} 
(compare to \eqref{eq_num_leg_exp} in Section~\ref{sec_legendre}).
Also, we evaluate $\psi_n'(x)$ via the formula
\begin{align}
\psi_n'(x) = \sum_{k=1}^{2N} P_k'(x) \cdot \alpha^{(n)}_k
          = \sum_{k=0}^{2N} P_k'(x) \cdot \beta^{(n)}_k \cdot \sqrt{k+1/2}.
\label{eq_evaluate_dpsi}
\end{align} 
\begin{remark}
Due to the combination of
Remark~\ref{rem_tridiagonal}
in Section~\ref{sec_legendre} and
Remark~\ref{rem_beta_absolute} above,
both $\psi_n(x)$ and $\psi_n'(x)$
are evaluated via \eqref{eq_evaluate_psi},
\eqref{eq_evaluate_dpsi} essentially to machine precision,
for any real $-1 \leq x \leq 1$
(also see \cite{RokhlinXiaoProlate} for more details).
\label{rem_psi_accuracy}
\end{remark}
\begin{remark}
Due to Remarks~\ref{rem_sturm_cost}, \ref{rem_inverse_cost} above,
the cost of the evaluation of $\chi_n$ and 
$\beta^{(n)}_0, \beta^{(n)}_1, \dots$ via Steps 1,2 is 
$O\left(n + c \log c\right)$ 
operations.
Once this precomputation has been performed, the cost of each
subsequent evaluation of $\psi_n(x)$, $\psi_n'(x)$, for any real
$-1 \leq x \leq 1$, is $O(n)$ operations,
according to \eqref{eq_evaluate_psi},
\eqref{eq_evaluate_dpsi} and Remark~\ref{rem_legendre_evaluate}
in Section~\ref{sec_legendre}.
\label{rem_total_cost}
\end{remark}
\begin{remark}
Once $\chi_n$ and $\beta^{(n)}_0, \beta^{(n)}_1, \dots$ have been
evaluated, one does not have to use \eqref{eq_evaluate_psi},
\eqref{eq_evaluate_dpsi}, to compute $\psi_n(x)$, $\psi_n'(x)$
at an arbitrary point $x$ in $[-1,1]$. Instead, the cost of evaluating,
say, $\psi_n(x)$ can be brought down from $O(n)$ to $O(1)$
(see Remark~\ref{rem_newton_interpolation} 
in Section~\ref{sec_evaluate_nodes}).
\label{rem_o1_interpolation}
\end{remark}

%%%%%%%%%%%%%%%%%%%%%%%%%%%%%%%%%%%%%%
\subsection{Evaluation of $\lambda_n$}
\label{sec_evaluate_lambda}
Suppose now that $n \geq 0$ is an integer, and that
one needs to evaluate the eigenvalue
$\lambda_n$ of the integral operator
$F_c$ (see \eqref{eq_pswf_fc} in Section~\ref{sec_pswf}).
Due to the combination of
\eqref{eq_pswf_fc} and Theorem~\ref{thm_pswf_main} in Section~\ref{sec_pswf}, 
if $n$ is even, then $\psi_n(0) \neq 0$,
and
\begin{align}
\lambda_n = \frac{1}{\psi_n(0)} \int_{-1}^1 \psi_n(t) \; dt;
\label{eq_eval_lambda_even}
\end{align}
for odd $n$,
\begin{align}
\lambda_n = \frac{ic}{\psi_n'(0)} \int_{-1}^1 t \cdot \psi_n(t) \; dt.
\label{eq_eval_lambda_odd}
\end{align}
The formulae \eqref{eq_eval_lambda_even} and \eqref{eq_eval_lambda_odd}
provide an obvious way to calculate $\lambda_n$ for
even and odd $n$, respectively, via numerical integration.
In fact, when $|\lambda_n|$ is relatively large, such
procedure is quite satisfactory. More specifically,
if $n < 2c/\pi$, then $|\lambda_n| \approx \sqrt{2\pi/c}$,
and $\lambda_n$ can be calculated via
 \eqref{eq_eval_lambda_even}, \eqref{eq_eval_lambda_odd}
to high relative precision
(see Theorems~\ref{thm_mu_spectrum}, \ref{thm_crude_inequality}
in Section~\ref{sec_pswf} and Remark~\ref{rem_psi_accuracy}
in Section~\ref{sec_evaluate_beta};
see also \cite{RokhlinXiaoProlate} for more details).
On the other hand, we observe that 
$\| \psi_n \|_{L^2[-1,1]} = 1$,
due to Theorem~\ref{thm_pswf_main}
in Section~\ref{sec_pswf}. As a result, 
when $|\lambda_n|$ is small, the formulae 
\eqref{eq_eval_lambda_even}, \eqref{eq_eval_lambda_odd}
are unsuitable for the evaluation of $\lambda_n$
via numerical integration, due to catastrophic cancellation.
For example, if $|\lambda_n|<\varepsilon$, where $\varepsilon$
is the machine precision, the formulae
\eqref{eq_eval_lambda_even}, \eqref{eq_eval_lambda_odd}
produce no correct digits at all.

The standard way to overcome this obstacle for numerical
evaluation of small $\lambda_n'$s is 
to calculate all 
the ratios $\lambda_0/\lambda_1, \dots, \lambda_n/\lambda_{n-1}$
(see, for example,
\cite{ProlateLandau1},
\cite{ProlateSlepian1},
\cite{ProlateSlepian2}); this turns out to be a well-conditioned
numerical procedure (see \cite{RokhlinXiaoProlate} for more details).
Then, $\lambda_0$ is evaluated via \eqref{eq_eval_lambda_even} above,
and
the eigenvalues $\lambda_1, \dots, \lambda_n$ are evaluated via the formula
\begin{align}
\lambda_m = \lambda_0 \cdot \frac{\lambda_1}{\lambda_0} \cdot \dots
   \cdot \frac{\lambda_m}{\lambda_{m-1}},
\label{eq_lambda_as_ratios}
\end{align}
for every integer $m=1,\dots,n$.

Suppose, on the other hand, that one is interested in a single $\lambda_n$ only
(as opposed to all the first $n$ eigenvalues).
Obviously, $\lambda_n$ can be evaluated via \eqref{eq_lambda_as_ratios}
from the ratios $\lambda_{j+1}/\lambda_j$,
as described above; however, it requires
at least $O(n^2)$ operations
(see \cite{RokhlinXiaoProlate}).

Unexpectedly, it turns out that $\lambda_n$ can be obtained 
to high relative accuracy in $O(1)$ operations
as a by-product of the algorithm described in Section~\ref{sec_evaluate_beta}.
More specifically, suppose that
the coefficients $\beta_0^{(n)}, \beta_1^{(n)}, \dots$
are defined via \eqref{eq_num_leg_beta_knc}.
We combine \eqref{eq_eval_lambda_even}, 
\eqref{eq_eval_lambda_odd} above with
\eqref{eq_prolate_integral},
\eqref{eq_legendre_pol_0_1},
\eqref{eq_legendre_normalized},
\eqref{eq_num_leg_beta_knc},
\eqref{eq_num_leg_alpha_knc}
to make the following observation.

{\bf Observation 1.}
If $n$ is even, then $\psi_n(0) \neq 0$, and
\begin{align}
\lambda_n = \frac{1}{\psi_n(0)} \int_{-1}^1 \psi_n(t) \; dt = 
%            \frac{2\alpha_0^{(n)}}{\psi_n(0)} =
            \frac{\beta_0^{(n)}  \sqrt{2}}{\psi_n(0)}.
\label{eq_num_lambda_even}
\end{align}

If $n$ is odd, then $\psi_n'(0) \neq 0$,
and 
\begin{align}
\lambda_n = \frac{ic}{\psi_n'(0)} \int_{-1}^1 t \cdot \psi_n(t) \; dt =
%            \frac{2}{3} \cdot \frac{ic \alpha_1^{(n)}}{ \psi_n'(0)} =
            \sqrt{\frac{2}{3}} \cdot \frac{ic \beta_1^{(n)}}{ \psi_n'(0)}.
\label{eq_num_lambda_odd}
\end{align}
\begin{remark}
Obviously, the cost of evaluating $\lambda_n$ from 
$\psi_n(0), \beta^{(n)}_0$ via \eqref{eq_num_lambda_even} (for even $n$)
or from $\psi_n'(0), \beta^{(n)}_1$
via \eqref{eq_num_lambda_odd} (for odd $n$)
is $O(1)$ operations.
\label{rem_lambda_cost}
\end{remark}
\begin{remark}
Due to Remarks~\ref{rem_total_cost}, \ref{rem_lambda_cost} and
\eqref{eq_num_lambda_even}, \eqref{eq_num_lambda_odd},
a single $\lambda_n$ can be evaluated as a by-product
of the procedure described in Section~\ref{sec_evaluate_beta},
at the total cost of $O\left(n + c \log(c)\right)$ operations.
\label{rem_lambda_cost_alone}
\end{remark}
Remarks~\ref{rem_lambda_cost}, \ref{rem_lambda_cost_alone} describe
the {\it cost} of the evaluation of $\lambda_n$
via \eqref{eq_num_lambda_even}, \eqref{eq_num_lambda_odd}.
To describe the {\it accuracy} of this procedure, we start
with the following observation.

{\bf Observation 2.}
Due to Remark~\ref{rem_psi_accuracy},
$\lambda_n$ is evaluated to the same relative accuracy as $\beta^{(n)}_0$
(for even $n$) or as $\beta^{(n)}_1$ (for odd $n$).
According to \eqref{eq_beta_vector_acc} in Remark~\ref{rem_beta_absolute},
the algorithm of Section~\ref{sec_evaluate_beta}
evaluates the vector $\beta^{(n)}$ to relative accuracy $\varepsilon$,
where $\varepsilon$ is the machine precision.
However, this means
that a single {\bf coordinate} of $\beta^{(n)}$ is only guaranteed to be
evaluated to
{\bf absolute} accuracy $\varepsilon$.
More specifically, the inequality
\begin{align}
\left| \frac{\beta^{(n)}_k - \hat{\beta}^{(n)}_k}{\beta^{(n)}_k} \right|
\leq
\frac{\varepsilon}{\left| \beta^{(n)}_k \right|}
\label{eq_beta_absolute}
\end{align}
holds for every integer $k=0,\dots,N$,
where $N$ is defined via \eqref{eq_n_choice} 
in Section~\ref{sec_evaluate_beta}, and $\hat{\beta}^{(n)}_k$
is the numerical approximation to $\beta^{(n)}_k$.
In general, the inequality \eqref{eq_beta_absolute} can 
be rather tight; as a result,
% due to \eqref{eq_beta_absolute}, 
if, for example, 
$|\beta^{(n)}_0| \leq \varepsilon/10$, 
then, apriori, we cannot
expect $\hat{\beta}^{(n)}_0$ to
approximate $\beta^{(n)}_0$ to any digit at all!

In practical computations, it is sometimes desirable to evaluate
extremely small $\lambda_n$'s 
(e.g. $|\lambda_n| \approx \mbox{\text{\rm{1D-50}}}$).
Observation 2
seems to suggest that, in such cases, the evaluation of
$\lambda_n$ via the procedure described above is futile
due to disastrous loss of accuracy.

%from $\beta^{(n)}_0, \beta^{(n)}_1$ via
%\eqref{eq_num_lambda_even}, \eqref{eq_num_lambda_odd} is useless.

Fortunately, it turns out that the algorithm described 
in Section~\ref{sec_evaluate_beta}
{\bf always} evaluates $\beta^{(n)}_0, \beta^{(n)}_1$ to high relative 
accuracy, regardless of how small they are. 
This is a consequence of a more general (and somewhat surprising!) 
phenomenon studied in
detail in \cite{Report4}, \cite{Report4Arxiv}. We summarize the 
corresponding results in the following theorem.
\begin{thm}
For a certain class of real symmetric tridiagonal matrices, the coordinates
of their eigenvectors are defined to high relative precision.
Moreover, the matrices $A^{even}, A^{odd}$ 
defined, respectively, via 
\eqref{eq_a_even}, \eqref{eq_a_odd} in Section~\ref{sec_legendre},
belong to this class.
\label{thm_phenomenon}
\end{thm}
In the following theorem, we summarize implications of
Theorem~\ref{thm_phenomenon} for the evaluation of $\beta^{(n)}_0, 
\beta^{(n)}_1$ via the algorithm in Section~\ref{sec_evaluate_beta}
(the proof of a slightly modified version of this theorem appears
in \cite{Report4}, \cite{Report4Arxiv}).
\begin{thm}
Suppose that $c>0$ is a real number, that $n \geq 0$ is an integer,
and that $\beta^{(n)}_0, \beta^{(n)}_1$ are defined
via \eqref{eq_num_leg_beta_knc} in Section~\ref{sec_legendre}.
Then, the algorithm described in Section~\ref{sec_evaluate_beta}
evaluates $\beta^{(n)}_0, \beta^{(n)}_1$ to high relative accuracy.
More specifically, 
\begin{align}
\left| \frac{\beta^{(n)}_0 - \hat{\beta}^{(n)}_0}{\beta^{(n)}_0} \right|
\leq
10 \cdot \varepsilon \cdot c
\label{eq_beta0_relative}
\end{align}
for even $n$, and
\begin{align}
\left| \frac{\beta^{(n)}_1 - \hat{\beta}^{(n)}_1}{\beta^{(n)}_1} \right|
\leq
10 \cdot \varepsilon \cdot c
\label{eq_beta1_relative}
\end{align}
for odd $n$,
where $\hat{\beta}^{(n)}_0, \hat{\beta}^{(n)}_1$ are 
the numerical approximation to $\beta^{(n)}_0, \beta^{(n)}_1$,
respectively,
and $\varepsilon$ is the machine accuracy
(e.g. 
$\varepsilon \approx \mbox{\text{\rm{1D-16}}}$ for
double precision calculations).
\label{thm_beta01}
\end{thm}
\begin{remark}
The algorithm described in Section~\ref{sec_evaluate_beta}
evaluates the eigenvectors $\beta^{(n)}$ by
the Shifted Inverse Power Method (see Section~\ref{sec_power}).
It turns out that the choice of method
is important in this situation: 
if, for example, these eigenvectors are evaluated via
the standard and well known Jacobi Rotations
(rather than Inverse Power),
the small coordinates exhibit the loss of accuracy expected from
\eqref{eq_beta_absolute} (see \cite{Report4}, \cite{Report4Arxiv}
for more details about this and related issues).
\label{rem_method}
\end{remark}
\begin{remark}
Due to the combination of Remark~\ref{rem_psi_accuracy} 
in Section~\ref{sec_evaluate_beta},
Observation 2 above and Theorem~\ref{thm_beta01},
the algorithm of this section evaluates $\lambda_n$ to high
relative accuracy. More specifically, 
at most $1+\log_{10}\left(c\right)$ decimal digits
are lost in the evaluation of $\lambda_n$.
\label{rem_lambda_acc}
\end{remark}
%
\begin{comment}
conduct additional $K$ iterations of inverse power method,
until the convergence of the first coordinate of $\hat{\beta}^{(n)}$. \\
{\bf Comment.} Both analysis and numerical experiments 
(to be reported at a later date) suggest that
\begin{align}
K = 1 + 
\text{ceil}\left(
\frac{\log\left( \left|\beta^{(n)}_0\right| + \left|\beta^{(n)}_1\right| 
 \right)}
{\log\left( \varepsilon \right)}
\right),
\label{eq_k_power}
\end{align}
where $\varepsilon$ is the machine precision (e.g. 
$\varepsilon \approx \mbox{\text{\rm{1D-16}}}$
for double precision calculations),
and $\text{ceil}(a)$ is the minimal integer greater than $a$,
for a real number $a$. For example, if 
$|\beta^{(n)}_0| \approx \mbox{\text{\rm{1D-99}}}$, and
$\varepsilon \approx \mbox{\text{\rm{1D-16}}}$,
then
$K = 8$.
In practice, $K$ does not to 
be known in advance; rather, we iterate until convergence.
\end{comment}

%%%%%%%%%%%%%%%%%%%%%%%%%%%%%%%%%%%%%%%%%%%%%%%
\subsection{Evaluation of the Quadrature Nodes}
\label{sec_evaluate_nodes}
Suppose that $n > 0$ is an integer, and that
the quadrature rule $S_n$ is defined via \eqref{eq_quad_sn}
in Section~\ref{sec_quad}. According to \eqref{eq_quad_t},
the nodes of $S_n$ are precisely the $n$ roots
$t_1,\dots,t_n$ of $\psi_n$ in the interval $(-1,1)$.
% In this subsection, we describe a numerical algorithm for
% the evaluation of $t_1,\dots,t_n$.
% Since $\psi_n$ is symmetric
% about the origin (see Theorem~\ref{thm_pswf_main}
% in Section~\ref{sec_pswf}), it suffices to evaluate the roots
% of $\psi_n$ in the interval $(0,1)$.

In this section, we describe a numerical procedure
for the evaluation of these quadrature nodes. This procedure
is based on the fast algorithm
for the calculation of the roots of special functions
described in \cite{Glaser}. It combines
Pr\"ufer's transformation (see Section~\ref{sec_prufer}),
Runge-Kutta method (see Section~\ref{sec_runge_kutta})
and Taylor's method (see Section~\ref{sec_taylor}).
This algorithm also evaluates
$\psi_n'(t_1), \dots, \psi_n'(t_n)$.
It requires $O(n)$ operations to 
compute all roots of $\psi_n$ in $(-1,1)$ as well as
the derivative of $\psi_n$ at these roots.

% The algorithm consists of the following steps.
% Suppose that $t_{\min}$ is the minimal root of $\psi_n$ in
% $[0,1)$. We start with evaluating $t_{\min}$ and 
% $\psi_n'\left(t_{\min}\right)$, using 

% We start with evaluating the minimal root $t_{\min}$ of
% $\psi_n$ in $[0,1)$. Then, we evaluate $\psi_n'(t_{\min})$.
% Next, we iteratively evaluate all the roots of $\psi_n$ in $(0,1)$
% and the derivative of $\psi_n$ at these roots. Finally, we use
% the symmetry of $\psi_n$ about the origin
% (see Theorem~\ref{thm_pswf_main} in Section~\ref{sec_pswf})
% to complete the procedure.

A short outline of the principal steps of the algorithm is provided below.
For a more detailed description of the algorithm and its properties,
the reader is referred to \cite{Glaser}.

Suppose that $t_{\min}$ is the minimal root of $\psi_n$ in $[0,1)$.

\paragraph{Step 1 (evaluation of $t_{\min}$).}
If $n$ is odd, then
\begin{align}
t_{\min} = t_{(n+1)/2} = 0,
\label{eq_t_min_odd}
\end{align}
due to Theorem~\ref{thm_pswf_main} in Section~\ref{sec_pswf}.
On the other hand, if $n$ is even, then
\begin{align}
t_{\min} = t_{(n+2)/2} > 0.
\label{eq_t_min_even}
\end{align}
To compute $t_{\min}$ in the case of even $n$,
we numerically solve the ODE \eqref{eq_ds_ode} with the initial condition
\eqref{eq_s_at_0} in the interval $\left[\pi n/2, \pi\cdot(n+1)/2\right]$,
by using 20 steps of Runge-Kutta method described in
Section~\ref{sec_runge_kutta}.
The rightmost value $\tilde{t}_{\min}$ of the solution is 
a low-order approximation of $t_{\min} = t_{(n+2)/2}$ (see
\eqref{eq_s_at_ti}, \eqref{eq_t_min_even}).
Then, we evaluate $t_{\min}$ to machine precision
via Newton's method (see Section~\ref{sec_newton}), using
$\tilde{t}_{\min}$ as an initial approximation to $t_{\min}$.
On each Newton iteration, we evaluate $\psi_n$ and $\psi_n'$
by using the algorithm of Section~\ref{sec_evaluate_beta}
(see \eqref{eq_evaluate_psi}, \eqref{eq_evaluate_dpsi}).

{\bf Observation 1.} The point $\tilde{t}_{\min}$ approximates
$t_{\min}$ to at least three decimal digits
(see Section~\ref{sec_runge_kutta}).
Since Newton's method converges quadratically in the vicinity
of the solution,
only several Newton iterations are required to obtain $t_{\min}$ from
$\tilde{t}_{\min}$
to essentially machine precision (see \cite{Glaser} for more details).
In our experience, the number of Newton iterations in
this step never exceeds four in double precision calculations
(and never exceeds six in extended precision calculations).
We combine this observation with 
Remark~\ref{rem_total_cost} in Section~\ref{sec_evaluate_beta}
to conclude that the total cost of Step 1 is $O(n)$ operations.

\paragraph{Step 2 (evaluation of $\psi_n'(t_{\min})$).}
We evaluate $\psi_n'(t_{\min})$ to machine precision 
via \eqref{eq_evaluate_dpsi} in Section~\ref{sec_evaluate_beta}.

{\bf Observation 2.} 
Due to Remark~\ref{rem_total_cost} in Section~\ref{sec_evaluate_beta}, 
the cost of Step 2 is $O(n)$ operations.

The remaining roots of $\psi_n$ in $(t_{\min}, 1)$ are computed 
one by one, as follows. Suppose that $n/2 < j < n$ is an integer,
and both $t_j$ and $\psi_n'(t_j)$ have already been evaluated.

\paragraph{Step 3 (evaluation of $t_{j+1}$ and $\psi_n'(t_{j+1})$,
given $t_j$ and $\psi_n'(t_j)$).}

\begin{itemize}
\item evaluate $\psi_n^{(2)}(t_j), \dots, \psi_n^{(M)}(t_j)$
via the recurrence relation \eqref{eq_dpsi_reck} in Section~\ref{sec_pswf}
(in double precision calculations, $M=30$; in extended precision calculations,
$M=60$).
% (see Theorem~\ref{thm_dpsi_reck} in Section~\ref{sec_pswf}).
\item use 20 steps of Runge-Kutta method (see Section~\ref{sec_runge_kutta}),
to solve the ODE \eqref{eq_ds_ode} with the initial condition
\begin{align}
s\left( \pi \cdot \left(j-\frac{1}{2}\right) \right) = t_j
\label{eq_s_eq_tj}
\end{align}
in the interval $\left[\pi \cdot(j- 1/2), \pi \cdot(j+ 1/2) \right]$
(see \eqref{eq_s_at_ti}).
The rightmost value $\tilde{t}_{j+1}$ of the solution is a low-order
approximation of $t_{j+1}$.
\item compute $t_{j+1}$ via Newton's method 
(see Section~\ref{sec_newton}), using
$\tilde{t}_{j+1}$ as the initial approximation to $t_{j+1}$.
On each Newton iteration, we evaluate $\psi_n$ and $\psi_n'$
via Taylor's method (see Section~\ref{sec_taylor}).
The Taylor expansion of appropriate order $M$ about $t_j$ is used, i.e.
\begin{align}
\psi_n(t) = \sum_{k=0}^{M} \frac{\psi_n^{(k)}(t_j)}{k!} \cdot (t-t_j)^k +
O\left( (t-t_j)^{M+1} \right).
\label{eq_t_taylor}
\end{align}
\item evaluate $\psi_n'(t_{j+1})$ via Taylor's method. The
Taylor expansion of order $M-1$ is used, i.e.
\begin{align}
\psi_n'(t_{j+1}) = 
\sum_{k=0}^{M-1} \frac{\psi_n^{(k+1)}(t_j)}{k!} \cdot (t_{j+1}-t_j)^k
+ O\left( (t_{j+1}-t_j)^M \right).
\label{eq_dpsi_taylor}
\end{align}
In both \eqref{eq_t_taylor} and \eqref{eq_dpsi_taylor},
we set
$M=30$ for double precision calculations, and $M=60$ for extended
precision calculations.
\end{itemize}
{\bf Observation 3.} The point $\tilde{t}_{j+1}$ approximates
$t_{j+1}$ to at least three decimal digits
(see Section~\ref{sec_runge_kutta}). Subsequently,
only several Newton iterations are required to obtain $t_{j+1}$
to essentially machine precision (see 
Observation 1 above, and also \cite{Glaser} for more details).
Thus the cost of Step 3 is $O(1)$ operations, for every 
integer $n/2 < j < n$.
\begin{remark}
Obviously, on each Newton iteration 
one can evaluate $\psi_n$ and $\psi_n'$ via \eqref{eq_evaluate_psi},
\eqref{eq_evaluate_dpsi} in Section~\ref{sec_evaluate_beta}
rather than via \eqref{eq_t_taylor}, \eqref{eq_dpsi_taylor}.
However, this would increase the cost 
of each such evaluation from $O(1)$ to $O(n)$,
and the total cost of the procedure from $O(n)$ to $O(n^2)$
(see Remark~\ref{rem_total_cost} in Section~\ref{sec_evaluate_beta}).
\label{rem_newton}
\end{remark}

\paragraph{Step 4 (evaluation of $t_j$ and $\psi_n'(t_j)$ for 
all $j \leq n/2$).}
Step 3 is repeated  for every integer $n/2 < j < n$. 
To evaluate $t_j$ and $\psi_n'(t_j)$ for $-1 < t_j < 0$,
we use the symmetry
of $\psi_n$ about zero (see Theorem~\ref{thm_pswf_main}
in Section~\ref{sec_pswf}). More specifically,
for every integer $1 \leq j \leq n/2$, we compute $t_j$ and $\psi_n'(t_j)$,
respectively,
via the formulae
\begin{align}
t_j = t_{n+1-j}
\label{eq_tj_symmetry}
\end{align} 
and
\begin{align}
\psi_n'(t_j) = (-1)^{n+1} \cdot \psi_n'(t_{n+1-j}).
\end{align}

\paragraph{Summary (evaluation of $t_j$ and $\psi_n'(t_j)$, for
all $j=1,\dots,n$).}
To summarize, 
the procedure for the evaluation of all roots of $\psi_n$ in $(-1,1)$
(as well as the derivative of $\psi_n$ at these roots)
is as follows:
\begin{itemize}
\item Evaluate $t_{\min}$ defined via \eqref{eq_t_min_odd},
\eqref{eq_t_min_even} (see Step 1). Cost: $O(n)$ operations.
\item Evaluate $\psi_n'(t_{\min})$ (see Step 2). Cost: $O(n)$ operations.
\item For every integer $n/2 < j < n$, 
evaluate $t_{j+1}$ and $\psi_n'(t_{j+1})$ (see Step 3). 
Cost: $O(n)$ operations.
\item For every integer $1 \leq j \leq n/2$, 
evaluate $t_j$ and $\psi_n'(t_j)$ (see Step 4).
Cost: $O(n)$ operations.
\end{itemize}
\begin{remark}
We observe that the algorithm described in this section not only
computes the roots $t_1, \dots, t_n$ of $\psi_n$ in $(-1,1)$,
but also evaluates $\psi_n'$ at all these roots.
The total cost of this algorithm is $O(n)$ operations,
and all the quantities are evaluated essentially to machine precision
(see Observations 1,2,3 above).
\label{rem_evaluate_nodes}
\end{remark}
\begin{remark}
The algorithm described in this section uses the quantities
$\chi_n$ and $\beta^{(n)}_0, \beta^{(n)}_1, \dots$ 
computed via the procedure of Section~\ref{sec_evaluate_beta}.
If $n < 2c/\pi$, then these quantities are obtained at the cost
of $O(n + c \log(c))$ operations;
if $n > 2c/\pi$, then these quantities are obtained at the cost
of $O(n)$ operations
(see Remarks~\ref{rem_sturm_cost}, \ref{rem_total_cost}
in Section~\ref{sec_evaluate_beta}).
\label{rem_evaluate_nodes2}
\end{remark}
\begin{remark}
As a by-product of the algorithm described in this section,
we obtain
a table of all the derivatives of $\psi_n$ up to order $M$ at all roots
of $\psi_n$ in $(-1,1)$ (here $M=30$ in double precision calculation,
and $M=60$ in extended precision calculations). In other words,
$\psi_n^{(k)}(t_j)$ are calculated for every $k=1,\dots,M$ and every
$j=1,\dots,n$ (see Step 3 above). This table can be used to evaluate
$\psi_n(x), \psi_n'(x)$ at an arbitrary
point $t_1 \leq x \leq t_n$ to essentially machine precision
in $O(1)$ operations via interpolation,
using the formulae
\eqref{eq_t_taylor}, \eqref{eq_dpsi_taylor} 
(see also Remark~\ref{rem_o1_interpolation} in Section~\ref{sec_evaluate_beta}).
% For $t_n \leq x \leq 1$, a similar Taylor expansion about $1$ can be used.
\label{rem_newton_interpolation}
\end{remark}

%%%%%%%%%%%%%%%%%%%%%%%%%%%%%%%%%%%%%%%%%%%%%%%%%
\subsection{Evaluation of the Quadrature Weights}
\label{sec_evaluate_weights}
Suppose now that $n > 0$ is an integer, and that
the quadrature rule $S_n$ is defined via \eqref{eq_quad_sn}
in Section~\ref{sec_quad}. 
In this subsection, we describe an algorithm for the evaluation
of the weights $W_1,\dots,W_n$ of this quadrature rule
(see \eqref{eq_quad_w} in Section~\ref{sec_quad}).
The results of this subsection are 
illustrated in Table~\ref{t:test96} and in Figure~\ref{fig:test96}
(see Experiment 4 in Section~\ref{sec_exp15}).

In the description of the algorithms below, we assume that
the coefficients $\beta^{(n)}_0, \beta^{(n)}_1, \dots$
(defined via \eqref{eq_num_leg_beta_knc} in Section~\ref{sec_legendre})
have already been evaluated (for example, by the algorithm
in Section~\ref{sec_evaluate_beta}).
In addition, we assume that the quadrature nodes $t_1, \dots, t_n$ as well as
$\psi_n'(t_1), \dots, \psi_n'(t_n)$ have also been computed
(for example, by the algorithm of Section~\ref{sec_evaluate_nodes}).

An obvious way to compute $W_1, \dots, W_n$ is to evaluate
\eqref{eq_quad_w}
numerically. However, due to \eqref{eq_quad_phi}, 
the integrand
$\varphi_j$ in \eqref{eq_quad_w} has
$n-1$ roots in $(-1,1)$,
for every $j=1,\dots,n$.
In particular, such approach is unlikely to require
less that $O(n^2)$ operations.
% (see also Section~\ref{sec_evaluate_beta}). 
% In addition, each $\varphi_j$
% has a singularity (albeit, removable) at $t_j$, which might be
% a nuisance for numerical integration, especially if high precision
% is required.

Rather than computing \eqref{eq_quad_w} directly, we evaluate
$W_1,\dots,W_n$ by using the results of Section~\ref{sec_weights}.
In the rest of this subsection, we describe two such algorithms;
both evaluate $W_1,\dots,W_n$ essentially to machine precision.
One of these algorithms (based on
Theorem~\ref{thm_tilde_phi}) 
is fairly straightforward; however, its cost 
is $O(n^2)$ operations. The other algorithm
(based on Theorem~\ref{lem_tilde_phi_ode}),
while still rather simple,
is also computationally efficient: its cost is $O(n)$ operations.

\paragraph{Algorithm 1: evaluation of $W_1, \dots, W_n$ in $O(n^2)$ operations.}
Suppose that the integer $N$ is defined via \eqref{eq_n_choice}
in Section~\ref{sec_evaluate_beta}.
For every integer $j=1,\dots,n$, we compute an approximation
$\widetilde{W}_j$ to $W_j$ via
the formula
\begin{align}
\widetilde{W}_j = 
-\frac{2}{\psi_n'(t_j)} 
\sum_{k = 0}^{2N} \alpha_k^{(n)} \cdot Q_k(t_j) =
-\frac{2}{\psi_n'(t_j)}
\sum_{k = 0}^{2N} \beta_k^{(n)} \cdot Q_k(t_j) \cdot \sqrt{k+1/2},
\label{eq_wj_as_sum}
\end{align}
where $Q_k(t)$ and $\alpha_k^{(n)}$ are defined, respectively,
via \eqref{eq_legendre_fun_0_1}, \eqref{eq_legendre_fun_rec}
and \eqref{eq_num_leg_alpha_knc}
in Section~\ref{sec_legendre}.
We observe that \eqref{eq_wj_as_sum} is obtained from
the identity
\eqref{eq_tilde_phi_w} in Theorem~\ref{thm_tilde_phi}
in Section~\ref{sec_weights}
by truncating the infinite series at $2N$ terms.
\begin{remark}
Due to the combination of 
Remarks~\ref{rem_tridiagonal}, \ref{rem_legendre_evaluate}
in Section~\ref{sec_legendre},
Remark~\ref{rem_beta_absolute}
in Section~\ref{sec_evaluate_beta}, \eqref{eq_n_choice} 
and Theorem~\ref{thm_tilde_phi},
each weight $W_j$ is evaluated via \eqref{eq_wj_as_sum}
essentially to machine precision
(see also Experiment 4 in Section~\ref{sec_exp15}).
\label{rem_algorithm1_acc}
\end{remark}
\begin{remark}
Due to the combination of Remark~\ref{rem_legendre_evaluate}
in Section~\ref{sec_legendre}
and
\eqref{eq_n_choice} in Section~\ref{sec_evaluate_beta},
the overall cost of computing $W_1,\dots,W_n$ via \eqref{eq_wj_as_sum}
is $O(n^2)$ operations.
\label{rem_algorithm1_cost}
\end{remark}

\paragraph{Algorithm 2: evaluation of $W_1, \dots, W_n$ in $O(n)$ operations.}
This algorithm is somewhat similar to the procedure
for the evaluation of the roots of $\psi_n$ in $(-1,1)$
described in Section~\ref{sec_evaluate_nodes}.

Suppose first that $t_{\min}$ is the minimal root of $\psi_n$ in $[0,1)$.
In other words,
\begin{align}
t_{\min} = 
\begin{cases}
t_{(n+1)/2} = 0& \text{ if } n \text{ is odd}, \\
t_{(n+2)/2} > 0& \text{ if } n \text{ is even}
\end{cases}
\label{eq_tmin_both}
\end{align}
(see 
\eqref{eq_t_min_odd}, \eqref{eq_t_min_even}
in Section~\ref{sec_evaluate_nodes}).
Suppose also that the function $\tilde{\Phi}_n: (-1,1) \to \Rc$ is defined
via \eqref{eq_num_tilde_phi_def} in Theorem~\ref{thm_tilde_phi}
in Section~\ref{sec_weights}.
\paragraph{Step 1 (evaluation of 
$\tilde{\Phi}_n(t_{\min})$ and $\tilde{\Phi}_n'(t_{\min})$).}
We evaluate $\tilde{\Phi}_n(t_{\min})$ 
and
$\tilde{\Phi}'_n(t_{\min})$ via the formulae
\begin{align}
\tilde{\Phi}_n(t_{\min}) = 
\sum_{k = 0}^{2N} \alpha_k^{(n)} \cdot Q_k(t_{\min}) =
\sum_{k = 0}^{2N} \beta_k^{(n)} \cdot Q_k(t_{\min}) \cdot \sqrt{k+1/2}
\label{eq_tilde_phi_as_sum}
\end{align}
and
\begin{align}
\tilde{\Phi}'_n(t_{\min}) =
\sum_{k = 0}^{2N} \alpha_k^{(n)} \cdot Q_k'(t_{\min}) =
\sum_{k = 0}^{2N} \beta_k^{(n)} \cdot Q_k'(t_{\min}) \cdot \sqrt{k+1/2},
\label{eq_dtilde_phi_as_sum}
\end{align}
respectively
(see \eqref{eq_wj_as_sum} in the description of
Algorithm 1 above).
Observe that \eqref{eq_tilde_phi_as_sum},
\eqref{eq_dtilde_phi_as_sum} are obtained from the infinite
expansion \eqref{eq_num_tilde_phi_def} in Theorem~\ref{thm_tilde_phi}
by truncation.
\begin{remark}
Due to Remarks~\ref{rem_algorithm1_acc}, \ref{rem_algorithm1_cost},
the cost of Step 1 is $O(n)$ operations; moreover,
$\tilde{\Phi}_n(t_{\min})$ and $\tilde{\Phi}_n'(t_{\min})$
are evaluated via \eqref{eq_tilde_phi_as_sum},
\eqref{eq_dtilde_phi_as_sum} essentially to machine precision.
\label{rem_phi_tmin}
\end{remark}

We evaluate $\tilde{\Phi}_n$ at all but the last four remaining roots
of $\psi_n$ in $[0,1)$ as follows.
Suppose that $n/2 < j < n$ is an integer,
and both 
$\tilde{\Phi}_n(t_j)$ and $\tilde{\Phi}_n'(t_j)$
have already been evaluated.

\paragraph{Step 2 (evaluation of
$\tilde{\Phi}_n(t_{j+1})$ and $\tilde{\Phi}_n'(t_{j+1})$,
given $\tilde{\Phi}_n(t_j)$ and $\tilde{\Phi}_n'(t_j)$).}
\begin{itemize}
\item use the recurrence relation 
\eqref{eq_num_dtilde_phi_3}, \eqref{eq_dphi_reck}
(see Theorem~\ref{lem_tilde_phi_ode} in Section~\ref{sec_weights})
to evaluate $\tilde{\Phi}_n^{(2)}(t_j), \dots, \tilde{\Phi}_n^{(M)}(t_j)$
(here $M=60$ in double precision calculations, and
$M=120$ in extended precision calculations).
\item evaluate $\tilde{\Phi}_n(t_{j+1})$ via Taylor's method
(see Section~\ref{sec_taylor}). The Taylor expansion of appropriate
order $M$ is used, i.e.
\begin{align}
\tilde{\Phi}_n(t_{j+1}) = 
\sum_{k=0}^{M} \frac{\tilde{\Phi}_n^{(k)}(t_j)}{k!} \cdot (t_{j+1}-t_j)^k
+
O\left( (t_{j+1}-t_j)^{M+1} \right)
\label{eq_phi_newton}
\end{align}
(compare to \eqref{eq_t_taylor} in Section~\ref{sec_evaluate_nodes}).
\item evaluate $\tilde{\Phi}'_n(t_{j+1})$ via Taylor's method.
The Taylor expansion of order $M-1$ is used, i.e.
\begin{align}
\tilde{\Phi}'_n(t_{j+1}) =
\sum_{k=0}^{M-1} \frac{\tilde{\Phi}_n^{(k+1)}(t_j)}{k!} \cdot (t_{j+1}-t_j)^k
+ O\left( (t_{j+1}-t_j)^M \right)
\label{eq_dphi_newton}
\end{align}
(compare to \eqref{eq_dpsi_taylor} in Section~\ref{sec_evaluate_nodes}).
In both \eqref{eq_phi_newton} and \eqref{eq_dphi_newton}, we set $M=60$
for double precision calculations and $M=120$ for 
extended precision calculations.
\end{itemize}
\begin{remark}
For each $j$, the cost of Step 2 
is $O(1)$ operations (i.e. does not depend on $n$). 
Also, it turns out that $\tilde{\Phi}_n(t_j)$ and $\tilde{\Phi}_n'(t_j)$
are evaluated via
\eqref{eq_phi_newton},
\eqref{eq_dphi_newton}
respectively,
essentially to machine precision
(compare to \eqref{eq_t_taylor}, \eqref{eq_dpsi_taylor}
in Section~\ref{sec_evaluate_nodes}).
 For a detailed discussion
of the accuracy and stability of this step, the reader
is referred to \cite{Glaser}.
\label{rem_algorithm2_step2}
\end{remark}

\paragraph{Step 3 (evaluation of
$\tilde{\Phi}_n(t_j)$ for $n-3 \leq j \leq n$).}
For $j=n-3, n-2, n-1, n$, we evaluate $\tilde{\Phi}_n(t_j)$
via the formula 
\begin{align}
\tilde{\Phi}_n(t_j) = 
\sum_{k = 0}^{2N} \alpha_k^{(n)} \cdot Q_k(t_j) =
\sum_{k = 0}^{2N} \beta_k^{(n)} \cdot Q_k(t_j) \cdot \sqrt{k+1/2}
\label{eq_tilde_phi_j_as_sum}
\end{align}
(as in \eqref{eq_tilde_phi_as_sum} in Step 1; 
see also
\eqref{eq_wj_as_sum} in the description of
Algorithm 1 above).
\begin{remark}
We compute $\tilde{\Phi}_n$ at the last four nodes via
\eqref{eq_tilde_phi_j_as_sum} rather than \eqref{eq_phi_newton},
since the accuracy of the latter deteriorates when 
$t_j$ is too close to $1$
(interestingly, the evaluation of $\psi_n(t_j)$ via
\eqref{eq_t_taylor} in Section~\ref{sec_evaluate_nodes}
for any $j=1,\dots,n$ does not have this unpleasant feature).
Since this approach works in practice, is cheap in terms of
the number of operations and eliminates the accuracy problem,
there was no need in a detailed analysis of the issue
(see, however, \cite{Glaser} for more details).
\label{rem_last_weights}
\end{remark}

\paragraph{Step 4 (evaluation of $\tilde{\Phi}_n(t_j)$ for
$1 \leq j \leq n/2$).}
Due to the combination of 
Theorem~\ref{thm_tilde_phi} in Section~\ref{sec_weights}
and \eqref{eq_legendre_fun_rec} in Section~\ref{sec_legendre},
the function 
$\tilde{\Phi}_n$ is symmetric about the origin. We use this observation
to evaluate $\tilde{\Phi}_n(t_j)$ via the formula
\begin{align}
\tilde{\Phi}_n(t_j) = (-1)^{n+1} \cdot \tilde{\Phi}_n(t_{n+1-j}),
\label{eq_tilde_phi_sym}
\end{align}
for every $j = 1, 2, \dots, n/2$.

\paragraph{Step 5 (evaluation of $W_1, \dots, W_n$).}
For every $j=1,\dots,n$,
we compute an approximation $\widehat{W}_j$ 
to $W_j$ from $\tilde{\Phi}_n(t_j)$ 
and $\psi_n'(t_j)$ via the formula
\begin{align}
\widehat{W}_j =
- 2 \cdot \frac{ \tilde{\Phi}_n(t_j) }{ \psi_n'(t_j) }
\label{eq_tilde_phi_wshort}
\end{align}
(see 
\eqref{eq_tilde_phi_w}
in Theorem~\ref{thm_tilde_phi}
in Section~\ref{sec_weights}).

\begin{remark}
Due to the combination of Remarks~\ref{rem_phi_tmin},
\ref{rem_algorithm2_step2}, \ref{rem_last_weights},
Algorithm 2 evaluates all $W_1,\dots,W_n$
essentially to machine precision.
This algorithms requires
$O(n)$ operations 
(compare to Remark~\ref{rem_algorithm1_cost}).
\label{rem_weights_cost}
\end{remark}
\begin{remark}
Algorithm 2 described in this section uses some of the quantities
evaluated by the procedures of Sections~\ref{sec_evaluate_beta},
\ref{sec_evaluate_nodes}. If $n<2c/\pi$, then 
the cost of obtaining these quantities
is $O\left(n + c \log(c)\right)$ operations; 
if $n>2c/\pi$, then the cost of obtaining these quantities
is $O(n)$ operations
(see Remarks~\ref{rem_evaluate_nodes}, 
\ref{rem_evaluate_nodes2} in Section~\ref{sec_evaluate_nodes}).
\label{rem_weights_cost2}
\end{remark}

%%%%%%%%%%%%%%%%%%%%%%%%%%%%%%%%%%%%%%%%%%%%%
%%%%%%%%%%%%%%%%%%%%%%%%%%%%%%%%%%%%%%%%%%%%%
\section{Numerical Results}
\label{sec_num_res}

In this section, we demonstrate
the performance of the quadrature rules
from Section~\ref{sec_quad}.
All the calculations were
implemented in FORTRAN (the Lahey 95 LINUX version),
and carried out in double precision.
Extended precision calculations
were used for comparison and verification
(in extended precision, the floating point numbers are
128 bits long, as opposed to 64 bits in double precision).

\paragraph{Experiment 1.}
\label{exp_0}

%%%%%%%%%%%%%%%%%%%%%
\begin{figure} [htbp]
\begin{center}
\includegraphics[width=11.5cm, bb=81   227   529   564, clip=true]
{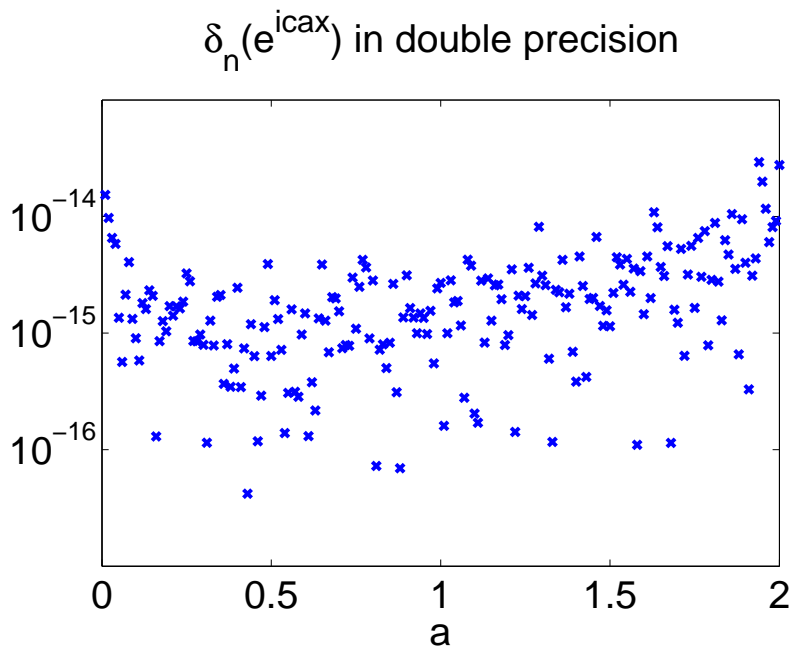}
\caption
{\it The quadrature error vs $|\lambda_n|$, with $c=1000$ and $n=682$. 
Here $\lambda_n = \mbox{\text{\rm{-.60352E-15}}}$.}
\label{fig:test253a}
\end{center}
\end{figure}
%%%%%%%%%%%%
%%%%%%%%%%%%%%%%%%%%%
\begin{figure} [htbp]
\begin{center}
\includegraphics[width=11.5cm, bb=81   227   529   564, clip=true]
{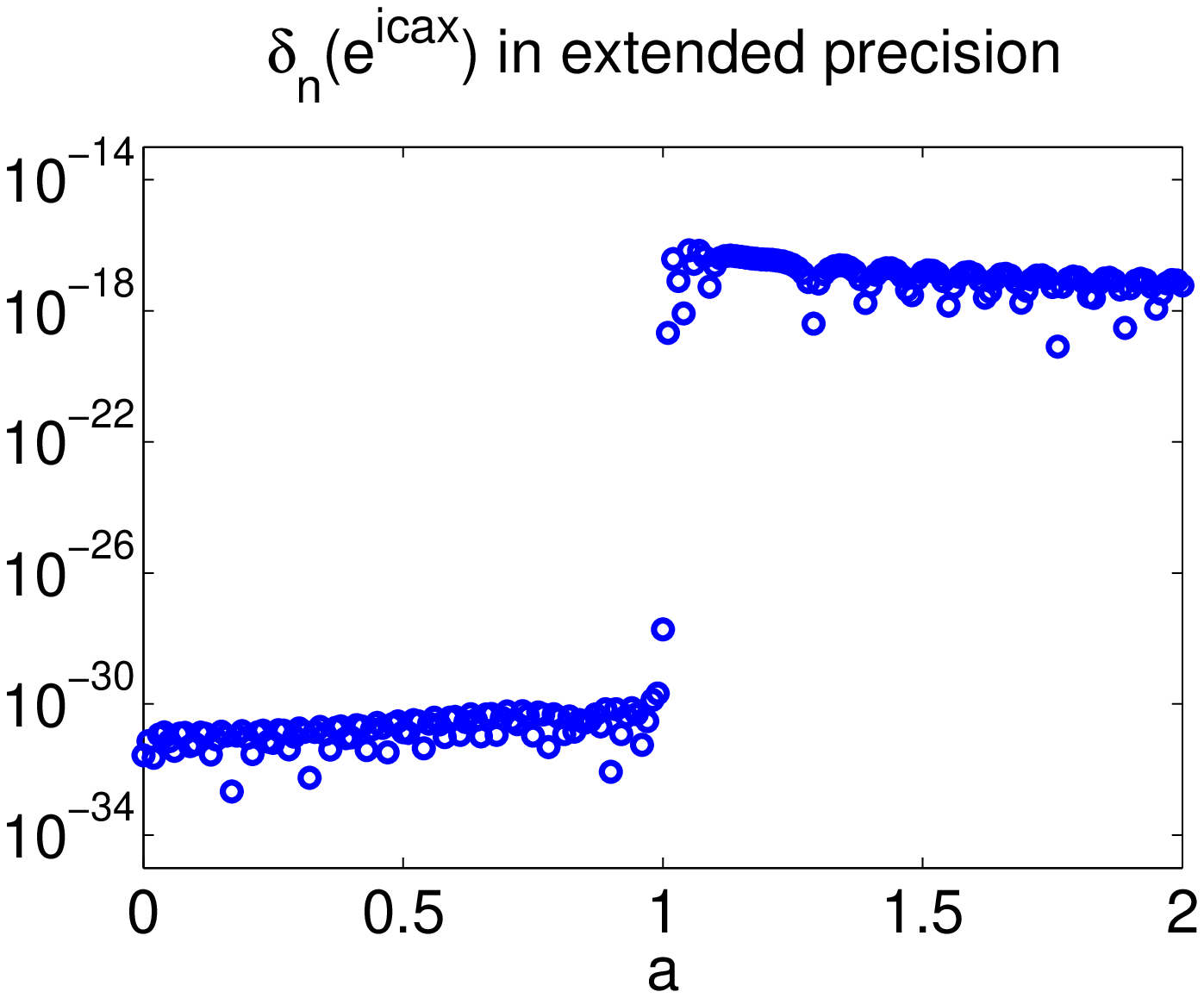}
\caption
{\it The quadrature error vs $|\lambda_n|$, with $c=1000$ and $n=682$. 
Here $\lambda_n = \mbox{\text{\rm{-.60352E-15}}}$.}
\label{fig:test253b}
\end{center}
\end{figure}
%%%%%%%%%%%%

%%%%%%%%%%%%%%%%%%%%%
\begin{figure} [htbp]
\begin{center}
\includegraphics[width=11.5cm, bb=81   227   529   564, clip=true]
{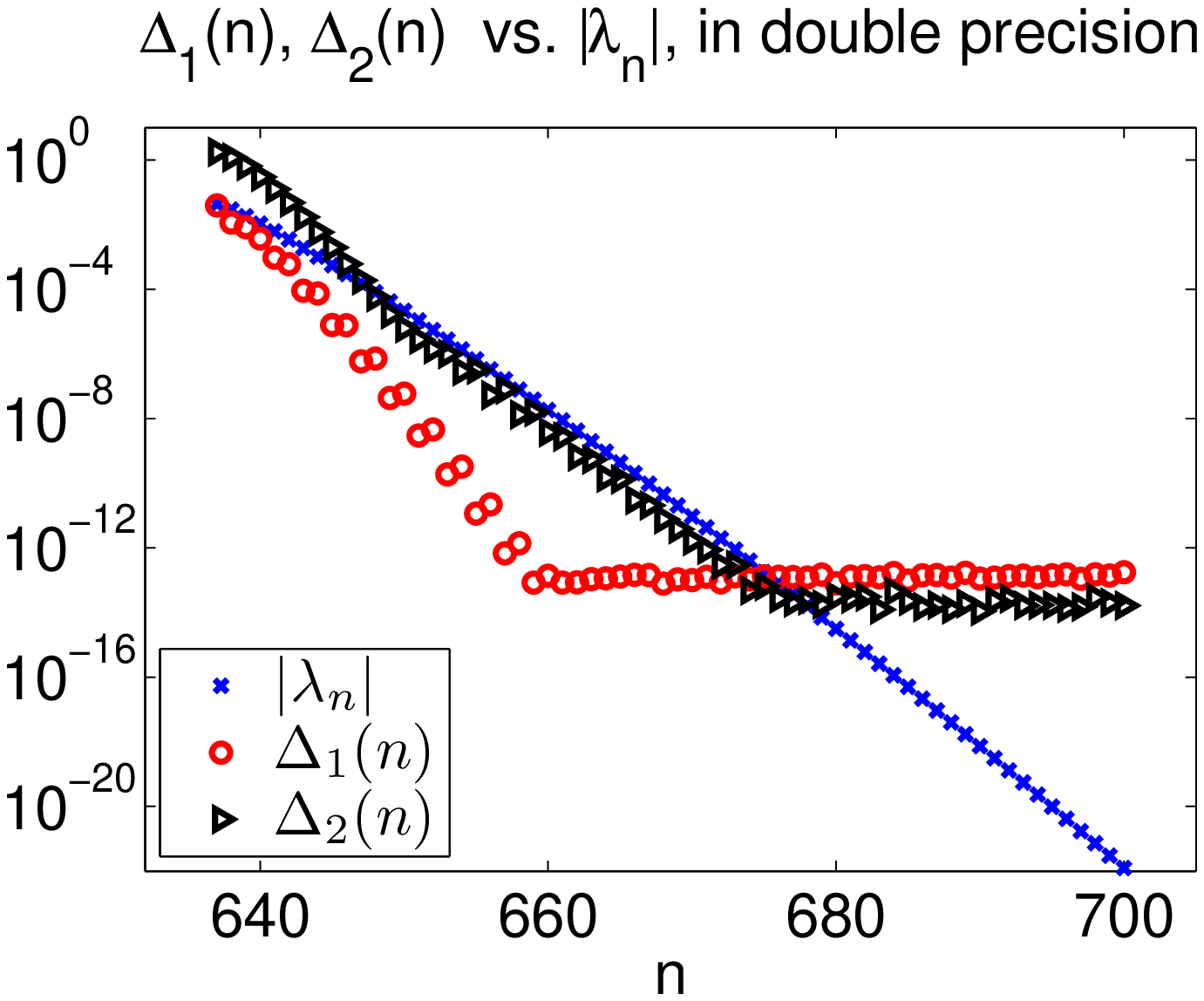}
\caption
{\it The maximal quadrature errors $\Delta_1(n), \Delta_2(n)$ 
vs $|\lambda_n|$, with $c=1000$. }
\label{fig:test254a}
\end{center}
\end{figure}
%%%%%%%%%%%%
%%%%%%%%%%%%%%%%%%%%%
\begin{figure} [htbp]
\begin{center}
\includegraphics[width=11.5cm, bb=81   227   529   564, clip=true]
{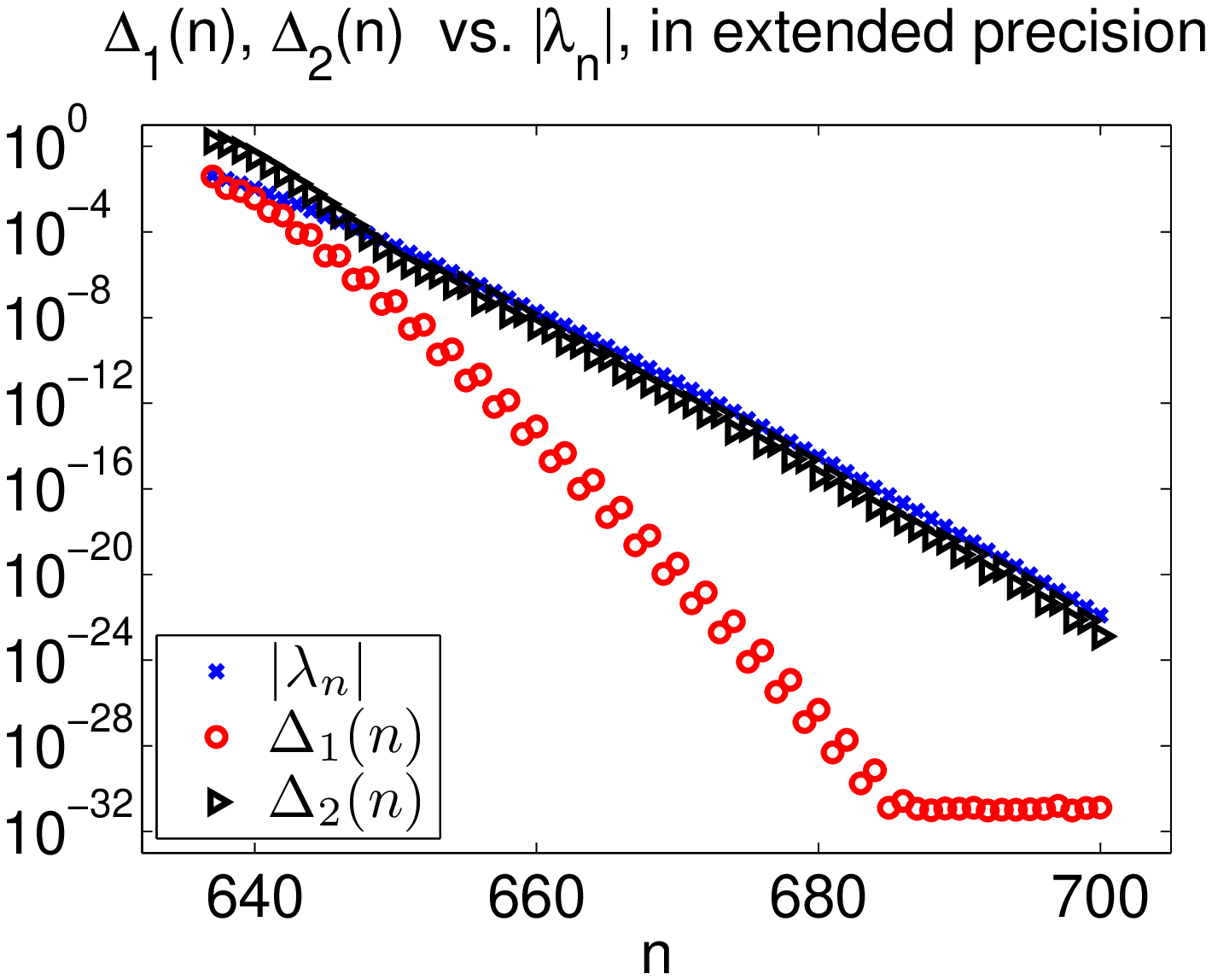}
\caption
{\it The maximal quadrature errors $\Delta_1(n), \Delta_2(n)$ 
vs $|\lambda_n|$, with $c=1000$. }
\label{fig:test254b}
\end{center}
\end{figure}
%%%%%%%%%%%%

Here we demonstrate the performance of the quadrature rule
$S_n$ (see \eqref{eq_quad_sn} in Section~\ref{sec_quad})
on exponential functions. We proceed as follows.
We choose, more or less arbitrarily, the band limit $c$
and the prolate index $n$. Next, we evaluate the quadrature
nodes $t_1, \dots, t_n$ and the quadrature weights $W_1,\dots,W_n$
via
the algorithms of Sections~\ref{sec_evaluate_nodes},
~\ref{sec_evaluate_weights}, respectively.
Also, we evaluate $|\lambda_n|$
via the algorithm in Section~\ref{sec_evaluate_lambda}.
Then, we choose a real number $0 \leq a \leq 2$, and evaluate the 
integral of $e^{icax}$ over $-1 \leq x \leq 1$ via
the formula
\begin{align}
\int_{-1}^1 e^{iacx} \; dx = 
\int_{-1}^1 \cos(acx) \; dx = \frac{2 \sin(ac)}{ac}.
\label{eq_exp0_int}
\end{align}
Also, we use $S_n$ to approximate \eqref{eq_exp0_int}
via the formula
\begin{align}
\int_{-1}^1 e^{iacx} \; dx \approx
\sum_{j=1}^n e^{icat_j} \cdot W_j
\label{eq_exp0_quad}
\end{align}
(see \eqref{eq_quad_quad} in Section~\ref{sec_quad}).
Finally, we evaluate the quadrature error $\delta_n(e^{iacx})$ via 
the formula
\begin{align}
\delta_n(e^{iacx}) =
\left|
\frac{2 \sin(ac)}{ac} - \sum_{j=1}^n e^{icat_j} \cdot W_j \right|
\label{eq_exp0_error}
\end{align}
(see \eqref{eq_quad_error_def} in Section~\ref{sec_quad}).

In Figure~\ref{fig:test253a}, we display the results
of this experiment. The band limit and the prolate index were
chosen to be, respectively, $c=1000$ and $n=682$. 
For this choice of parameters, $\lambda_n =$ -.60352E-15.
In this figure, we plot the quadrature error \eqref{eq_exp0_error}
as a function of the real parameter $a$, for $0 \leq a \leq 2$,
on the logarithmic scale.
The calculations are carried out in double precision.

We make the following observations
from Figure~\ref{fig:test253a}.
The quadrature error is essentially zero up to machine precision
$\varepsilon$,
for all real $0 \leq a \leq 2$. In other words,
for this choice of parameters,
the quadrature rule $S_n$ integrates the functions of the form
$f(x) = e^{icax}$ with $0 \leq a \leq 1$ exactly, for all 
practical purposes. It is perhaps surprising, however, 
that such functions are integrated exactly via $S_n$ even when 
$1 < a \leq 2$. In other words, the quadrature rule 
$S_n$ (corresponding to band limit $c$ and $|\lambda_n| \approx \varepsilon$)
integrates exactly the exponential functions with the band limit
up to $2c$.

To get a clearer picture, we repeat this experiment in extended precision. 
In Figure~\ref{fig:test253b}, we plot the
quadrature error \eqref{eq_exp0_error} as a function
of the real parameter $a$, for $0 \leq a \leq 2$,
on the logarithmic scale.
In other words, Figure~\ref{fig:test253b} is a version
of Figure~\ref{fig:test253a} in extended precision.

We make the following observations from Figure~\ref{fig:test253b}.
If $0 \leq a \leq 1$, then the quadrature rule $S_n$ integrates
the functions of the form $f(x) = e^{icax}$ up to the error
of order  $|\lambda_n|^2$
(in Figure~\ref{fig:test253a} we used double precision calculations
and thus did not have enough digits
to see this phenomenon). On the other hand, for $1 < a \leq 2$
the quadrature rule $S_n$ integrates such functions
up to the error roughly $|\lambda_n|$.
In other words, the quadrature rule $S_n$ (corresponding
to band limit $c$ and $|\lambda_n| \approx \varepsilon$)
integrates the functions of band limit up to $c$ up to $\varepsilon^2$
(rather than $\varepsilon$); on the other hand,
the functions of band limit 
between $c$ and $2c$ are integrated up to $\varepsilon$.

{\bf Explanation.} These observations admit the following (somewhat imprecise)
explanation
(see \cite{Report3}, \cite{Report3Arxiv} for more details).
 Suppose that $a \geq 0$ is a real number.
Due to \eqref{eq_prolate_integral} and
Theorem~\ref{thm_pswf_main} in Section~\ref{sec_pswf},
\begin{align}
e^{iacx} = \sum_{m=0}^{\infty} \lambda_m \psi_m(a) \psi_m(x),
\label{eq_exp0_exp}
\end{align}
for all real $-1 \leq x \leq 1$.
% (we note that while $f(x) = e^{iacx}$ is
% not a bandlimited function of $-1\leq x \leq 1$, 
% it does belong to $L^2\left[-1,1\right]$).
Moreover,
\begin{align}
\int_{-1}^1 e^{iacx} \; dx =
\frac{2\sin(ac)}{ac} = \sum_{m=0}^{\infty} \lambda_m^2 \psi_m(a) \psi_m(0).
\label{eq_exp0_sinc}
\end{align}
We combine \eqref{eq_exp0_error}, \eqref{eq_exp0_exp}, 
\eqref{eq_exp0_sinc} to obtain
\begin{align}
\frac{2 \sin(ac)}{ac} - \sum_{j=1}^n W_j \cdot e^{icat_j} = 
\sum_{m = 0}^{\infty} \lambda_m \psi_m(a) 
\left(\lambda_m \psi_m(0) - \sum_{j=1}^n W_j \psi_m(t_j) \right).
\label{eq_exp0_exp_error}
\end{align}
Obviously, the quadrature error $\delta_n(\psi_m)$ (see \eqref{eq_exp12_dif})
is zero for odd $m$. Also, $\delta_n(\psi_m)$ rapidly increases
as a function 
of even $0 \leq m < n$; moreover, 
 $\delta_n(\psi_m)$ 
is of order $|\lambda_n|$ when $m<n$ is an even integer  close to $n$
(see Conjectures~\ref{conj_deltan_psim},~\ref{conj_quad_error}
in Section~\ref{sec_quad_err_num} and Theorem~\ref{thm_quad_simple}
in Section~\ref{sec_quad_error}).
Therefore, roughly speaking, 
\begin{align}
\sum_{m = 0}^{n-1} \lambda_m \psi_m(a) 
\left(\lambda_m \psi_m(0) - \sum_{j=1}^n W_j \psi_m(t_j) \right) 
% = O\left( |\lambda_n|^2 \cdot \psi_{n-1}(a) \right).
\approx |\lambda_n|^2 \cdot \psi_{n-1}(a).
\label{eq_exp0_head}
\end{align}
On the other hand, due to the fast decay of $|\lambda_m|$
(see Theorems~\ref{thm_mu_spectrum}, \ref{thm_crude_inequality}
in Section~\ref{sec_pswf}), 
\begin{align}
\sum_{m = n}^{\infty} \lambda_m \psi_m(a) 
\left(\lambda_m \psi_m(0) - \sum_{j=1}^n W_j \psi_m(t_j) \right) 
% = O\left( |\lambda_n|^2 \cdot \psi_{n}(a) \right).
\approx |\lambda_n|^2.
\label{eq_exp0_tail}
\end{align}
Finally, the following approximate formula appears 
in \cite{Report3}, \cite{Report3Arxiv}, in a slightly different form:
suppose that $n>0$ is an integer, that $\chi_n > c^2$, and that
$0 \leq a \leq 2$ is a real number.
Then,
\begin{align}
|\psi_n(a)| = 
\begin{cases}
O\left(\sqrt{n}\right), & 0 \leq a \leq 1, \\
O\left( |\lambda_n|^{-1} \right), & 1 < a \leq 2.
\end{cases}
\label{eq_psin_a}
\end{align}
It follows from the combination of \eqref{eq_exp0_head},
\eqref{eq_exp0_tail}, \eqref{eq_psin_a} that
the quadrature error \eqref{eq_exp0_error} is expected 
to be of the order
$|\lambda_n|^2 \cdot \sqrt{n}$, if $0 \leq a \leq 1$. On the other hand,
the quadrature error \eqref{eq_exp0_error} is expected 
to be of the order
$|\lambda_n|$, if $1 < a \leq 2$.
Figures~\ref{fig:test253a}, \ref{fig:test253b},
\ref{fig:test254a}, \ref{fig:test254b} support
these somewhat vague conclusions.

We summarize this crude analysis, supported by the observations above,
in the following conjecture about the quadrature error
\eqref{eq_exp0_error} for $0 \leq a \leq 2$.
\begin{conjecture}
Suppose that $c>0$ and $a \geq 0$ are real numbers, 
and that $n>2c/\pi$ is an integer. 
Suppose also that $\delta_n(e^{icax})$ is defined
via \eqref{eq_quad_error_def} in Definition~\ref{def_quad}
in Section~\ref{sec_quad}. If $0 \leq a \leq 1$,
then
\begin{align}
\delta_n\left(e^{icax}\right) =
\left|
\int_{-1}^1 e^{icax} \; dx - \sum_{j=1}^n e^{icat_j} \cdot W_j 
\right| 
% = O\left( |\lambda_n|^2 \cdot \sqrt{n} \right),
\approx |\lambda_n|^2 \cdot \sqrt{n},
\label{eq_exp_conj}
\end{align}
where $\lambda_n$ is that of \eqref{eq_prolate_integral}
in Section~\ref{sec_pswf}. If, on the other hand,
$1 < a \leq 2$, then
\begin{align}
\delta_n\left(e^{icax}\right) =
\left|
\int_{-1}^1 e^{icax} \; dx - \sum_{j=1}^n e^{icat_j} \cdot W_j 
\right| 
% = O\left( |\lambda_n| \right).
\approx |\lambda_n|.
\label{eq_exp_conj2}
\end{align}
\label{conj_exp}
\end{conjecture}

We repeat the above experiment with various values of $n$, and 
plot the results 
in Figure~\ref{fig:test254a}. This figure also corresponds
to band limit $c=1000$. We plot the following three quantities
as functions of the prolate index $n$ that varies
between $637 \approx 2c/\pi$ and $700$.
First, we plot $|\lambda_n|$. Second, we plot 
the maximal quadrature error
$\Delta_1(n)$
defined via the formula
\begin{align}
\Delta_1(n) = \max_{0 \leq a \leq 1} \delta_n(e^{icax}) 
            = \max_{0 \leq a \leq 1} 
\left|
\frac{2 \sin(ac)}{ac} - \sum_{j=1}^n e^{icat^{(n)}_j} \cdot W^{(n)}_j \right|,
\label{eq_big_delta1}
\end{align}
where $t_1^{(n)},\dots,t_n^{(n)}$ and $W_1^{(n)},\dots,W_n^{(n)}$
are, respectively, the notes and weights of the quadrature rule $S_n$
(see \eqref{eq_quad_sn} in Section~\ref{sec_quad}). 
Finally, we plot the maximal quadrature error 
$\Delta_2(n)$ defined via the formula
\begin{align}
\Delta_2(n) = \max_{1 < a \leq 2} \delta_n(e^{icax}) 
            = \max_{1 < a \leq 2} 
\left|
\frac{2 \sin(ac)}{ac} - \sum_{j=1}^n e^{icat^{(n)}_j} \cdot W^{(n)}_j \right|.
\label{eq_big_delta2}
\end{align}
We observe that in \eqref{eq_big_delta1} the parameter $a$ varies between
$0$ and $1$, and in \eqref{eq_big_delta2} the parameter $a$
varies between $1$ and $2$. In other words,
$\Delta_1(n)$ is the maximal quadrature errors of $S_n$ for 
the exponential functions of band limits up to $c$, and
$\Delta_2(n)$ is the maximal quadrature error of $S_n$
for the exponential functions of band limit between $c$ and $2c$.

We make the following observations from Figure~\ref{fig:test254a}.
As long as $|\lambda_n|$ is less than roughly 
$10^{-7} \approx \sqrt{\varepsilon}$ 
(with $\varepsilon$ the machine precision), $\Delta_1(n)$ 
is roughly equal to $|\lambda_n|^2$. On the other hand,
$\Delta_1(n)$ is zero up to machine precision once
$|\lambda_n| > 10^{-7}$.
These observations are in agreement with Conjecture~\ref{conj_exp}
above.

We also observe that $\Delta_2(n)$ is roughly of order $|\lambda_n|$,
as long as $|\lambda_n| > \varepsilon$. On the other hand,
when $\lambda_n$ is zero to machine precision, so is $\Delta_2(n)$
(see Conjecture~\ref{conj_exp}).

We repeat this experiment in extended precision, and plot
the results in Figure~\ref{fig:test254b}. In other words,
Figure~\ref{fig:test254b} is a version of Figure~\ref{fig:test254a}
in extended precision. 
We observe the same phenomenon: $\Delta_1(n)$ is of order $|\lambda_n|^2$,
and $\Delta_2(n)$ is of order $|\lambda_n|$
(as long as we do not run out of digits to see it; if, for example,
 $|\lambda_n|$ is below the machine zero so are both $\Delta_1(n)$
and $\Delta_2(n)$).
In other words, the quadrature error of $S_n$
for exponential functions with band limit up to $c$ is
of order $|\lambda_n|^2$,
and the quadrature error of $S_n$
for exponential functions with band limit between $c$ and $2c$
is of order $|\lambda_n|$,
which supports Conjecture~\ref{conj_exp}.

%%%%%%%%%%%%%%%%%%%%%%%%%%%%%%%%%%%%%%%%%%%%%
%%%%%%%%%%%%%%%%%%%%%%%%%%%%%%%%%%%%%%%%%%%%%
\section{Numerical Illustration of Analysis in Section~\ref{sec_quad}}
\label{sec_num_ill}

In this section, we illustrate the analytical results
from Section~\ref{sec_quad} and 
the performance of the algorithms described in
Section~\ref{sec_num_algo}. All the calculations were
implemented in FORTRAN (the Lahey 95 LINUX version),
and carried out in double precision.
Extended precision calculations
were used for comparison and verification
(in extended precision, the floating point numbers are
128 bits long, as opposed to 64 bits in double precision).

%%%%%%%%%%%%%%%%%%%%%%%%%%%%%%
\subsection{Quadrature Error and its Relation to $|\lambda_n|$}
\label{sec_quad_err_num}
In this section, we describe several numerical experiments
that illustrate the quadrature error 
(see \eqref{eq_quad_sn}, \eqref{eq_quad_error_def}
in Section~\ref{sec_quad}) and its relation
to $|\lambda_n|$.

%%%%%%%%%%%%%%%%%%%
\begin{table}[htbp]
\begin{center}
\begin{tabular}{c|c|c|c}
$m$   &
$ \displaystyle \lambda_m \psi_m(0)$ &
$ \displaystyle \delta_n(\psi_m), \text{ double precision } $ &
$ \displaystyle \delta_n(\psi_m), \text{ extended precision }$
\\[2ex]
\hline
 0
& 0.70669E+00 & 0.44409E-15  & 0.33258E-26 \\
 2
& 0.49581E+00 & 0.16653E-15  & 0.22426E-25 \\
 4
& 0.42581E+00 & 0.13323E-14  & 0.26756E-23 \\
 6
& 0.38527E+00 & 0.21649E-14  & 0.19692E-21 \\
 8
& 0.35695E+00 & 0.22760E-14  & 0.91546E-20 \\
 10
& 0.33516E+00 & 0.16653E-14  & 0.29148E-18 \\
 12
& 0.31730E+00 & 0.23870E-14  & 0.88165E-17 \\
 14
& 0.30201E+00 & 0.24980E-14  & 0.21007E-15 \\
 16
& 0.28844E+00 & 0.11102E-14  & 0.35574E-14 \\
 18
& 0.27604E+00 & 0.59230E-13  & 0.57028E-13 \\
 20
& 0.26435E+00 & 0.83716E-12  & 0.83954E-12 \\
 22
& 0.25299E+00 & 0.89038E-11  & 0.89011E-11 \\
 24
& 0.24150E+00 & 0.76862E-10  & 0.76864E-10 \\
 26
& 0.22919E+00 & 0.65870E-09  & 0.65870E-09 \\
 28
& 0.21377E+00 & 0.45239E-08  & 0.45239E-08 \\
 30
& 0.18075E+00 & 0.19826E-07  & 0.19826E-07 \\
 32
& 0.10038E+00 & 0.68548E-07  & 0.68548E-07 \\
 34
& 0.27988E-01 & 0.33810E-06  & 0.33810E-06 \\
 36
& 0.49822E-02 & 0.27232E-05  & 0.27232E-05 \\
 38
& 0.70503E-03 & 0.22754E-04  & 0.22754E-04 \\
\end{tabular}
\end{center}
\caption{\it
Illustration of Theorem~\ref{thm_quad_simple}
with $c = 50$ and $n = 40$. For these parameters,
$\lambda_n = \mbox{\text{\rm{0.12915E-03}}}$.
%Here $|\lambda_n| =$ 0.12915E-03.
See Experiment 2. 
}
\label{t:test90}
\end{table}
%%%%%%%%%%%
\paragraph{Experiment 2.}
\label{sec_exp12}
Here we illustrate 
Theorem~~\ref{thm_quad_simple}
in Section~\ref{sec_quad_error}. 
We choose, more or less arbitrarily, band limit $c$ and
prolate index $n$. We evaluate $\chi_n$, $\lambda_n$
and the quadrature rule $S_n$ defined via \eqref{eq_quad_sn}
in Section~\ref{sec_quad} via the algorithms of
Sections~\ref{sec_evaluate_beta}, \ref{sec_evaluate_lambda},
\ref{sec_evaluate_nodes}, \ref{sec_evaluate_weights}, respectively.
Then, we choose an even integer $0 \leq m < n$, and
evaluate $\lambda_m$, $\psi_m(0)$, and $\psi_m(t_j)$
for all $j=1,\dots,n$, via the algorithms
of Sections~\ref{sec_evaluate_beta}, \ref{sec_evaluate_lambda}.
All the calculations are carried out in double precision.

We display the results of this experiment in
Table~\ref{t:test90}.
The data in this table correspond to $c=50$ and $n=40$.
Table~\ref{t:test90} has the following structure.
The first column contains the even integer $m$,
that varies between $0$ and $n-2$.
The second column contains $\lambda_m \psi_m(0)$
(we observe that
\begin{align}
\lambda_m \psi_m(0) = \int_{-1}^1 \psi_m(t) \; dt,
\label{eq_exp11_psim}
\end{align}
due to \eqref{eq_prolate_integral} in Section~\ref{sec_pswf}).
The third column contains the quadrature error 
\begin{align}
\delta_n(\psi_m) = \left|
\lambda_m \psi_m(0) - \sum_{j=1}^n \psi_m(t_j) \cdot W_j \right|
\label{eq_exp12_dif}
\end{align}
(see \eqref{eq_quad_error_def} in Section~\ref{sec_quad}),
computed in double precision. 

Then, we repeat all the calculations in extended precision;
the last column of Table~\ref{t:test90} contains $\delta_n(\psi_m)$
defined via \eqref{eq_exp12_dif}
(the same quantity as in the third column evaluated in extended precision).

We make the following observations from Table~\ref{t:test90}.
We note that $\lambda_m \psi_m(0)$ is always positive and
monotonically decreases with $m$.
We also note that 
$\delta_n(\psi_m)$ (evaluated in double precision)
is close to the machine accuracy for small $m$,
and grows up to $\approx 2 \cdot 10^{-5}$ for $m=38$.
Also, $\delta_n(\psi_m)$ 
is bounded by $|\lambda_n|$, for all values of $m$ in Table~\ref{t:test90}
(in this case, $|\lambda_n| =$ 0.12915E-03). Finally, 
$\delta_n(\psi_m)$ (evaluated in extended precision)
is a monotonically increasing function of even $0 \leq m < n$
(obviously, $\delta_n(\psi_m)=0$ for odd $m$).

We summarize these observations in the following conjecture.
We have not fully investigated the phenomenon described in
this conjecture; see, however, Theorem~\ref{thm_quad_simple}
in Section~\ref{sec_quad_error}, Conjecture~\ref{conj_quad_error} below,
Figure~\ref{fig:test92} and Table~\ref{t:test91}
(see also \cite{Report3}, \cite{Report3Arxiv} for additional
details and analysis).
\begin{conjecture}
Suppose that $c>1$ is a real number, that $n > 2c/\pi$ is an integer,
and that the quadrature rule $S_n$ is defined via \eqref{eq_quad_sn}
in Section~\ref{sec_quad}.
Then, the quadrature error $\delta_n(\psi_m)$ defined via
\eqref{eq_exp12_dif} above
is a monotonically increasing function of even $0 \leq m < n$.
Moreover, in double precision calculations
$\delta_n(\psi_m)$ is zero up to machine precision
for all $0 \leq m < 2c/\pi$.
\label{conj_deltan_psim}
\end{conjecture}

In \eqref{eq_quad_simple_thm} in Theorem~\ref{thm_quad_simple}, 
we provide an upper bound on $\delta_n(\psi_m)$.
This bound does not depend on $m$; more specifically,
for every $m = 0,\dots,n-1$,
\begin{align}
\delta_n(\psi_m) = \left|
\lambda_m \psi_m(0) - \sum_{j=1}^n \psi_m(t_j) \cdot W_j
\right| 
\leq |\lambda_n| \cdot
\left(
24 \cdot \log\left( \frac{1}{|\lambda_n|} \right) +
6 \cdot \chi_n
\right).
\label{eq_exp12_simple}
\end{align}
On the other hand, according to Table~\ref{t:test90}
the quadrature error $\delta_n(\psi_m)$ is bounded by
$|\lambda_n|$ alone, for all even $0 \leq m < n$
(obviously, $\delta_n(\psi_m)=0$ for all odd $m$).

%%%%%%%%%%%%%%%%%%%%%
\begin{figure} [htbp]
\begin{center}
\includegraphics[width=11.5cm, bb=81   227   529   564, clip=true]
{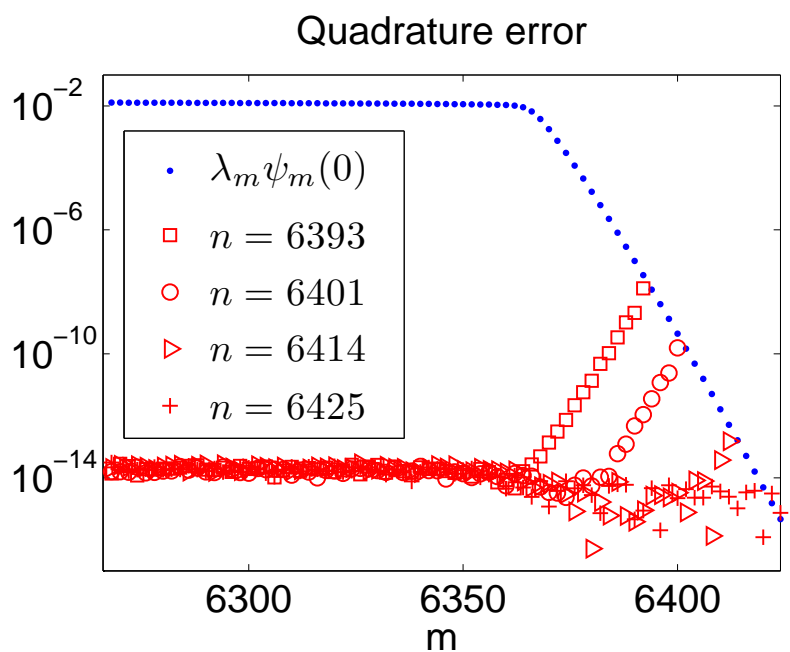}
\caption
{\it
The quadrature error 
% $\delta_n(\psi_m)$ (see \eqref{eq_quad_error_def})
$\displaystyle 
\delta_n(\psi_m) = \left| \int_{-1}^1 \psi_m(t) \; dt - 
\sum_{j=1}^n \psi_m(t_j) \cdot W_j \right|$
as a function of even $m < n$, for four different values
of $n$ and $c = 10000$, vs. 
$\displaystyle \lambda_m \psi_m(0)$. See Experiment 2.
}
\label{fig:test92}
\end{center}
\end{figure}
%%%%%%%%%%%%
In 
Figure~\ref{fig:test92}, we display the results of the same experiment
with different choice of parameters $c$ and $n$.
Namely, we choose $c=10000$ and plot $\lambda_m \psi_m(0)$
as a function of even $0 \leq m < 6425$, on the logarithmic scale
(solid line). In addition, we plot 
the quadrature error $\delta_n(\psi_m)$
as a function of $m$, for four
different values of $n$: $n=6393$ (dashed line),
$n=6401$ (circles), $n=6414$ (triangles),
and $n=6425$ (pluses). The corresponding values of $|\lambda_n|$
are displayed in Table~\ref{t:test92}.
%%%%%%%%%%%%%%%%%%%
\begin{table}[htbp]
\begin{center}
\begin{tabular}{c|c|c|c|c}
$n$ & 6393 & 6401 & 6414 & 6425 
\\
\hline
$|\lambda_n|$ & 0.43299E-07 & 0.54119E-09 & 0.33602E-12 & 0.52616E-15
\end{tabular}
\end{center}
\caption{\it
Values of $|\lambda_n|$ for $c=10000$ and different choices of $n$.
}
\label{t:test92}
\end{table}

We make the following observations from Figure~\ref{fig:test92}.
First, the quantities $\lambda_m \psi_m(0)$ are of 
the same order of magnitude for all $m < 2c/\pi$, 
and decay rapidly with $m$
for $m > 2c/\pi$.
Also, for each value of $n$,
the quadrature error $\delta_n(\psi_m)$ 
is essentially zero for all $m < 2c/\pi$
and increases rapidly with $m$ for $m > 2c/\pi$.
Nevertheless, $\delta_n(\psi_m)$
is always bounded from above by $|\lambda_n|$, for each $n$
and each $m<n$.
See also Tables~\ref{t:test90},~\ref{t:test91} and
Conjecture~\ref{conj_quad_error} below.

%%%%%%%%%%%%%%%%%%%%%%%%%%%%%%
%%%%%%%%%%%%%%%%%%%
\begin{table}[htbp]
\begin{center}
\begin{tabular}{c|c|c|c|c|c}
$c$   &
$n$   &
$m$   &
$\lambda_m \psi_m(0)$ &
$ \displaystyle \int_{-1}^1 \psi_m(t) \; dt 
  - \sum_{j=1}^n \psi_m(t_j) \cdot W_j$ &
$|\lambda_n|$ \\[1ex]
\hline
250 &  179  &  178
& 0.28699E-07 & -.52496E-08 & 0.18854E-07  \\
250 &  184  &  182
& 0.68573E-09 & -.38341E-10 & 0.16130E-09  \\
250 &  188  &  186
& 0.14108E-10 & -.68758E-12 & 0.30500E-11  \\
\hline
500 &  339  &  338
& 0.52368E-07 & -.13473E-07 & 0.40938E-07  \\
500 &  345  &  344
& 0.37412E-09 & -.86136E-10 & 0.27418E-09  \\
500 &  350  &  348
& 0.12148E-10 & -.99816E-12 & 0.35537E-11  \\
\hline
1000 & 659  &  658
& 0.42709E-07 & -.14354E-07 & 0.38241E-07  \\
1000 & 665  &  664
& 0.51665E-09 & -.15924E-09 & 0.43991E-09  \\
1000 &  671  &  670
& 0.52494E-11 & -.15024E-11 & 0.42815E-11  \\ 
\hline
2000 & 1297  &  1296
& 0.41418E-07 & -.17547E-07 & 0.41740E-07  \\
2000 &  1304  &  1302
& 0.77185E-09 & -.15036E-09 & 0.37721E-09  \\
2000 &  1311  &  1310
& 0.31078E-11 & -.11386E-11 & 0.28754E-11  \\
\hline
4000 &  2572  &  2570
& 0.54840E-07 & -.15493E-07 & 0.33682E-07  \\
4000 &  2579  &  2578
& 0.43032E-09 & -.20771E-09 & 0.46141E-09  \\
4000 &  2587  &  2586
& 0.28193E-11 & -.12805E-11 & 0.29164E-11  \\
\hline
8000 &  5119  &  5118
& 0.43268E-07 & -.26751E-07 & 0.52899E-07  \\
8000 &  5128  &  5126
& 0.50230E-09 & -.16395E-09 & 0.33442E-09  \\
8000 &  5136  &  5134
& 0.50508E-11 & -.15448E-11 & 0.32132E-11  \\
\hline
16000 &  10213  &  10212
& 0.42725E-07 & -.30880E-07 & 0.56568E-07  \\
16000 &  10222  &  10220
& 0.69663E-09 & -.28201E-09 & 0.52821E-09  \\
16000 &  10231  &  10230
& 0.34472E-11 & -.22162E-11 & 0.42902E-11  \\
\end{tabular}
\end{center}
\caption{\it
Relation between the quadrature error and $|\lambda_n|$.
See Experiment 2.
}
\label{t:test91}
\end{table}
%%%%%%%%%%%

We repeat the experiment above
with several other values of band limit $c$ and prolate index $n$.
The results are displayed in Table~\ref{t:test91}.
This table has the following structure.
The first and second column contain, respectively,
the band limit $c$ and the prolate index $n$.
The third column contains the even integer $0 \leq m < n$
(the values of $m$ were chosen to be close to $n$).
The fourth column contains $\lambda_m \psi_m(0)$.
The fifth column contains 
\eqref{eq_exp12_dif}. The last column contains $|\lambda_n|$.

We make the following observations from Table~\ref{t:test91}.
First, for each of the seven
values of $c$, the three indices $n$ were chosen in such
a way that $|\lambda_n|$ is between
$10^{-12}$ and $10^{-7}$. The values of the band limit 
$c$ vary between $250$
(the first three rows) and $16000$ (the last three rows).
For each $n$, the value of $m$ is chosen to be the largest
even integer below $n$. This choice of $m$ yields
the smallest $\lambda_m \psi_m(0)$ and the largest quadrature
error $\delta_n(\psi_m)$ among all $m < n$ 
(see also Table~\ref{t:test90} and Figure~\ref{fig:test92}). 
Obviously, for this choice of $m$ the eigenvalues $\lambda_m$ and
$\lambda_n$ are roughly of the same order of magnitude.
We also observe that for all the values of $c,n,m$, 
the quadrature error $\delta_n(\psi_m)$
is bounded from above by $|\lambda_n|$ (and is roughly equal
to $|\lambda_n|/2$).
In other words, the upper bound on $\delta_n(\psi_m)$ provided
by Theorem~\ref{thm_quad_simple} (see \eqref{eq_exp12_simple})
is somewhat overcautious.

We summarize these observations in the following conjecture.
\begin{conjecture}
Suppose that $c>0$ is a positive real number, and that $n>2c/\pi$ is
an integer. Suppose also that $0 \leq m < n$ is an integer.
Suppose furthermore that $\delta_n(\psi_m)$ is defined
via \eqref{eq_quad_error_def} in Definition~\ref{def_quad}
in Section~\ref{sec_quad}.
Then,
\begin{align}
\delta_n(\psi_m) = 
\left| \int_{-1}^1 \psi_m(s) \; ds - 
\sum_{j=1}^n \psi_m(t_j) \cdot W_j \right| \leq
|\lambda_n|,
\label{eq_quad_error_conj}
\end{align}
where $\lambda_n$ is that of \eqref{eq_prolate_integral}
in Section~\ref{sec_pswf}.
\label{conj_quad_error}
\end{conjecture}
\begin{remark}
The inequality \eqref{eq_quad_error_conj} in
Conjecture~\ref{conj_quad_error} is stronger than
the inequality \eqref{eq_quad_simple_thm} 
in Theorem~\ref{thm_quad_simple}. On the other hand,
as opposed to Theorem~\ref{thm_quad_simple},
Conjecture~\ref{conj_quad_error} has been only supported
by numerical evidence.
\label{rem_conj}
\end{remark}

\paragraph{Experiment 3.}
\label{sec_exp14}
Here we illustrate
Theorems~\ref{thm_quad_eps_large},~\ref{thm_quad_eps_simple}
in Section~\ref{sec_main_result}.
We proceed as follows. We choose, more or less arbitrarily,
the band limit $c > 0$ and the accuracy parameter $\varepsilon > 0$.
Then, we use the algorithm of Section~\ref{sec_evaluate_lambda}
to find the minimal integer $m$ 
such that $|\lambda_m| < \varepsilon$. 
In other words, we define the integer $n_1(\varepsilon)$
via the formula
\begin{align}
n_1(\varepsilon) = 
\min\left\{ m \geq 0 \; : \; |\lambda_m| < \varepsilon \right\}.
\label{eq_exp14_n1}
\end{align}
Also, we find the minimal integer such that the
right-hand side of \eqref{eq_quad_simple_thm} in
Theorem~\ref{thm_quad_simple} in Section~\ref{sec_quad_error}
is less that $\varepsilon$.
In other words, we define the integer $n_2(\varepsilon)$ via the formula
\begin{align}
n_2(\varepsilon) = 
\min\left\{ m \geq 0 \; : \; 
|\lambda_m| \cdot
\left(
24 \cdot \log\left( \frac{1}{|\lambda_m|} \right) +
6 \cdot \chi_m
\right)
< \varepsilon \right\}.
\label{eq_exp14_n2}
\end{align}
Next, we evaluate the integer $n_3(\varepsilon)$ via the formula
\eqref{eq_quad_eps_large_nu} in Theorem~\ref{thm_quad_eps_large}.
In other words,
\begin{align}
n_3(\varepsilon) = \text{floor}\left(
\frac{2c}{\pi} + \frac{\alpha(\varepsilon)}{2\pi} \cdot 
\log\left(
\frac{16ec}{\alpha(\varepsilon)},
\right)
\right)
\label{eq_exp14_n3}
\end{align}
where $\alpha(\varepsilon)$ is defined via
\eqref{eq_quad_eps_large_alpha} in Theorem~\ref{thm_quad_eps_large}.
Finally, we evaluate the integer $n_4(\varepsilon)$ via
the right-hand side of \eqref{eq_quad_eps_simple_n}
in Theorem~\ref{thm_quad_eps_simple}. In other words,
\begin{align}
n_4(\varepsilon) = \text{floor}\left(
\frac{2c}{\pi} +
\left(10 + \frac{3}{2} \cdot \log(c) + 
   \frac{1}{2} \cdot \log\frac{1}{\varepsilon}
\right) \cdot \log\left( \frac{c}{2} \right) \right).
\label{eq_exp14_n4}
\end{align}
In both \eqref{eq_exp14_n3} and \eqref{eq_exp14_n4}, $\text{floor}(a)$
denotes the integer part of a real number $a$.

%%%%%%%%%%%%%%%%%%%
\begin{table}[htbp]
\begin{center}
\begin{tabular}{c|c|c|c|c|c|c|c}
$c$     &
$\varepsilon$   &
$n_1(\varepsilon)$ &
$n_2(\varepsilon)$ &
$n_3(\varepsilon)$ &
$n_4(\varepsilon)$ &
$|\lambda_{n_1(\varepsilon)}|$ &
$|\lambda_{n_2(\varepsilon)}|$ \\[1ex]
\hline
   250 & $10^{-10}$ &   184&   198&   277&   303 & 0.60576E-10 & 0.86791E-16\\
   250 & $10^{-25}$ &   216&   227&   326&   386 & 0.31798E-25 & 0.14863E-30\\
   250 & $10^{-50}$ &   260&   270&   393&   525 & 0.28910E-50 & 0.75155E-56\\
\hline
   500 & $10^{-10}$ &   346&   362&   460&   488 & 0.49076E-10 & 0.60092E-16\\
   500 & $10^{-25}$ &   382&   397&   520&   583 & 0.54529E-25 & 0.19622E-31\\
   500 & $10^{-50}$ &   433&   446&   607&   742 & 0.82391E-50 & 0.38217E-56\\
\hline
  1000 & $10^{-10}$ &   666&   687&   803&   834 & 0.95582E-10 & 0.92947E-17\\
  1000 & $10^{-25}$ &   707&   725&   875&   942 & 0.97844E-25 & 0.14241E-31\\
  1000 & $10^{-50}$ &   767&   783&   981&  1120 & 0.39772E-50 & 0.56698E-57\\
\hline
  2000 & $10^{-10}$ &  1305&  1330&  1467&  1500 & 0.95177E-10 & 0.25349E-17\\
  2000 & $10^{-25}$ &  1351&  1373&  1550&  1619 & 0.86694E-25 & 0.27321E-32\\
  2000 & $10^{-50}$ &  1418&  1438&  1675&  1818 & 0.88841E-50 & 0.22795E-57\\
\hline
  4000 & $10^{-10}$ &  2581&  2610&  2768&  2804 & 0.70386E-10 & 0.64396E-18\\
  4000 & $10^{-25}$ &  2632&  2658&  2862&  2935 & 0.57213E-25 & 0.53827E-33\\
  4000 & $10^{-50}$ &  2707&  2730&  3007&  3154 & 0.56712E-50 & 0.88819E-58\\
\hline
  8000 & $10^{-10}$ &  5130&  5163&  5344&  5383 & 0.59447E-10 & 0.22821E-18\\
  8000 & $10^{-25}$ &  5185&  5216&  5450&  5526 & 0.87242E-25 & 0.16237E-33\\
  8000 & $10^{-50}$ &  5268&  5296&  5614&  5765 & 0.95784E-50 & 0.23927E-58\\
\hline
 16000 & $10^{-10}$ & 10225& 10264& 10468& 10509 & 0.63183E-10 & 0.37516E-19\\
 16000 & $10^{-25}$ & 10285& 10321& 10585& 10664 & 0.85910E-25 & 0.41416E-34\\
 16000 & $10^{-50}$ & 10377& 10409& 10769& 10923 & 0.51912E-50 & 0.56250E-59\\
\hline
 32000 & $10^{-10}$ & 20413& 20457& 20686& 20730 & 0.62113E-10 & 0.12818E-19\\
 32000 & $10^{-25}$ & 20478& 20519& 20815& 20897 & 0.78699E-25 & 0.12197E-34\\
 32000 & $10^{-50}$ & 20577& 20615& 21018& 21176 & 0.96802E-50 & 0.15816E-59\\
\hline
 64000 & $10^{-10}$ & 40786& 40837& 41092& 41139 & 0.89344E-10 & 0.28169E-20\\
 64000 & $10^{-25}$ & 40857& 40903& 41232& 41318 & 0.66605E-25 & 0.39212E-35\\
 64000 & $10^{-50}$ & 40964& 41008& 41454& 41616 & 0.85451E-50 & 0.28036E-60\\
\hline
$10^6$ & $10^{-10}$ &636669&636747&637115&637174 & 0.79326E-10 & 0.13385E-22\\
$10^6$ & $10^{-25}$ &636759&636832&637301&637400 & 0.77413E-25 & 0.15758E-37\\
$10^6$ & $10^{-50}$ &636899&636968&637600&637778 & 0.69235E-50 & 0.15801E-62\\
\end{tabular}
\end{center}
\caption{\it
Illustration of
Theorems~\ref{thm_quad_eps_large},~\ref{thm_quad_eps_simple}.
See Experiment 3.
}
\label{t:test178}
\end{table}
%%%%%%%%%%%
We display the results of this experiment in Table~\ref{t:test178}.
This table has the following structure.
The first column contains the band limit $c$.
The second column contains the accuracy parameter $\varepsilon$.
The third column contains $n_1(\varepsilon)$ defined via
\eqref{eq_exp14_n1}. 
The fourth column contains $n_2(\varepsilon)$ defined via
\eqref{eq_exp14_n2}. 
The fifth column contains $n_3(\varepsilon)$ defined via
\eqref{eq_exp14_n3}. 
The sixth column contains $n_4(\varepsilon)$ defined via
\eqref{eq_exp14_n4}. 
The seventh column contains $|\lambda_{n_1(\varepsilon)}|$.
The last column contains $|\lambda_{n_2(\varepsilon)}|$.

Suppose that $c>0$ is a band limit, and $n>0$ is an integer.
We define the real number $Q(c,n)$ via the formula
\begin{align}
Q(c,n) = \max\left\{
\delta_n(\psi_m) = 
\left| 
\int_{-1}^1 \psi_m(t) \; dt - \sum_{j=1}^n \psi_m(t_j) \cdot W_j
\right| \; : \; 0 \leq m \leq n-1
\right\},
\label{eq_exp14_qnc}
\end{align}
where $t_1,\dots,t_n$ and $W_1,\dots,W_n$ are
defined, respectively, via \eqref{eq_quad_t}, \eqref{eq_quad_w} in
Definition~\ref{def_quad} in Section~\ref{sec_quad}.
In other words, $Q(c,n)$ is the maximal error
to which the quadrature rule $S_n$ defined via \eqref{eq_quad_sn}
integrates the first $n$ PSWFs.

We make the following observations from Table~\ref{t:test178}.
We observe that $Q(c,n_1(\varepsilon)) < \varepsilon$,
due to the combination of Conjecture~\ref{conj_quad_error} 
in Section~\ref{sec_exp12} and
\eqref{eq_exp14_n1}, \eqref{eq_exp14_qnc}. In other words,
numerical evidence suggests that the quadrature rule
$S_{n_1(\varepsilon)}$ integrates
the first $n_1(\varepsilon)$ PSWFs up to an error less than $\varepsilon$
(see Remark~\ref{rem_conj}).
On the other hand, we combine Theorem~\ref{thm_quad_simple}
in Section~\ref{sec_quad_error}
with \eqref{eq_exp14_n2}, \eqref{eq_exp14_qnc}, to conclude that
the quadrature rule $S_{n_2(\varepsilon)}$ has been
\emph{rigorously proven} to integrate the first $n_2(\varepsilon)$ PSWFs
up to an error less than $\varepsilon$. In both Theorem~\ref{thm_quad_simple}
and Conjecture~\ref{conj_quad_error}, we establish upper bounds
on $Q(c,n)$ in terms of $|\lambda_n|$. The ratio of 
$|\lambda_{n_1(\varepsilon)}|$ to $|\lambda_{n_2(\varepsilon)}|$ is quite
large: from about $10^6$ for $c=250$ and $\varepsilon = 10^{-10}, 10^{-25},
10^{-50}$ (see the first three rows in Table~\ref{t:test178}),
to about $10^{10}$ for $c=64000$
and $\varepsilon = 10^{-10}, 10^{-25}, 10^{-50}$,
to about $5 \cdot 10^{12}$ for $c=10^6$ 
and $\varepsilon = 10^{-10}, 10^{-25}, 10^{-50}$,
(see the last six rows in Table~\ref{t:test178}).
On the other hand, the difference between $n_2(\varepsilon)$ and
$n_1(\varepsilon)$ is fairly small; for example, for $\varepsilon=10^{-50}$,
this difference varies from 10 for $c=250$
to $23$ for $c=4000$, to merely $44$ for $c=64000$
and $69$ for as large $c$ as $c=10^6$.

As opposed to $n_1(\varepsilon)$ and $n_2(\varepsilon)$, 
the integer $n_3(\varepsilon)$
is computed via the explicit formula \eqref{eq_exp14_n3}
that depends only on $c$ and $\varepsilon$
(rather than on $|\lambda_n|$ and $\chi_n$, that need to be
evaluated numerically); this formula appears in
Theorem~\ref{thm_quad_eps_large}.
The convenience of \eqref{eq_exp14_n3} vs. \eqref{eq_exp14_n1},
\eqref{eq_exp14_n2} comes at a price: for example, 
for $\varepsilon = 10^{-50}$, the difference between $n_3(\varepsilon)$
and $n_2(\varepsilon)$ is equal to 123 for $c=250$, to 446 for $c=64000$,
and to 632 for $c=10^6$.
However, the difference 
$n_3(\varepsilon)-n_2(\varepsilon)$ is rather small compared to $c$:
for example, for $\varepsilon=10^{-50}$, this difference is
roughly $4 \cdot \left(\log(c)\right)^2$, for all values
of $c$ in Table~\ref{t:test178}.

Furthermore, we observe that 
$n_4(\varepsilon)$ is computed via the explicit
formula \eqref{eq_exp14_n4} that depends only on $c$ and $\varepsilon$.
This formula can be viewed as a simplified version of \eqref{eq_exp14_n3}
(see Theorems~\ref{thm_quad_eps_large}, \ref{thm_quad_eps_simple});
in particular,
$n_4(\varepsilon)$ is greater than $n_3(\varepsilon)$, for
all $c$ and $\varepsilon$.

We summarize these observations as follows. 
Suppose that the band limit $c>0$ and 
the accuracy parameter $\varepsilon>0$ are given.
According to Theorem~\ref{thm_quad_eps_large}, 
for any $n \geq n_3(\varepsilon)$
the quadrature
rule $S_n$ defined via \eqref{eq_quad_sn} in Section~\ref{sec_quad}
is {\it guaranteed} to integrate the first $n$ PSWFs to
precision $\varepsilon$
(see Definition~\ref{def_quad_n_eps} in Section~\ref{sec_outline_pswf}).
On the other hand, numerical evidence  (see Experiment 2) suggests that 
the choice $n \geq n_3(\varepsilon)$ is overly cautious for this purpose;
more specifically, $S_n$ integrates the first $n$ PSWFs
to precision $\varepsilon$
for every $n$ between $n_1(\varepsilon)$ and
$n_3(\varepsilon)$ as well. 
In this experiment, we 
observed that
the difference between the "theoretical" bound
$n_3(\varepsilon)$ and "empirical" bound $n_1(\varepsilon)$
is of order $\left(\log(c)\right)^2$, and, in particular,
is relatively small compared to 
both $n_1(\varepsilon)$ and $n_3(\varepsilon)$
(which are of order $c$).

Finally, we observe that
\begin{align}
n_1(\varepsilon) < \frac{2c}{\pi} + 
\frac{2}{\pi^2} \cdot \left( \log c \right) \cdot \log \frac{1}{\varepsilon},
\label{eq_n1_ineq}
\end{align}
for all the values of $c$ and $\varepsilon$
in Table~\ref{t:test178}.
Combined with some additional numerical experiments
by the authors, this observation leads to the following conjecture
(see also Theorem~\ref{thm_crude_inequality} in Section~\ref{sec_pswf}
for a rigorously proven
and more precise statement).
\begin{conjecture}
Suppose that $c>1$ and $0<\varepsilon<1$ are real numbers. Suppose also
that $n > 0$ is an integer, and that
\begin{align}
n > \frac{2c}{\pi} + 10 + 
\frac{2}{\pi^2} \cdot \left(\log c\right)  \cdot \log\frac{1}{\varepsilon}.
\label{eq_conj_n1_a}
\end{align}
Then,
\begin{align}
|\lambda_n| < \varepsilon,
\label{eq_conj_n1}
\end{align}
where $\lambda_n$ is that of \eqref{eq_prolate_integral}
in Section~\ref{sec_pswf}.
\label{conj_n1}
\end{conjecture}

%%%%%%%%%%%%%%%%%%%%%%%%%%%%%%
\subsection{Quadrature Weights}
\label{sec_quad_w_num}
In this section, we illustrate the results of Section~\ref{sec_weights}
and the algorithms of Section~\ref{sec_evaluate_weights}.
%%%%%%%%%%%%%%%%%%%
\begin{table}[htbp]
\begin{center}
\begin{tabular}{c|c|c|c}
$j$     &
$\widehat{W}_j$   &
$\widehat{W_j} - \widetilde{W}_j $ &
$\widehat{W}_j - \frac{\widehat{W}_{(n+1)/2} \left( \psi_n'(0) \right)^2}
      {\left( \psi_n'(t_j) \right)^2 \cdot \left(1 - t_j^2\right) }$ \\[1ex]
\hline
 1
&  0.7602931556894E-02 & 0.00000E+00 & -.55796E-11  \\
 2
&  0.1716167229714E-01 & 0.00000E+00 & -.55504E-10  \\
 3
&  0.2563684665002E-01 & 0.00000E+00 & -.21825E-12  \\
 4
&  0.3278512460580E-01 & 0.00000E+00 & -.11959E-09  \\
 5
&  0.3863462966166E-01 & 0.16653E-15 & 0.82238E-11  \\
 6
&  0.4334940472363E-01 & 0.22204E-15 & -.16247E-09  \\
 7
&  0.4713107235981E-01 & 0.22204E-15 & 0.11270E-10  \\
 8
&  0.5016785516291E-01 & 0.19429E-15 & -.18720E-09  \\
 9
&  0.5261660773966E-01 & 0.26368E-15 & 0.10495E-10  \\
 10
&  0.5460119701692E-01 & 0.29837E-15 & -.20097E-09  \\
 11
&  0.5621699326080E-01 & 0.17347E-15 & 0.81464E-11  \\
 12
&  0.5753664411864E-01 & 0.12490E-15 & -.20866E-09  \\
 13
&  0.5861531690539E-01 & 0.10408E-15 & 0.55098E-11  \\
 14
&  0.5949490764741E-01 & 0.23592E-15 & -.21301E-09  \\
 15
&  0.6020725336886E-01 & 0.13184E-15 & 0.31869E-11  \\
 16
&  0.6077650804037E-01 & 0.18041E-15 & -.21545E-09  \\
 17
&  0.6122088420703E-01 & 0.48572E-16 & 0.14361E-11  \\
 18
&  0.6155390478472E-01 & 0.83267E-16 & -.21675E-09  \\
 19
&  0.6178529976346E-01 & 0.11102E-15 & 0.36146E-12  \\
 20
&  0.6192162112196E-01 & 0.48572E-16 & -.21732E-09  \\
 21
&  0.6196665001384E-01 & 0.00000E+00 & 0.00000E+00  \\
\end{tabular}
\end{center}
\caption{\it
Quadrature weights \eqref{eq_quad_w} with $c = 40$, $n = 41$.
$\lambda_n = \mbox{\text{\rm{i0.69857E-08}}}$.
See Experiment 4.
}
\label{t:test96}
\end{table}
%%%%%%%%%%%
%%%%%%%%%%%%%%%%%%%%%%%%%%%%%%
\paragraph{Experiment 4.}
\label{sec_exp15}
In this experiment, 
% we illustrate Theorem~\ref{thm_tilde_phi}
% and Remark~\ref{rem_w_approx} in Section~\ref{sec_weights}.
we choose, more or less arbitrarily,
band limit $c$ and prolate index $n$.
Then, we compute $t_1, \dots, t_n$ (see \eqref{eq_quad_t}) and
$\psi_n'(t_1), \dots, \psi_n'(t_n)$
via the algorithm of Section~\ref{sec_evaluate_nodes}. 
Also, we evaluate $\psi_n(0), \psi_n'(0)$ via
the algorithm of Section~\ref{sec_evaluate_beta}.
Next, compute approximations $\widetilde{W}_1, \dots, \widetilde{W}_n$
to $W_1,\dots,W_n$
via Algorithm 1 in Section~\ref{sec_evaluate_weights}
(in particular, $\widetilde{W}_j$ is evaluated via
\eqref{eq_wj_as_sum} for every $j=1,\dots,n$).
Also, we compute approximations
$\widehat{W}_1, \dots, \widehat{W}_n$ to $W_1,\dots,W_n$
via
Algorithm 2 in Section~\ref{sec_evaluate_weights}.
All the calculations are carried out in double precision.

We display the results of this experiment
in Table~\ref{t:test96}.
The data in this table correspond to
$c = 40$ and $n = 41$.
Table~\ref{t:test96} has the following structure.
The first column contains the weight index $j$, that varies
between 1 and $(n+1)/2=21$.
The second column contains $\widehat{W}_j$ (
an approximation to $W_j$ evaluated by Algorithm 2
in Section~\ref{sec_evaluate_weights}).
The third column contains the difference between $\widehat{W}_j$
and $\widetilde{W}_j$ (evaluated via \eqref{eq_wj_as_sum} by Algorithm 1).
The last column contains the difference 
\begin{align}
\widehat{W}_j- 
\frac{\widehat{W}_{(n+1)/2} \cdot \left( \psi_n'(0) \right)^2}
      {\left( \psi_n'(t_j) \right)^2 \cdot \left(1 - t_j^2\right) }
\label{eq_exp15_dif}
\end{align}
(see \eqref{eq_w_approx} in Remark~\ref{rem_w_approx}).

\begin{figure} [htbp]
\begin{center}
\includegraphics[width=11.5cm, bb=81   227   529   564, clip=true]
{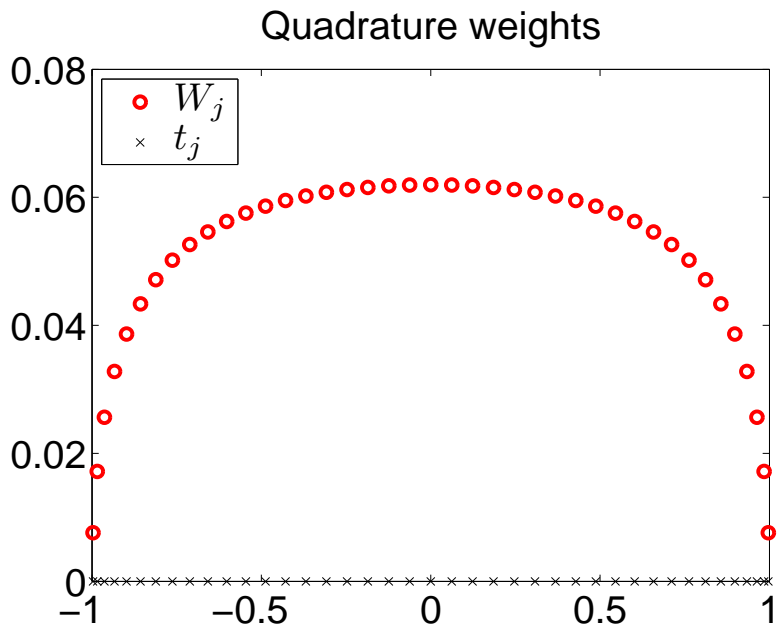}
\caption
{\it
The quadrature weights $W_1,\dots,W_n$ with $c = 40$, $n = 41$.
See Experiment 4.
}
\label{fig:test96}
\end{center}
\end{figure}

In Figure~\ref{fig:test96}, we plot
the weights $W_j$ as a function
of $j=1,\dots,n$.
Each $W_j$ is plotted as a circle above the corresponding node
$t_j$.

We make the following observations from Table~\ref{t:test96}.
First, due to the combination of Theorems~\ref{thm_tilde_phi},
\ref{lem_tilde_phi_ode} in Section~\ref{sec_weights},
the value in the third column would be zero
in exact arithmetic. We observe that, indeed, this value
is zero up to the machine precision, which confirms
Remarks~\ref{rem_algorithm1_acc}, \ref{rem_weights_cost}
in Section~\ref{sec_evaluate_weights}.
(We note that, for $j=1,2,3,4$ and $j=21$, both $\widehat{W}_j$
and $\widetilde{W}_j$ are evaluated via \eqref{eq_wj_as_sum}, and hence
the value in the corresponding rows is exactly zero).
In particular, either of the two approximations $\widetilde{W}_j, \widehat{W}_j$
can be used to evaluate $W_j$ to essentially machine precision.

We also observe that all of the weights $W_1, \dots, W_n$ are positive
(see Theorem~\ref{thm_positive_w} and Remark~\ref{rem_w_always_pos}).
Moreover, $W_j$ grow
monotonically as $j$ increases to $(n+1)/2$.
Finally, we observe that, for all $j=1,\dots,n$,
the value \eqref{eq_exp15_dif} in the last column
is of the order $|\lambda_n|$ (see
Remark~\ref{rem_w_approx}).
%%%%%%%%%%%%%%%%%%%%%
%%%%%%%%%%%%

%%%%%%%%%%%%%%%%%%%%%%%%%%%

\begin{comment}
\begin{figure}
\begin{center}
\begin{tabular}{c}
\includegraphics[width=10cm, bb=68 218 542 574, clip=true ]{test82a.eps} \\
(a) \\
\includegraphics[width=10cm, bb=68 218 542 574, clip=true ]{test82b.eps} \\
(b) \\
\includegraphics[width=8cm, bb=68 218 542 574, clip=true ]{test82c.eps} \\
\includegraphics[width=8cm, bb=68 218 542 574, clip=true ]{test82d.eps} \\
(c) \\
(d) \\
\end{tabular}
\caption{
A
}\label{fig:dist_example}
\end{center}
%\vspace{-0.6cm}
\end{figure}
\end{comment}

\end{document}